\documentclass[11pt]{article}
\usepackage{hyperref, amsfonts, amsmath, amssymb, amsthm, fullpage, enumitem, bm, mathtools, microtype, tikz, subcaption, alphalph, pgf, pgffor}
\usepackage[noabbrev, capitalize]{cleveref}
\usepackage[linesnumbered]{algorithm2e}
\SetAlCapSty{} 
\crefname{algocf}{Algorithm}{Algorithms}
\crefname{appendix}{Appendix}{Appendices}
\Crefname{appendix}{Appendix}{Appendices}

\tikzset{every picture/.append style={scale=.8, thick}}
\definecolor{litegray}{RGB}{192,192,192}
\definecolor{white}{rgb}{1.0, 1.0, 1.0}
\definecolor{g-red}{RGB}{165,14,14}
\definecolor{g-light-red}{RGB}{250,210,207}
\definecolor{g-blue}{RGB}{23,78,166}
\definecolor{g-light-blue}{RGB}{210,227,252}
\usetikzlibrary{decorations.pathreplacing, calc}
\tikzstyle{vertex}=[circle, draw, fill=litegray, inner sep=0pt, minimum width=4pt]
\tikzstyle{red-vertex}=[circle, draw=white, fill=g-red, text=white, font=\footnotesize, inner sep=0pt, minimum width=12pt]
\tikzstyle{reddish-vertex}=[circle, draw=white, fill=g-light-red, text=g-red, font=\footnotesize, inner sep=0pt, minimum width=12pt]
\tikzstyle{a-vertex}=[circle, draw=white, double=g-red, fill=white, text=g-red, font=\footnotesize, inner sep=0pt, minimum width=11pt]
\tikzstyle{blue-vertex}=[rectangle, rounded corners=2pt, draw=white, fill=g-blue, text=white, font=\footnotesize, inner sep=0pt, minimum width=12pt, minimum height=12pt]
\tikzstyle{bluish-vertex}=[rectangle, rounded corners=2pt, draw=white, fill=g-light-blue, text=g-blue, font=\footnotesize, inner sep=0pt,
minimum width=12pt, minimum height=12pt]
\tikzstyle{b-vertex}=[rectangle, rounded corners=2pt, draw=white, double=g-blue, fill=white, text=g-blue, font=\footnotesize, inner sep=0pt, minimum width=11pt, minimum height=11pt]
\tikzstyle{red-edge}=[draw=g-red, very thick]
\tikzstyle{blue-edge}=[draw=g-blue, very thick]

\pgfdeclareshape{degree-one-circle}{
  \inheritsavedanchors[from=circle]
  \inheritanchorborder[from=circle]
  \inheritanchor[from=circle]{text}
  \inheritanchor[from=circle]{base}
  \inheritanchor[from=circle]{mid}
  \inheritanchor[from=circle]{center}
  \inheritanchor[from=circle]{north}
  \inheritanchor[from=circle]{south}
  \inheritanchor[from=circle]{east}
  \inheritanchor[from=circle]{west}
  \inheritanchor[from=circle]{north east}
  \inheritanchor[from=circle]{north west}
  \inheritanchor[from=circle]{south east}
  \inheritanchor[from=circle]{south west}
  \backgroundpath{
    \pgfpathmoveto{\pgfpointadd{\centerpoint}{\pgfpoint{0pt}{6pt}}}
    \pgfpathcurveto
      {\pgfpointadd{\centerpoint}{\pgfpoint{3.314pt}{6pt}}}
      {\pgfpointadd{\centerpoint}{\pgfpoint{6pt}{3.314pt}}}
      {\pgfpointadd{\centerpoint}{\pgfpoint{6pt}{0pt}}}
    \pgfpathcurveto
      {\pgfpointadd{\centerpoint}{\pgfpoint{6pt}{-3.314pt}}}
      {\pgfpointadd{\centerpoint}{\pgfpoint{3.314pt}{-6pt}}}
      {\pgfpointadd{\centerpoint}{\pgfpoint{0pt}{-6pt}}}
    \pgfpathlineto{\pgfpointadd{\centerpoint}{\pgfpoint{-4pt}{-6pt}}}
    \pgfpathcurveto
      {\pgfpointadd{\centerpoint}{\pgfpoint{-5.105pt}{-6pt}}}
      {\pgfpointadd{\centerpoint}{\pgfpoint{-6pt}{-5.105pt}}}
      {\pgfpointadd{\centerpoint}{\pgfpoint{-6pt}{-4pt}}}
    \pgfpathlineto{\pgfpointadd{\centerpoint}{\pgfpoint{-6pt}{0pt}}}
    \pgfpathcurveto
      {\pgfpointadd{\centerpoint}{\pgfpoint{-6pt}{3.314pt}}}
      {\pgfpointadd{\centerpoint}{\pgfpoint{-3.314pt}{6pt}}}
      {\pgfpointadd{\centerpoint}{\pgfpoint{0pt}{6pt}}}
    \pgfpathclose
  }
}

\pgfdeclareshape{degree-two-circle}{
  \inheritsavedanchors[from=circle]
  \inheritanchorborder[from=circle]
  \inheritanchor[from=circle]{text}
  \inheritanchor[from=circle]{base}
  \inheritanchor[from=circle]{mid}
  \inheritanchor[from=circle]{center}
  \inheritanchor[from=circle]{north}
  \inheritanchor[from=circle]{south}
  \inheritanchor[from=circle]{east}
  \inheritanchor[from=circle]{west}
  \inheritanchor[from=circle]{north east}
  \inheritanchor[from=circle]{north west}
  \inheritanchor[from=circle]{south east}
  \inheritanchor[from=circle]{south west}
  \backgroundpath{
    \pgfpathmoveto{\pgfpointadd{\centerpoint}{\pgfpoint{0pt}{6pt}}}
    \pgfpathcurveto
      {\pgfpointadd{\centerpoint}{\pgfpoint{3.314pt}{6pt}}}
      {\pgfpointadd{\centerpoint}{\pgfpoint{6pt}{3.314pt}}}
      {\pgfpointadd{\centerpoint}{\pgfpoint{6pt}{0pt}}}
    \pgfpathlineto{\pgfpointadd{\centerpoint}{\pgfpoint{6pt}{-4pt}}}
    \pgfpathcurveto
      {\pgfpointadd{\centerpoint}{\pgfpoint{6pt}{-5.105pt}}}
      {\pgfpointadd{\centerpoint}{\pgfpoint{5.105pt}{-6pt}}}
      {\pgfpointadd{\centerpoint}{\pgfpoint{4pt}{-6pt}}}
    \pgfpathlineto{\pgfpointadd{\centerpoint}{\pgfpoint{-4pt}{-6pt}}}
    \pgfpathcurveto
      {\pgfpointadd{\centerpoint}{\pgfpoint{-5.105pt}{-6pt}}}
      {\pgfpointadd{\centerpoint}{\pgfpoint{-6pt}{-5.105pt}}}
      {\pgfpointadd{\centerpoint}{\pgfpoint{-6pt}{-4pt}}}
    \pgfpathlineto{\pgfpointadd{\centerpoint}{\pgfpoint{-6pt}{0pt}}}
    \pgfpathcurveto
      {\pgfpointadd{\centerpoint}{\pgfpoint{-6pt}{3.314pt}}}
      {\pgfpointadd{\centerpoint}{\pgfpoint{-3.314pt}{6pt}}}
      {\pgfpointadd{\centerpoint}{\pgfpoint{0pt}{6pt}}}
    \pgfpathclose
  }
}

\tikzstyle{red-degree-one-vertex}=[red-vertex, shape=degree-one-circle]
\tikzstyle{red-degree-two-vertex}=[red-vertex, shape=degree-two-circle]
\tikzstyle{reddish-degree-one-vertex}=[reddish-vertex, shape=degree-one-circle]
\tikzstyle{reddish-degree-two-vertex}=[reddish-vertex, shape=degree-two-circle]

\newcommand{\defc}[1]{%
  \foreach \v/\x/\y in {#1} {\coordinate (\v) at (\x,\y);}
}
\newcommand{\drawv}[4]{%
  \foreach \v in {#1} {\node[red-vertex] at (\v) {$\v$};}
  \foreach \v in {#2} {\node[reddish-vertex] at (\v) {$\v$};}
  \foreach \v in {#3} {\node[blue-vertex] at (\v) {$\v$};}
  \foreach \v in {#4} {\node[bluish-vertex] at (\v) {$\v$};}
}
\newcommand{\drawvnl}[4]{%
  \foreach \v in {#1} {\node[red-vertex] at (\v) {};}
  \foreach \v in {#2} {\node[reddish-vertex] at (\v) {};}
  \foreach \v in {#3} {\node[blue-vertex] at (\v) {};}
  \foreach \v in {#4} {\node[bluish-vertex] at (\v) {};}
}
\newcommand{\drawe}[3]{%
  \foreach \u/\v in {#1} {\draw[red-edge] (\u) -- (\v);}
  \foreach \u/\v in {#2} {\draw[blue-edge] (\u) -- (\v);}
  \foreach \u/\v in {#3} {\draw (\u) -- (\v);}
}
\newcommand{\arcthreeone}[1]{\draw (#1) arc [start angle=172.380146, end angle=44.489752, radius=1.760000];}
\newcommand{\arcthreeminusone}[1]{\draw (#1) arc [start angle=135.510248, end angle=7.619854, radius=1.760000];}
\newcommand{\arcthreeminusoneneg}[1]{\draw (#1) arc [start angle=-172.380146, end angle=-44.489752, radius=1.760000];}
\newcommand{\arcthreeoneneg}[1]{\draw (#1) arc [start angle=-135.510248, end angle=-7.619854, radius=1.760000];}
\newcommand{\arcfourone}[1]{\draw (#1) arc [start angle=165.800136, end angle=42.272351, radius=2.340000];}
\newcommand{\arcfourminusone}[1]{\draw (#1) arc [start angle=137.727649, end angle=14.199864, radius=2.340000];}
\newcommand{\arcfiveminusone}[1]{\draw (#1) arc [start angle=139.513295, end angle=17.866840, radius=2.920000];}
\newcommand{\arcfiveone}[1]{\draw (#1) arc [start angle=162.133160, end angle=40.486705, radius=2.920000];}

\newlist{kase}{enumerate}{7} 
\setlistdepth{7}
\setlist[kase,1]{label=(\arabic*)}
\setlist[kase,2]{label=(\arabic{kasei}.\arabic*)}
\setlist[kase,3]{label=(\arabic{kasei}.\arabic{kaseii}.\arabic*)}
\setlist[kase,4]{label=(\arabic{kasei}.\arabic{kaseii}.\arabic{kaseiii}.\arabic*)}
\setlist[kase,5]{label=(\arabic{kasei}.\arabic{kaseii}.\arabic{kaseiii}.\arabic{kaseiv}.\arabic*)}
\setlist[kase,6]{label=(\arabic{kasei}.\arabic{kaseii}.\arabic{kaseiii}.\arabic{kaseiv}.\arabic{kasev}.\arabic*)}
\setlist[kase,7]{label=(\arabic{kasei}.\arabic{kaseii}.\arabic{kaseiii}.\arabic{kaseiv}.\arabic{kasev}.\arabic{kasevi}.\arabic*)}

\newtheorem{theorem}{Theorem}[section]

\newtheorem{lemma}[theorem]{Lemma}
\newtheorem{proposition}[theorem]{Proposition}
\newtheorem{problem}[theorem]{Problem}

\theoremstyle{definition}
\newtheorem{definition}[theorem]{Definition}

\AddToHook{env/corollary/begin}{\crefalias{theorem}{corollary}}
\AddToHook{env/lemma/begin}{\crefalias{theorem}{lemma}}
\AddToHook{env/proposition/begin}{\crefalias{theorem}{proposition}}
\AddToHook{env/problem/begin}{\crefalias{theorem}{problem}}
\AddToHook{env/definition/begin}{\crefalias{theorem}{definition}}

\theoremstyle{remark}
\newtheorem*{remark}{Remark}

\newtheorem*{claim*}{Claim}
\newtheorem{claim}{Claim}

\crefname{claim}{Claim}{Claims}

\newenvironment{claimproof}[1][Proof]{\begin{proof}[#1]}{\end{proof}}
\newenvironment{ctikz}{\begin{center}\begin{tikzpicture}}{\end{tikzpicture}\end{center}}
\newenvironment{ctikzq}{\[\begin{tikzpicture}}{\end{tikzpicture}\qedhere\]}
\newenvironment{ctikzpq}{\[\pushQED{\qed}\begin{tikzpicture}}{\end{tikzpicture}\qedhere\popQED\]}
\NewDocumentEnvironment{csubfig}{m m}{\begin{subfigure}{#1}\centering\begin{tikzpicture}}{\end{tikzpicture}\caption{} \label{fig:#2}\end{subfigure}}

\newread\tabledataread
\newcount\tabledatacol
\newcommand{\tabledataline}{}
\newcommand{\tablenextline}{}
\newcommand{\inputtablecontents}{}
\newcommand{\outputtablecontents}{}
\newcommand{\appendtabledata}[2]{%
  \begingroup
  \toks0=\expandafter{#1}%
  \edef\x{\endgroup\global\def\noexpand#1{\the\toks0 #2}}%
  \x
}
\newcommand{\buildtabledata}[3]{%
  \global\let#1\empty
  \tabledatacol=0
  \openin\tabledataread=#2\relax
  \read\tabledataread to \tabledataline
  \ifeof\tabledataread
  \else
    \buildtabledataread{#1}{#3}%
  \fi
  \closein\tabledataread
}
\newcommand{\buildtabledataread}[2]{%
  \advance\tabledatacol by 1
  \appendtabledata{#1}{\tabledataline}%
  \read\tabledataread to \tablenextline
  \ifeof\tabledataread
  \else
    \ifnum\tabledatacol=#2\relax
      \tabledatacol=0
      \appendtabledata{#1}{\noexpand\tabularnewline}%
    \else
      \appendtabledata{#1}{&}%
    \fi
    \let\tabledataline\tablenextline
    \buildtabledataread{#1}{#2}%
  \fi
}

\DeclarePairedDelimiter\abs{\lvert}{\rvert}
\DeclarePairedDelimiter\floor{\lfloor}{\rfloor}
\DeclarePairedDelimiter\ceil{\lceil}{\rceil}

\newcommand{\dset}[2]{\left\{{#1}\colon{#2}\right\}}
\newcommand{\sset}[1]{\left\{{#1}\right\}}

\newcommand{\la}{\lambda}
\newcommand{\eps}{\varepsilon}

\newcommand{\ce}[1]{\cref{fig:cut-enhancer-#1}}

\makeatletter
\renewcommand{\nsim}{\mathrel{\mathpalette\n@sim\relax}}
\newcommand{\n@sim}[2]{%
 \ooalign{%
  $\m@th#1\sim$\cr
  \hidewidth$\m@th#1\rotatebox[origin=c]{50}{$#1-$}$\hidewidth\cr
 }%
}
\makeatother

\pgfmathtruncatemacro{\postail}{16}
\newcommand{\pathfourneg}{10}
\newcommand{\pathfoursubgraphneg}{10}
\pgfmathtruncatemacro{\paththreeneg}{\pathfoursubgraphneg+2}
\pgfmathtruncatemacro{\paththreesubgraphneg}{max(\pathfourneg, \paththreeneg)}
\pgfmathtruncatemacro{\pathtwoneg}{\paththreesubgraphneg+2}
\pgfmathtruncatemacro{\negtail}{max(\paththreeneg,\pathtwoneg)}
\pgfmathtruncatemacro{\abstail}{max(\postail,\negtail)}

\begin{document}

\title{Median eigenvalues of subcubic graphs}
\author{
  Hricha Acharya\thanks{School of Mathematical and Statistical Sciences, Arizona State University, Tempe, AZ 85281, USA. Email: \texttt{\char`\{hachary3, bjeter1, zilinj\char`\}@asu.edu}.}
  \and Benjamin Jeter\footnotemark[1]
  \and Zilin Jiang\footnotemark[1]${\ }^{,}$\thanks{School of Computing and Augmented Intelligence, Arizona State University, Tempe, AZ 85281, USA. Supported in part by the Simons Foundation through its Travel Support for Mathematicians program and by U.S.\ taxpayers through NSF grant 2451581.}
}
\date{}

\maketitle

\begin{abstract}
  We show that the median eigenvalues of every connected graph of maximum degree at most three, except for the Heawood graph, are at most $1$ in absolute value, resolving open problems posed by Fowler and Pisanski, and by Mohar.
\end{abstract}

\section{Introduction}

A graph is \emph{subcubic} if its maximum degree is at most $3$. In mathematical chemistry, connected subcubic graphs are often referred to as \emph{chemical graphs} because they naturally model the skeletons of organic molecules, where atoms such as carbon typically form bonds with up to three neighboring atoms.

The study of connected subcubic graphs has significant applications in theoretical chemistry, particularly in the analysis of molecular orbital models. For instance, in the H\"uckel molecular orbital model, the eigenvalues of a molecule's skeleton graph correspond to its $\pi$-electron energy levels. Among these, the \emph{median eigenvalues} --- associated with the highest occupied molecular orbital (HOMO) and the lowest unoccupied molecular orbital (LUMO) --- are of particular interest (see \cite{LLSG13,WYY18} and the references therein).

To formalize, let $G$ be an $n$-vertex graph, and let $\la_1(G) \ge \la_2(G) \ge \dots \ge \la_n(G)$ denote the eigenvalues of its adjacency matrix. The \emph{median eigenvalues} are defined as $\la_h$ and $\la_l$, where $h = \floor{(n+1)/2}$ and $l = \ceil{(n+1)/2}$.

Computational studies of Fowler and Pisanski \cite{FP10-b,FP10-a} suggest that the median eigenvalues of most connected subcubic graphs lie in the interval $[-1, 1]$. However, there is a single notable exception: the Heawood graph, whose median eigenvalues are $\pm\sqrt2$ (see \cref{fig:heawood}). The Heawood graph is the incidence graph of the Fano plane, that is, the projective plane of order $2$.

\begin{figure}[t]
  \centering
  \begin{tikzpicture}[scale=2]
    \foreach \i in {0,...,13} {
      \node[vertex] (\i) at ({360/14 * \i}:1) {};
    }
    \foreach \i in {0,...,13} {
      \pgfmathtruncatemacro{\next}{Mod(\i + 1, 14)}
      \draw (\i) -- (\next);
    }
    \foreach \i in {0,...,13} {
      \pgfmathtruncatemacro{\iseven}{Mod(\i, 2)}
      \ifnum\iseven=0
        \pgfmathtruncatemacro{\target}{Mod(\i + 5, 14)}
      \else
        \pgfmathtruncatemacro{\target}{Mod(\i - 5, 14)}
      \fi
      \draw (\i) -- (\target);
    }
  \end{tikzpicture}
  \caption{The Heawood graph.} \label{fig:heawood}
\end{figure}

Fowler and Pisanski \cite{FP10-b} conjectured that the median eigenvalues of all but finitely many connected subcubic graphs lie in $[-1,1]$, and they confirmed their conjecture for all connected subcubic trees. The interval $[-1,1]$ would be optimal because of the construction of an infinite family of bipartite connected subcubic graphs whose median eigenvalues are $\pm 1$ by Guo and Mohar \cite{GM14}.

Subsequently, Mohar~\cite{M15} conjectured that the median eigenvalues of all \emph{planar} subcubic graphs lie in $[-1,1]$, and he confirmed in \cite{M13} his conjecture for all \emph{bipartite planar} subcubic graphs. We point out that the Heawood graph is bipartite but not planar.

Later, Mohar~\cite{M16} greatly extended his result to all \emph{bipartite} connected subcubic graphs, except for the Heawood graph. In addition, he refined his proof to show that a positive fraction of the eigenvalues around the median eigenvalues belong to $[-1,1]$ (see \cite[Theorem~1.4]{M16}). Recently, Wang and Zhang \cite{WZ24} supplemented Mohar's main result in \cite{M16} by treating connected subcubic graphs that contain $K_{2,3}$ as a subgraph, as well as those that contain no $K_4$ as a minor. The latter generalizes an earlier result of Benediktovich~\cite{B14} for outerplanar connected subcubic graphs. Concurrent with our work, Chen, Wang and Zhang \cite{CWZ25} proved that the median eigenvalues lie in $[-1,1]$ for connected subcubic graphs in which every subgraph has average degree less than $44/17$, and also for planar connected subcubic graphs of girth at least $8$.

We completely settle the aforementioned conjectures in a strong form.

\begin{theorem} \label{thm:main1}
  The median eigenvalues of every connected subcubic graph, except for the Heawood graph, lie in the interval $[-1,1]$.
\end{theorem}

Just like Mohar's result \cite[Theorem~1.4]{M16}, our proof can also be extended to a positive fraction of the eigenvalues around the medians for large connected subcubic graphs.

\begin{theorem} \label{thm:main2}
  There exists a constant $\eta > 0$ such that for every connected subcubic graph $G$ on $n$ vertices, the $i$-th largest eigenvalue $\la_i(G)$ of $G$ lies in $[-1,1]$ for every integer $i \in [\tfrac12(1-\eta)(n+1), \tfrac12(1+\eta)(n+1)]$.
\end{theorem}

\begin{remark}
  When $G$ is the Heawood graph, or more generally, when $G$ is a small graph on an even number of vertices, \cref{thm:main2} holds for the trivial reason that $[\tfrac12 (1-\eta)(n+1), \tfrac12 (1+\eta)(n+1)]$ contains no integers.
\end{remark}

\section{Proof ideas}

The use of the Cauchy interlacing theorem in the context of median eigenvalues of connected subcubic graphs first appeared in \cite{M13}. We use the following immediate consequence of the Cauchy interlacing theorem to control graph eigenvalues.

\begin{definition}[Tails and their estimates]
  Given a graph $G$, its \emph{positive tail}, denoted by $t^+_G$, is the number of eigenvalues greater than $1$, and its \emph{negative tail}, denoted by $t^-_G$, is the number of eigenvalues less than $-1$. Furthermore, given a vertex subset $A$ of $G$, the \emph{positive tail estimate} of $G$ with respect to $A$ is defined by
  \[
    t^+_G(A) := \abs{A} + t^+_{G-A},
  \]
  and the \emph{negative tail estimate} of $G$ with respect to $A$ is defined by
  \[
    t^-_G(A) := \abs{A} + t^-_{G-A}.
  \]
\end{definition}

\begin{lemma} \label{lem:key}
  For every vertex subset $A$ of a graph $G$, the positive tail $t^+_G$ of $G$ is at most the positive tail estimate $t^+_G(A)$, and the negative tail $t^-_G$ of $G$ is at most the negative tail estimate $t^-_G(A)$.
\end{lemma}

\begin{proof}
  Set $k = t_{G-A}^+$. Thus $\la_{k+1}(G-A) \le 1$. By the Cauchy interlacing theorem, $\la_{\abs{A} + k + 1}(G) \le 1$, and so $t_G^+ \le \abs{A} + k = t_G^+(A)$. Similarly, one can also show that $t^-_G \le t^-_G(A)$.
\end{proof}

The maximum cut naturally provides a good starting point for tail estimates.

\begin{proposition} \label{lem:max-cut}
  For every maximum cut $(A, B)$ of a subcubic graph $G$, the edges of $G[A] \cup G[B]$ form a matching in $G$, and, in particular, $t_G^\pm(A) = \abs{A}$ and $t_G^\pm(B) = \abs{B}$.
\end{proposition}

\begin{proof}
  Suppose that some vertex $v \in A$ has at least two neighbors in $A$. Since $G$ is subcubic, $v$ has at most one neighbor in $B$. Moving $v$ from $A$ to $B$ therefore increases the size of the cut-set, contradicting the maximality of $(A,B)$. The same argument applies to vertices in $B$, and so every vertex has degree at most $1$ in $G[A] \cup G[B]$. Hence the edges of $G[A] \cup G[B]$ form a matching in $G$.

  In particular, $G[A]$ and $G[B]$ are disjoint unions of isolated vertices and edges, so all their eigenvalues lie in $[-1,1]$. Therefore $t^\pm_{G-A} = t^\pm_{G[B]} = 0$ and $t^\pm_{G-B} = t^\pm_{G[A]} = 0$, which gives $t_G^\pm(A) = \abs{A}$ and $t_G^\pm(B) = \abs{B}$.
\end{proof}

If $\abs{A} \neq \abs{B}$ for some maximum cut $(A, B)$ of $G$, then \cref{lem:key,lem:max-cut} imply that $t_G^\pm \le \min(\abs{A}, \abs{B}) < \abs{G}/2$, and so the median eigenvalues of $G$ lie in $[-1,1]$. It remains to handle graphs whose maximum cuts are all balanced, that is, $\abs{A} = \abs{B}$ for every maximum cut $(A, B)$ of $G$. The remaining case requires more work.

The strategy is to adjust the maximum cut to improve the tail estimates of $G$. This follows the philosophy of Mohar's imbalance argument \cite{M16}. We package the adjustments that decrease tail estimates in the following definitions.

\begin{definition}[Tail reducers] \label{def:tail-reducer}
  Given a cut $(A, B)$ of a graph $G$, a vertex subset $D$ of $G$ is a \emph{positive tail reducer} (respectively, \emph{negative tail reducer}) of $G$ with respect to $(A, B)$ if either $D \subseteq A$ or $D \subseteq B$, and $t_Q^+ < \abs{D}$ (respectively, $t_Q^- < \abs{D}$), where $Q$ is the union of the connected components of $G[D']$ that contain vertices in $D$, and
  \[
    D' = \begin{cases}
      D \cup B & \text{if }D \subseteq A;\\
      D \cup A & \text{if }D \subseteq B.\\
    \end{cases}
  \]
  We say that $D$ is an \emph{absolute tail reducer} if it is both a positive tail reducer and a negative tail reducer.
\end{definition}

As the name suggests, a tail reducer decreases the corresponding tail estimate.

\begin{lemma}[Lemma~2.3 of Mohar~\cite{M16}] \label{lem:q}
  For every cut $(A, B)$ of a graph $G$, and every vertex subset $D$ of $G$, if $D$ is a positive tail reducer with respect to $(A, B)$, then $$t^+_G(A \setminus D) < t^+_G(A) \text{ or } t^+_G(B \setminus D) < t^+_G(B).$$
  Likewise, if $D$ is a negative tail reducer with respect to $(A, B)$, then $$t^-_G(A \setminus D) < t^-_G(A) \text{ or } t^-_G(B \setminus D) < t^-_G(B).$$
\end{lemma}

Although \cite[Lemma~2.3]{M16} only concerns positive tail reducers, the proof for negative tail reducers is almost verbatim. Since the terminology in \cite{M16} is different from ours, we include the proof here for the reader's convenience.

\begin{proof}
  We only prove the case where $D \subseteq A$ is a positive tail reducer with respect to $(A, B)$. Let $Q$ be the union of the connected components of $G[D \cup B]$ that contain vertices in $D$. By \cref{def:tail-reducer}, we have $t^+_Q < \abs{D}$, and therefore
  \[
    t_G^+(A \setminus D) = \abs{A \setminus D} + t^+_{G - (A \setminus D)} = \abs{A} - \abs{D} + t^+_{G[D \cup B]} \le \abs{A} - \abs{D} + t^+_{Q} + t^+_{G[B]} < \abs{A} + t^+_{G[B]} = t^+_G(A).
  \]
  The other cases can be proved similarly.
\end{proof}

To decrease the positive tail estimate of a \emph{bipartite} connected subcubic graph, Mohar proved the following result. We denote by $G(x, r)$ the subgraph of $G$ on the vertices within distance $r$ from $x$, where $x$ could be a vertex, an edge, or a vertex subset.

\begin{lemma}[Lemma~3.1 of Mohar~\cite{M16}] \label{lem:mohar}
  For every bipartite connected subcubic graph $G$ other than the Heawood graph, and every vertex $v$ of $G$, if $(A, B)$ is the bipartition of $G$, then there exists a vertex subset $C$ of $G(v, 17)$ such that $C$ is an absolute tail reducer of $G$ with respect to $(A,B)$. \qed
\end{lemma}

\begin{remark}
  In \cref{lem:mohar}, since $G$ is bipartite, the spectrum of an induced subgraph of $G$ is symmetric with respect to $0$, and so a vertex subset $C$ is a positive tail reducer if and only if $C$ is a negative tail reducer.
\end{remark}

The non-bipartite setting introduces two difficulties that do not occur in \cite{M16}. First, a maximum cut need not be unique or canonical, so the proof must exploit only local consequences of maximality. Second, without bipartite spectral symmetry, a set that reduces the positive tail estimate need not reduce the negative tail estimate. This is why our main local work comes in two parallel forms, one for $t_G^+$ and one for $t_G^-$.

\begin{lemma} \label{lem:pos-tail-16}
  For every maximum cut $(A, B)$ of a subcubic graph $G$, and for every edge $\alpha$ of $G[A] \cup G[B]$, there exists a vertex subset of $G(\alpha,\postail)$ that is a positive tail reducer of $G$ with respect to $(A, B)$.
\end{lemma}

\begin{lemma} \label{lem:neg-tail-14}
  For every maximum cut $(A, B)$ of a subcubic graph $G$, and for every edge $\alpha$ of $G[A] \cup G[B]$, there exists a vertex subset of $G(\alpha,\negtail)$ that is a negative tail reducer of $G$ with respect to $(A, B)$.
\end{lemma}

We are ready to prove our first main theorem.

\begin{proof}[Proof of \cref{thm:main1} assuming \cref{lem:pos-tail-16,lem:neg-tail-14}]
  Let $G$ be an $n$-vertex connected subcubic graph that is not the Heawood graph, and let $(A, B)$ be a maximum cut of $G$; if $G$ is bipartite, choose $(A,B)$ to be a bipartition. By \cref{lem:max-cut}, we have
  \[
    t^\pm_G(A) \le \abs{A} \text{ and } t^\pm_G(B) \le \abs{B}.
  \]
  If $\abs{A} \neq \abs{B}$, \cref{lem:key} implies that \[
    t^\pm_G \le \min(t^\pm_G(A), t^\pm_G(B)) \le \min(\abs{A}, \abs{B}) < \abs{G}/2.
  \]

  From now on, assume that $\abs{A} = \abs{B}$. If $G$ is bipartite, then \cref{lem:mohar} provides an absolute tail reducer with respect to its bipartition. If $G$ is not bipartite, then $G[A] \cup G[B]$ contains an edge, and \cref{lem:pos-tail-16,lem:neg-tail-14} provide a positive tail reducer and a negative tail reducer with respect to $(A,B)$. In either case, \cref{lem:q} gives a set whose removal from $A$ or from $B$ strictly decreases the corresponding tail estimate. Applying \cref{lem:key}, we obtain
  \[
    t^+_G < \abs{G}/2 \quad \text{and} \quad t^-_G < \abs{G}/2,
  \]
  and hence the median eigenvalues of $G$ lie in $[-1,1]$.
\end{proof}

It remains to prove the two tail-reducer lemmas by analyzing, case by case, the neighborhoods of an edge of $G[A] \cup G[B]$. We begin with a maximum cut $(A,B)$ and repeatedly use local consequences of its maximality. Throughout the remainder of the paper, $\oplus$ denotes symmetric difference.

\begin{definition}[Cut enhancer]
  Given a cut $(A,B)$ of a graph $G$, a vertex subset $C$ is a \emph{cut enhancer} if the cut-set of $(A \oplus C, B \oplus C)$ is larger than that of $(A,B)$.
\end{definition}

A maximum cut therefore admits no cut enhancer. This observation yields many of the local restrictions used below.

To streamline the case analysis, we often modify the initial maximum cut $(A,B)$ so that a given case reduces to one already considered. We formalize this operation as follows.

\begin{definition}[Cut preserver and flipping]
  Given a cut $(A, B)$ of a graph $G$, a vertex subset $C$ of $G$ is a \emph{cut preserver} if the cut-set of $(A \oplus C, B \oplus C)$ has the same size as that of $(A, B)$. In this case, passing from $(A, B)$ to $(A \oplus C, B \oplus C)$ is called \emph{flipping} the cut preserver $C$.
\end{definition}

We shall only make use of the following cut preserver.

\begin{proposition} \label{prop:cut-preserver}
  For every cut $(A, B)$ of a graph $G$, and every edge $ab$ of $G$ with $a \in A$ and $b \in B$, if both $a$ and $b$ have degree $3$, $a$ has exactly one neighbor in $A$, and $b$ has exactly one neighbor in $B$, then $\sset{a,b}$ is a cut preserver. \qed
\end{proposition}

In our application, $(A, B)$ is a maximum cut of a subcubic graph $G$, and the cut preserver $C$ satisfies $\abs{A \cap C} = \abs{B \cap C}$ (see \cref{prop:cut-preserver}). Under these conditions, the next tool produces a tail reducer with respect to $(A, B)$ from one with respect to $(A \oplus C, B \oplus C)$.

\begin{lemma} \label{lem:mod-cut-preserver}
  For every maximum cut $(A, B)$ of a subcubic graph $G$, and every pair of vertex subsets $C$ and $D$ of $G$, if $C$ is a cut preserver with $\abs{A \cap C} = \abs{B \cap C}$, and $D$ is a positive (respectively, negative) tail reducer with respect to $(A \oplus C, B \oplus C)$, then $(C \cup D) \cap A$ or $(C \cup D) \cap B$ is a positive (respectively, negative) tail reducer with respect to $(A, B)$.
\end{lemma}

\begin{figure}
  \centering
  \begin{tikzpicture}
    \fill[white] (-2.5,0) arc[start angle=180, end angle=0, radius=2.5] -- cycle;
    \fill[litegray] (2.5,0) arc[start angle=0, end angle=-180, radius=2.5] -- cycle;
    \fill[litegray] (1.5,0) arc[start angle=0, end angle=180, radius=1] -- cycle;
    \fill[white] (1.5,0) arc[start angle=0, end angle=-180, radius=1] -- cycle;

    \draw (0,0) circle (2.5);
    \draw (-2.5,0) -- (2.5,0);
    \draw (-1.5,0) arc[start angle=180, end angle=0, radius=1];
    \draw (0.5,0) circle (1);

    \node[font=\footnotesize] at (0, -2) {$B \oplus C$};
    \node[font=\footnotesize] at (0, 2) {$A \oplus C$};
    \node[font=\footnotesize] at (-.8, .5) {$D$};
    \node[font=\footnotesize] at (0.8, .5) {$C$};
  \end{tikzpicture}
  \caption{A Venn diagram for \cref{lem:mod-cut-preserver}, where white represents $A$ and gray represents $B$.}
  \label{fig:venn}
\end{figure}

\begin{proof}
  Without loss of generality, we may assume that $D \subseteq A \oplus C$, and $D$ is a positive tail reducer with respect to $(A \oplus C, B \oplus C)$. Let $Q$ be the union of the connected components of $G[D \cup (B \oplus C)]$ that contain vertices in $D$. Since $(A \oplus C, B \oplus C)$ is also a maximum cut, \cref{lem:max-cut} implies that the maximum degree of $G[B \oplus C]$ is at most $1$. Thus
  \[
    t^+_{G[D \cup (B \oplus C)]} = t^+_Q < \abs{D}.
  \]
  The Cauchy interlacing theorem implies that
  \[
    t^+_{G[D \cup (B \oplus C) \cup (B \cap C)]} \le t^+_{G[D \cup (B \oplus C)]} + \abs{(B \cap C) \setminus D} < \abs{D} + \abs{(B \cap C) \setminus D}.
  \]
  Set $D' = (C \cup D) \cap A$. From \cref{fig:venn} and the assumption $\abs{A \cap C} = \abs{B \cap C}$, we deduce that
  \begin{gather*}
    D \cup (B\oplus C) \cup (C \setminus (A \cup D)) = D' \cup B, \\
    \abs{D} + \abs{(B \cap C) \setminus D} = \abs{D \setminus C} + \abs{B \cap C} = \abs{D \setminus C} + \abs{A \cap C} = \abs{D'}.
  \end{gather*}
  Let $Q'$ be the union of the connected components of $G[D' \cup B]$ that contain vertices in $D'$. \cref{lem:max-cut} implies that the maximum degree of $G[B]$ is at most $1$, and so $t^+_{Q'} = t^+_{G[D' \cup B]} < \abs{D'}$. Thus $D'$ is a positive tail reducer with respect to $(A, B)$.
\end{proof}

\begin{remark} \label{rmk:mod-cut-preserver}
  The proof above actually shows that, in \cref{lem:mod-cut-preserver}, if in addition $D \subseteq A \oplus C$, then $(C \cup D) \cap A$ is a positive (respectively, negative) tail reducer with respect to $(A, B)$.
\end{remark}

To facilitate case analysis, we introduce the following organizational tool that tracks the interaction between the edges outside the cut-set.

\begin{definition}[Underlying multigraph]
  Given a cut $(A, B)$ of a graph $G$, the \emph{underlying multigraph} $M$ of $G$ with respect to $(A, B)$ is the multigraph on the edges of $G[A] \cup G[B]$, where $k$ edges are placed between $\alpha$ and $\beta$ in $M$ if the bipartite graph $G[\alpha, \beta]$ contains $k$ edges.
\end{definition}

Recall from \cref{lem:max-cut} that, for every maximum cut $(A,B)$ of a subcubic graph $G$, the edges of $G[A] \cup G[B]$ form a matching. Hence every edge of the underlying multigraph joins an edge of $G[A]$ to an edge of $G[B]$, so the underlying multigraph is bipartite.

\begin{proposition} \label{lem:multigraph-bipartite}
  For every maximum cut $(A, B)$ of a subcubic graph $G$, the underlying multigraph of $G$ with respect to $(A,B)$ is bipartite. \qed
\end{proposition}

The proofs of \cref{lem:pos-tail-16,lem:neg-tail-14} follow the same outline. Because $M$ is bipartite, any relevant neighborhood that is simple, has maximum degree at most $2$, and is not an isolated vertex must be a path or an even cycle. The path lemmas \cref{lem:path-5-8,lem:path-5-neg-8} handle long paths and long even cycles uniformly; separate lemmas in \cref{sec:pos,sec:neg} cover the remaining short cases. These arguments either directly exhibit one of the finitely many positive or negative tail reducers displayed in \cref{fig:ptr,fig:ntr}, or first flip cut preservers and then transfer the resulting tail reducer back to the original maximum cut using \cref{lem:mod-cut-preserver}.

The arguments in \cref{sec:pos,sec:neg} reduce the proof to a finite analysis of local configurations, with each case subject to precise distinctness, inducedness, and degree conditions. This makes this part of the argument particularly amenable to formal verification. As an independent check of the case analysis and its use in \cref{lem:pos-tail-16,lem:neg-tail-14}, we developed a Lean 4 formalization using Codex 0.146 with OpenAI's GPT-5.6 Sol model at high reasoning effort. The formalization verifies the relevant local statements and the completeness of the case analysis; it does not formalize the entire paper.

Developing the formalization also served as a useful audit of the mathematical argument: it identified several minor omissions in an earlier version of the manuscript. All of these issues have been corrected in the present version. We reviewed the formal statements and confirmed that the ancillary Lean project builds successfully.

The rest of the paper is organized as follows. In \cref{sec:tools}, we introduce common notation, the concrete cut enhancers used later, basic local degree constraints, and the elementary facts about the underlying multigraph. In \cref{sec:pos,sec:neg}, we prove the positive and negative tail-reducer lemmas, respectively. In \cref{sec:pos-frac}, we combine local tail reducers to prove the positive-fraction result, \cref{thm:main2}. In \cref{sec:app}, we verify the tail reducers displayed in \cref{fig:ptr,fig:ntr}. In \cref{sec:lean-details}, we state precisely what is formalized and describe the simplifications used in translating the mathematical argument into Lean.

\section{Other common tools} \label{sec:tools}

For readability, we use lowercase Latin letters for vertices and lowercase Greek letters for edges throughout the remainder of the paper. We begin with the following convenient terminology.

\begin{definition}[Red, reddish, blue and bluish] \label{def:red-blue}
  Given a cut $(A, B)$ of a graph $G$, we say that the edges of $G[A]$ are \emph{red}, and the edges of $G[B]$ are \emph{blue}, and moreover, a vertex in $A$ is \emph{red} if it is incident with a red edge, and \emph{reddish} otherwise, and a vertex in $B$ is \emph{blue} if it is incident with a blue edge, and \emph{bluish} otherwise.
\end{definition}

Pictorially, we draw
\scalebox{0.4}{%
\begin{tikzpicture}
  \draw[red-edge, line width=2pt] (0,0) -- (1,0);
  \node[red-vertex] at (0,0) {};
  \node[red-vertex] at (1,0) {};
\end{tikzpicture}
}
for a red edge whose endpoints are red, a pastel red circle
\scalebox{0.4}{%
\begin{tikzpicture}
  \node[reddish-vertex] at (0,0) {};
\end{tikzpicture}
}
for a reddish vertex, a red-outlined circle
\scalebox{0.4}{%
\begin{tikzpicture}
  \node[a-vertex] at (0,0) {};
\end{tikzpicture}
} for a vertex in $A$,
\scalebox{0.4}{%
\begin{tikzpicture}
  \draw[blue-edge, line width=2pt] (0,0) -- (1,0);
  \node[blue-vertex] at (0,0) {};
  \node[blue-vertex] at (1,0) {};
\end{tikzpicture}
}
for a blue edge whose endpoints are blue, a pastel blue square
\scalebox{0.4}{%
\begin{tikzpicture}
  \node[bluish-vertex] at (0,0) {};
\end{tikzpicture}
}
for a bluish vertex, and a blue-outlined square
\scalebox{0.4}{%
\begin{tikzpicture}
  \node[b-vertex] at (0,0) {};
\end{tikzpicture}
} for a vertex in $B$.

We emphasize that the above terminology is \emph{context dependent} --- as we flip a cut preserver $C$, some vertices and edges might get recolored according to the new cut $(A \oplus C, B \oplus C)$.

\begin{proposition} \label{prop:cut-enhancer}
  For every cut $(A, B)$ of a subcubic graph $G$, if any of the graphs in \cref{fig:cut-enhancers} is an induced subgraph $F$ of $G$, then the vertices of $F$ form a cut enhancer.
\end{proposition}

\begin{figure}
  \centering
  \begin{csubfig}{2.8cm}{cut-enhancer-a}
    \defc{a/0/1,b/-1/0,c/1/0}
    \drawe{}{}{a/b,a/c}
    \node[red-vertex] at (a) {$\uparrow$};
    \node[blue-vertex] at (b) {};
    \node[blue-vertex] at (c) {};
  \end{csubfig}%
  \begin{csubfig}{2.8cm}{cut-enhancer-b}
    \defc{a/0/1,b/-1/0,c/-1/1,d/1/0,e/1/1}
    \drawe{}{}{a/b,a/d,b/c,d/e}
    \node[reddish-vertex] at (a) {$\downarrow$};
    \node[blue-vertex] at (b) {$\uparrow$};
    \node[red-vertex] at (c) {};
    \node[blue-vertex] at (d) {$\uparrow$};
    \node[red-vertex] at (e) {};
  \end{csubfig}%
  \begin{csubfig}{2.8cm}{cut-enhancer-c}
    \defc{a/0/1,b/1/1,c/0/0,d/1/0,e/2/1}
    \drawe{}{}{a/c,a/d,b/d,d/e}
    \node[red-vertex] at (a) {$\uparrow$};
    \node[red-vertex] at (b) {};
    \node[blue-vertex] at (c) {};
    \node[bluish-vertex] at (d) {};
    \node[red-vertex] at (e) {};
  \end{csubfig}%
  \begin{csubfig}{2.8cm}{cut-enhancer-d}
    \defc{a/0/1,b/1/1,c/0/0,d/1/0,e/2/1}
    \drawe{}{}{a/c,a/d,b/c,b/d,d/e}
    \node[red-vertex] at (a) {$\uparrow$};
    \node[reddish-vertex] at (b) {$\downarrow$};
    \node[blue-vertex] at (c) {$\uparrow$};
    \node[bluish-vertex] at (d) {};
    \node[red-vertex] at (e) {};
  \end{csubfig}

  \par\bigskip

  \begin{csubfig}{3.6cm}{cut-enhancer-e}
    \defc{a/-1/1,b/0/1,c/1/1,d/2/1,e/0/0,f/1/0,g/2/0}
    \drawe{}{}{a/e,b/e,b/g,c/e,c/f,d/f,d/g}
    \node[red-vertex] at (a) {};
    \node[reddish-vertex] at (b) {$\downarrow$};
    \node[red-vertex] at (c) {$\uparrow$};
    \node[reddish-vertex] at (d) {$\downarrow$};
    \node[bluish-vertex] at (e) {};
    \node[blue-vertex] at (f) {$\uparrow$};
    \node[blue-vertex] at (g) {$\uparrow$};
  \end{csubfig}%
  \begin{csubfig}{3.6cm}{cut-enhancer-f}
    \defc{a/-1/1,b/0/1,c/1/1,d/2/1,e/0/0,f/1/0,g/2/0}
    \drawe{}{}{a/e,b/e,b/f,c/f,c/g,d/g}
    \node[red-vertex] at (a) {};
    \node[reddish-vertex] at (b) {$\downarrow$};
    \node[reddish-vertex] at (c) {$\downarrow$};
    \node[red-vertex] at (d) {};
    \node[blue-vertex] at (e) {$\uparrow$};
    \node[blue-vertex] at (f) {$\uparrow$};
    \node[blue-vertex] at (g) {$\uparrow$};
  \end{csubfig}%
  \begin{csubfig}{3.6cm}{cut-enhancer-g}
    \defc{a/-1/1,b/0/1,c/1/1,d/2/1,e/0/0,f/1/0,g/2/0}
    \drawe{}{}{a/e,b/e,c/e,c/f,c/g,d/g}
    \node[red-vertex] at (a) {};
    \node[red-vertex] at (b) {};
    \node[reddish-vertex] at (c) {};
    \node[red-vertex] at (d) {};
    \node[bluish-vertex] at (e) {};
    \node[blue-vertex] at (f) {};
    \node[blue-vertex] at (g) {$\uparrow$};
  \end{csubfig}
  \caption{Cut enhancers.} \label{fig:cut-enhancers}
\end{figure}

\begin{proof}
  For each vertex $v$, let $\delta(v)$ be the number of edges incident to $v$ in the cut-set of $(A \oplus C, B \oplus C)$ minus the number of those edges in the cut-set of $(A,B)$. We label each vertex $v$ above by $\uparrow$ if $\delta(v) \ge 1$, nothing if $\delta(v) \ge 0$, and $\downarrow$ if $\delta(v) \ge -1$. Notice that there are more $\uparrow$s than $\downarrow$s.
\end{proof}

We shall repeatedly use the following immediate consequence of \ce{a} and \cref{lem:max-cut} to deduce colors of vertices.

\begin{lemma} \label{lem:blue-edge}
  For every maximum cut $(A, B)$ of a subcubic graph $G$, if $a$ is a red vertex, $bc$ is a blue edge, and both $b$ and $d$ are the neighbors of $a$ in $B$, then $d$ is bluish.
  \begin{ctikzpq}
    \begin{scope}
      \defc{a/0/1,b/0/0,c/-1/0,d/1/0}
      \drawe{}{b/c}{a/b,a/d}
      \drawv{a}{}{b,c}{}
      \node[b-vertex] at (d) {$d$};
    \end{scope}
    \draw[double,-stealth,thin] (1.5,0.5) -- (2.5,0.5);
    \begin{scope}[shift={(4,0)}]
      \defc{a/0/1,b/0/0,c/-1/0,d/1/0}
      \drawe{}{b/c}{a/b,a/d}
      \drawv{a}{}{b,c}{d}
    \end{scope}
  \end{ctikzpq}
\end{lemma}

For certain cut enhancers used only in this section, we require degree constraints on vertices in the ambient graph $G$. To indicate such a constraint, we give the red circle lower corners: one lower-left corner indicates degree $1$, while two lower corners indicate degree $2$.

\begin{proposition} \label{lem:degree-cut-enhancers}
  For every cut $(A, B)$ of a subcubic graph $G$, if any of the graphs in \cref{fig:degree-cut-enhancers} is an induced subgraph $F$ of $G$ satisfying the indicated degree constraint, then the vertices of $F$ form a cut enhancer.
\end{proposition}

\begin{figure}[b]
  \centering
  \begin{csubfig}{2.8cm}{degree-cut-enhancer-a}
    \defc{a/0/1,b/0/0}
    \node[red-degree-one-vertex] at (a) {$\uparrow$};
    \node[blue-vertex, opacity=0] at (b) {};
  \end{csubfig}%
  \begin{csubfig}{2.8cm}{degree-cut-enhancer-b}
    \defc{a/0/1,b/0/0}
    \drawe{}{}{a/b}
    \node[red-degree-two-vertex] at (a) {$\uparrow$};
    \node[blue-vertex] at (b) {};
  \end{csubfig}
  \caption{Cut enhancers with degree constraints --- the red vertices are of degree $1$ and $2$, respectively.} \label{fig:degree-cut-enhancers}
\end{figure}

\begin{proof}
  For each vertex $v$ of $F$, let $\delta(v)$ be the change in the number of cut-set edges incident to $v$ when $V(F)$ is flipped. The symbol $\uparrow$ means that $\delta(v) \ge 1$, while no arrow means that $\delta(v) \ge 0$. Thus, flipping $V(F)$ increases the size of the cut-set by at least $1$.
\end{proof}

Only the following degree-constrained absolute tail reducer is used in the paper, and only in this section.

\begin{lemma} \label{lem:degree-absolute-tail-reducer}
  For every cut $(A, B)$ of a subcubic graph $G$, if the graph in \cref{fig:degree-absolute-tail-reducer} is an induced subgraph of $G$ satisfying the indicated degree constraint, then $\sset{a}$ is an absolute tail reducer.
\end{lemma}

\begin{figure}
  \centering
  \begin{tikzpicture}
    \defc{a/0/1,b/0/0}
    \drawe{}{}{a/b}
    \node[red-degree-two-vertex] at (a) {$a$};
    \node[bluish-vertex] at (b) {$b$};
  \end{tikzpicture}
  \caption{An absolute tail reducer with a degree constraint --- the red vertex $a$ is of degree $2$.} \label{fig:degree-absolute-tail-reducer}
\end{figure}

\begin{proof}
  The connected component of $G[\sset{a} \cup B]$ containing $a$ is the edge $ab$, and $t^\pm_{G[\sset{a,b}]}=0<1$. Thus $\sset{a}$ is an absolute tail reducer.
\end{proof}

The following two lemmas impose degree constraints on subcubic graphs that lack absolute tail reducers.

\begin{lemma} \label{lem:deg-3}
  For every maximum cut $(A, B)$ of a subcubic graph $G$, if $ab$ is a red edge, then $a$ is of degree $3$ in $G$, or $\sset{a}$ is an absolute tail reducer.
\end{lemma}

\begin{proof}
  Suppose $a$ is of degree at most $2$. \Cref{fig:degree-cut-enhancer-a} implies that $a$ is of degree $2$, because otherwise $\sset{a}$ is a cut enhancer. Let $c$ be the second neighbor of $a$. \cref{lem:max-cut} implies that $c \in B$, and \Cref{fig:degree-cut-enhancer-b} implies that $c$ is bluish, because otherwise $\sset{a,c}$ is a cut enhancer. Thus $\sset{a}$ is an absolute tail reducer by \cref{lem:degree-absolute-tail-reducer}.
  \begin{ctikzq}
    \defc{a/0/1,b/1/1,c/0/0}
    \drawe{a/b}{}{a/c}
    \drawv{a,b}{}{}{c}
  \end{ctikzq}
\end{proof}

\begin{lemma} \label{lem:deg-3-flip}
  For every maximum cut $(A, B)$ of a subcubic graph $G$, if $a$ is a red vertex of degree $3$, $b$ is a blue vertex of degree $3$, $c$ is a reddish vertex, $a \sim b$, and $b \sim c$, then $c$ is of degree $3$ in $G$, or $\sset{a, c}$ is an absolute tail reducer.
\end{lemma}

\begin{proof}
  \cref{prop:cut-preserver} implies that $\sset{a,b}$ is a cut preserver. Flip the cut preserver $\sset{a,b}$. Since $c$ becomes red, \cref{lem:deg-3} implies that $c$ is of degree $3$, or $\sset{c}$ becomes an absolute tail reducer. In the latter case, $\sset{a,c}$ is an absolute tail reducer with respect to $(A, B)$ by \cref{lem:mod-cut-preserver} and the remark following its proof on \cpageref{rmk:mod-cut-preserver}.
  \begin{ctikzq}
    \defc{a/0/1,b/0/0,c/1/1}
    \drawe{}{}{a/b,b/c}
    \drawv{a}{c}{b}{}
    \draw[-stealth, thin] (1.5,0.5) -- (3.5,0.5) node[midway, above, font=\footnotesize] {$\sset{a,b}$};
    \begin{scope}[shift={(4,0)}]
      \defc{a/0/1,b/0/0,c/1/1}
      \drawe{b/c}{}{a/b}
      \drawv{b,c}{}{a}{}
    \end{scope}
  \end{ctikzq}
\end{proof}

We shall repeatedly use the following results on induced paths and induced even cycles in $M$.

\begin{lemma} \label{lem:induced-path}
  For every maximum cut $(A,B)$ of a subcubic graph $G$, if the underlying multigraph contains an induced path $\alpha_1 \dots \alpha_n$, then $G$ contains an induced path $a_1b_1 \dots a_nb_n$ such that $\alpha_1 = a_1b_1, \dots, \alpha_n = a_nb_n$.
\end{lemma}

\begin{proof}
  Let $\alpha_1 = a_1b_1$ and $\alpha_2 = a_2b_2$ be such that $G[\alpha_1,\alpha_2]$ consists of the single edge $b_1a_2$. To find the next $a_i, b_i$, let $\alpha_i = a_ib_i$ such that either $a_i \sim a_{i-1}$ or $a_i \sim b_{i-1}$. \cref{lem:blue-edge} implies that $a_i \nsim a_{i-1}$ because $a_{i-1} \sim b_{i-2}$, and so $a_i \sim b_{i-1}$.
\end{proof}

\begin{lemma} \label{lem:induced-cycle}
  For every maximum cut $(A,B)$ of a subcubic graph $G$, if the underlying multigraph contains an induced even cycle $\alpha_1 \dots \alpha_{2n}$, then $G$ contains an induced cycle $a_1b_1 \dots a_{2n}b_{2n}$ such that $\alpha_1 = a_1b_1, \dots, \alpha_n = a_nb_n$.
\end{lemma}

\begin{proof}
  We find $a_1, b_1, \dots, a_{2n}, b_{2n}$ the same way as in the proof of \cref{lem:induced-path}. Finally, \cref{lem:blue-edge} implies that $a_{2n} \nsim a_1$ and $a_{2n} \nsim b_1$ because $a_{2n} \sim b_{2n-1}$, and so either $b_{2n} \sim a_1$ or $b_{2n} \sim b_1$, and furthermore \cref{lem:blue-edge} implies that $b_{2n} \nsim b_1$ because $b_1 \sim a_2$, and so $b_{2n} \sim a_1$.
\end{proof}

\section{Positive tail reducers} \label{sec:pos}

We collect all the positive tail reducers ever needed in the following proposition.

\begin{figure}
  \centering
  \begin{csubfig}{2.0cm}{ptr-a}
    \defc{a/0/1,b/1/1,c/0/0,d/1/0}
    \drawe{a/b}{c/d}{a/c,a/d,b/c,b/d}
    \drawv{a,b}{}{c,d}{}
  \end{csubfig}
  \begin{csubfig}{2.0cm}{ptr-b}
    \defc{a/0/1,b/1/1,c/0/0,d/1/0}
    \drawe{a/b}{}{a/c,a/d,b/c,b/d}
    \drawv{a,b}{}{}{c,d}
  \end{csubfig}
  \begin{csubfig}{2.8cm}{ptr-c}
    \defc{c/0/0,d/1/0,a/0/1,b/1/1,e/2/0}
    \drawe{}{c/d}{a/c,b/d,a/e,b/e}
    \drawv{a,b}{}{c,d}{e}
  \end{csubfig}
  \begin{csubfig}{2.8cm}{ptr-d}
    \defc{a/0/1,b/1/1,c/0/0,d/1/0,e/2/0}
    \drawe{a/b}{c/d}{a/c,a/d,b/d,b/e}
    \drawv{a,b}{}{c,d}{e}
  \end{csubfig}
  \begin{csubfig}{2.8cm}{ptr-e}
    \defc{c/0/0,d/1/0,a/0/1,b/1/1,e/2/0}
    \drawe{a/b}{c/d}{a/c,b/c,a/e,b/e}
    \drawv{a,b}{}{c,d}{e}
  \end{csubfig}
  \begin{csubfig}{2.8cm}{ptr-f}
    \defc{a/0/1,b/1/1,c/0/0,d/1/0,e/2/0}
    \drawe{a/b}{c/d}{a/e,a/c,b/e,b/d}
    \drawv{a,b}{}{c,d}{e}
  \end{csubfig}

  \par\bigskip

  \begin{csubfig}{2.8cm}{ptr-g}
    \defc{a/0/1,b/1/1,d/1/0,e/2/0,c/0/0}
    \drawe{}{}{a/d,a/c,b/d,b/e}
    \drawv{a,b}{}{}{c,d,e}
  \end{csubfig}
  \begin{csubfig}{2.8cm}{ptr-h}
    \defc{a/0/1,b/1/1,c/0/0,d/1/0,e/2/0}
    \drawe{a/b}{}{a/c,a/d,b/d,b/e}
    \drawv{a,b}{}{}{c,d,e}
  \end{csubfig}
  \begin{csubfig}{3.6cm}{ptr-i}
    \defc{a/0/1,b/1/1,d/0/0,e/1/0,f/2/0,c/-1/0}
    \drawe{a/b}{d/e}{a/d,a/c,b/f,b/e}
    \drawv{a,b}{}{d,e}{c,f}
  \end{csubfig}
  \begin{csubfig}{3.6cm}{ptr-j}
    \defc{d/0/0,e/1/0,a/0/1,b/1/1,f/2/0,c/-1/0}
    \drawe{a/b}{d/e}{a/d,b/d,a/c,b/f}
    \drawv{a,b}{}{d,e}{c,f}
  \end{csubfig}
  \begin{csubfig}{3.6cm}{ptr-k}
    \defc{a/0/1,b/1/1}
    \defc{c/-1/0,d/0/0,e/1/0,f/2/0}
    \drawe{a/b}{}{a/c,a/d,b/e,b/f}
    \drawv{a,b}{}{}{c,d,e,f}
  \end{csubfig}

  \par\bigskip

  \begin{csubfig}{3.6cm}{ptr-l}
    \defc{a/0/1,b/1/1,c/-1/0,d/0/0,e/1/0,f/2/0}
    \drawe{}{e/f}{a/c,a/d,a/e,b/d,b/e}
    \drawv{b}{a}{e,f}{c,d}
  \end{csubfig}
  \begin{csubfig}{3.6cm}{ptr-m}
    \defc{a/0/1,b/1/1,c/-1/0,d/0/0,e/1/0,f/2/0}
    \drawe{}{e/f}{a/c,a/d,a/e,b/d,b/f}
    \drawv{b}{a}{e,f}{c,d}
  \end{csubfig}
  \begin{csubfig}{3.6cm}{ptr-n}
    \defc{a/0/1,b/1/1,c/2/1,d/-1/0,e/0/0,f/1/0,g/2/0}
    \drawe{a/b}{f/g}{a/d,a/e,b/d,b/f,c/e,c/g}
    \drawv{a,b,c}{}{f,g}{d,e}
  \end{csubfig}
  \begin{csubfig}{3.6cm}{ptr-o}
    \defc{a/0/1,b/1/1,c/2/1,e/0/0,f/1/0,g/2/0,d/-1/0}
    \drawe{b/c}{d/e}{a/d,a/e,a/f,b/g,b/e,c/f,c/g}
    \drawv{b,c}{a}{d,e}{f,g}
  \end{csubfig}

  \par\bigskip

  \begin{csubfig}{4.4cm}{ptr-p}
    \defc{a/0/1,b/1/1,c/2/1,e/0/0,f/1/0,g/2/0,h/3/0,d/-1/0}
    \drawe{}{e/f}{a/d,a/e,a/f,b/e,b/g,c/f,c/h}
    \drawv{c,b}{a}{e,f}{d,g,h}
  \end{csubfig}
  \begin{csubfig}{4.4cm}{ptr-q}
    \defc{a/0/1,b/1/1,c/2/1,e/0/0,f/1/0,g/2/0,h/3/0,d/-1/0}
    \drawe{b/c}{e/f}{a/d,a/e,a/f,b/f,b/g,c/h,c/g}
    \drawv{b,c}{a}{e,f}{d,g,h}
  \end{csubfig}
  \begin{csubfig}{5.2cm}{ptr-r}
    \defc{b/1/1,a/0/1,c/2/1,e/0/0,d/-1/0,g/2/0,h/3/0,i/4/0,f/1/0}
    \drawe{}{g/h,f/e}{b/f,b/d,a/e,a/g,a/d,c/g,c/i}
    \drawv{b,c}{a}{f,e,g,h}{d,i}
  \end{csubfig}

  \par\bigskip

  \begin{csubfig}{5.2cm}{ptr-s}
    \defc{b/1/1,a/0/1,c/2/1,e/0/0,d/-1/0,g/2/0,h/3/0,i/4/0,f/1/0}
    \drawe{}{g/h,f/e}{b/f,b/d,a/e,a/g,a/d,c/h,c/i}
    \drawv{b,c}{a}{f,e,g,h}{d,i}
  \end{csubfig}
  \begin{csubfig}{5.2cm}{ptr-t}
    \defc{a/0/1,b/1/1,c/2/1,f/0/0,g/1/0,h/2/0,i/3/0,e/-1/0,d/-2/0}
    \drawe{b/c}{f/g}{a/d,a/e,a/g,b/g,b/h,c/h,c/i}
    \drawv{b,c}{a}{f,g}{d,e,h,i}
  \end{csubfig}
  \begin{csubfig}{6.0cm}{ptr-u}
    \defc{a/0/1,b/1/1,c/2/1,f/0/0,g/1/0,h/2/0,i/3/0,j/4/0,e/-1/0,d/-2/0}
    \drawe{b/c}{h/i}{a/d,a/e,a/h,b/f,b/g,c/j,c/h}
    \drawv{b,c}{a}{h,i}{d,e,f,g,j}
  \end{csubfig}

  \par\medskip

  \begin{csubfig}{5.2cm}{ptr-v}
    \defc{a/0/1,b/1/1,c/2/1,d/3/1,e/0/0,f/1/0,g/2/0,h/3/0,i/4/0,j/5/0}
    \drawe{a/b,c/d}{f/g,h/i}{b/f,c/g,d/h,a/e,b/e,d/j}
    \draw (a) to [bend right] (2,0.5) to [bend left] (i);
    \draw (c) .. controls (3.14,1.6) .. (j);
    \drawv{a,b,c,d}{}{f,g,h,i}{e,j}
  \end{csubfig}
  \begin{csubfig}{6.0cm}{ptr-w}
    \defc{a/0/1,b/1/1,c/2/1,d/3/1,e/-1/0,f/0/0,g/1/0,h/2/0,i/3/0,j/4/0,k/5/0}
    \drawe{a/b,c/d}{g/h,i/j}{b/g,c/h,d/i,a/e,b/f,d/k}
    \draw (a) to [bend right] (2,0.5) to [bend left] (j);
    \draw (c) .. controls (3.14,1.6) .. (k);
    \drawv{a,b,c,d}{}{g,h,i,j}{e,f,k}
  \end{csubfig}

  \par\medskip

  \begin{csubfig}{6.0cm}{ptr-x}
    \defc{a/0/1,b/1/1,c/2/1,d/3/1,h/1/0,i/2/0,j/3/0,k/4/0,g/0/0,f/-1/0,e/-2/0}
    \drawe{a/b}{h/i,j/k}{a/e,a/g,b/f,b/h,c/g,c/j,d/h,d/i,d/k}
    \drawv{a,b,c}{d}{h,i,j,k}{e,f,g}
  \end{csubfig}
  \begin{csubfig}{6.8cm}{ptr-y}
    \defc{a/0/1,b/1/1,c/2/1,d/3/1,e/-1/0,f/0/0,g/1/0,h/2/0,i/3/0,j/4/0,k/5/0,l/6/0}
    \drawe{a/b,c/d}{g/h,i/j}{b/g,c/h,d/i,a/e,b/f,d/k}
    \draw (a) to [bend right] (2,0.5) to [bend left] (j);
    \draw (c) .. controls (3,1.6) .. (l);
    \drawv{a,b,c,d}{}{g,h,i,j}{e,f,k,l}
  \end{csubfig}

  \caption{Positive tail reducers.} \label{fig:ptr}
\end{figure}

\begin{figure}
  \centering
  \begin{csubfig}{3.6cm}{degree-ptr-a}
    \defc{a/0/1,b/1/1,c/0/0,d/1/0,e/2/0}
    \drawe{}{d/e}{a/c,a/d,b/c,b/e}
    \node[reddish-degree-two-vertex] at (a) {$a$};
    \drawv{b}{}{d,e}{c}
  \end{csubfig}
  \caption{A positive tail reducer with a degree constraint --- the reddish vertex $a$ is of degree $2$.} \label{fig:degree-ptr}
\end{figure}

\begin{proposition} \label{prop:ptr}
  For every maximum cut $(A, B)$ of a subcubic graph $G$, if either one of the graphs in \cref{fig:ptr} is an induced subgraph $F$ of $G$, or the graph in \cref{fig:degree-ptr-a} is an induced subgraph $F$ of $G$ satisfying the indicated degree constraint, then $A \cap V(F)$ is a positive tail reducer.
\end{proposition}

\begin{proof}
  Let $D = A \cap V(F)$. Suppose that $F$ is one of the graphs in \cref{fig:ptr}. Each reddish vertex has degree $3$ in $F$, each red vertex has two neighbors in $B \cap V(F)$, and each blue vertex has its blue neighbor in $F$. By subcubicity and \cref{lem:max-cut}, no edge of $G[D \cup B]$ joins $V(F)$ to its complement. Thus $F$ is the connected component of $G[D \cup B]$ containing $D$. The inequality $t_F^+ < \abs{D}$ is verified in \cref{sec:app}.

  Finally suppose that $F$ is the graph in \cref{fig:degree-ptr-a}. The degree constraint implies, by the same argument, that $F$ is the connected component of $G[D \cup B]$ containing $D$. Moreover, $F$ is obtained from the graph $F'$ in \cref{fig:ptr-m} by deleting vertex $c$. By the computation in \cref{sec:app} and the Cauchy interlacing theorem, $t_F^+ \le t_{F'}^+ < 2 = \abs{D}$.
\end{proof}

Throughout the rest of this section, we assume that $(A, B)$ is a maximum cut of a subcubic graph $G$, and $M$ is the underlying multigraph of $G$ with respect to $(A, B)$. When the conclusion of a result in this section is of the form ``a subset of $U$ is a positive tail reducer'' for a specific vertex subset $U$ of $G$, in view of \cref{lem:deg-3}, we always assume from the start of the proof that every red or blue vertex in $U$ is of degree $3$.

Roughly speaking, we establish, through the next three lemmas, that every connected component of $M$ is a path or a cycle of even length if $G$ lacks positive tail reducers.

\begin{lemma} \label{lem:no-multi}
  If $M$ contains multiple edges between $\alpha$ and $\beta$, then a subset of $\alpha \cup \beta$ is a positive tail reducer.
\end{lemma}

\begin{proof}
  Let $\alpha = ab$ and $\beta = cd$. In view of \cref{lem:multigraph-bipartite}, we may assume that $ab$ is red, and $cd$ is blue up to symmetry.
  \begin{kase}
    \item If $M$ contains $4$ edges between $ab$ and $cd$, then $\sset{a,b}$ is a positive tail reducer by \cref{fig:ptr-a}.
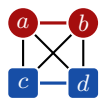
    \begin{ctikz}
      \defc{a/0/1,b/1/1,c/0/0,d/1/0}
      \drawe{a/b}{c/d}{a/c,a/d,b/c,b/d}
      \drawv{a,b}{}{c,d}{}
    \end{ctikz}
    \item Suppose $M$ contains $3$ edges between $ab$ and $cd$. We may assume that $a \sim c, a \sim d$, and $b \sim d$ up to symmetry. Let $e$ be the third neighbor of $b$. \cref{lem:blue-edge} implies that $e$ is bluish, and so $\sset{a,b}$ is a positive tail reducer by \cref{fig:ptr-d}.
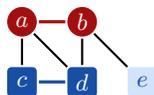
    \begin{ctikz}
      \defc{a/0/1,b/1/1,c/0/0,d/1/0,e/2/0}
      \drawe{a/b}{c/d}{a/c,a/d,b/d,b/e}
      \drawv{a,b}{}{c,d}{e}
    \end{ctikz}
    \item Suppose $M$ contains $2$ edges between $ab$ and $cd$.
    \begin{kase}
      \item Suppose the $2$ edges between $ab$ and $cd$ are vertex-disjoint. We may assume that $a \sim c$ and $b \sim d$ up to symmetry. \cref{lem:blue-edge} implies that the third neighbor of $a$ and that of $b$ are bluish. Depending on whether $a$ and $b$ share a bluish neighbor, there are two possibilities. In either case, $\sset{a,b}$ is a positive tail reducer by \cref{fig:ptr-f,fig:ptr-i}.
      \begin{ctikz}
        \begin{scope}
          \defc{a/0/1,b/1/1,c/0/0,d/1/0,e/2/0}
          \drawe{a/b}{c/d}{a/e,a/c,b/e,b/d}
          \drawv{a,b}{}{c,d}{e}
        \end{scope}
        \begin{scope}[shift={(4.5,0)}]
          \defc{a/0/1,b/1/1,c/0/0,d/1/0,e/2/0,f/-1/0}
          \drawe{a/b}{c/d}{a/f,a/c,b/e,b/d}
          \drawv{a,b}{}{c,d}{e,f}
        \end{scope}
        \end{ctikz}
      \item Suppose the $2$ edges between $ab$ and $cd$ share a vertex. We may assume that $a \sim c$ and $a \sim d$ up to symmetry. \cref{lem:blue-edge} implies that the third neighbor of $c$ and that of $d$ are reddish. Depending on whether $c$ and $d$ share a reddish neighbor, there are two possibilities. In either case, $\sset{c,d}$ is a positive tail reducer by \cref{fig:ptr-e,fig:ptr-j}.
      \begin{ctikzq}
        \begin{scope}
          \defc{a/0/1,b/1/1,c/0/0,d/1/0,e/2/1}
          \drawe{a/b}{c/d}{a/c,a/d,c/e,d/e}
          \drawv{a,b}{e}{c,d}{}
        \end{scope}
        \begin{scope}[shift={(4.5,0)}]
          \defc{a/0/1,b/1/1,c/0/0,d/1/0,e/2/1,f/-1/1}
          \drawe{a/b}{c/d}{a/c,a/d,c/f,d/e}
          \drawv{a,b}{e,f}{c,d}{}
        \end{scope}
      \end{ctikzq}
    \end{kase}
  \end{kase}
\end{proof}

\begin{lemma} \label{lem:m-deg-2}
  If $\alpha$ is a vertex of $M$, then $\alpha$ is of degree at most $2$ in $M$, or a subset of $G(\alpha,1)$ is a positive tail reducer.
\end{lemma}

\begin{proof}
  Without loss of generality, we may assume that $\alpha = ab$ is red. \cref{lem:blue-edge} implies that neither $a$ nor $b$ can be adjacent to two blue edges in $G$. If $\alpha$ is of degree at least $3$ in $M$, then $M$ contains multiple edges between $\alpha$ and its neighbor $\beta$, and so a subset of $\alpha \cup \beta$ is a positive tail reducer by \cref{lem:no-multi}. Finally, notice that the distance between every vertex of $\beta$ and $\alpha$ is $1$.
\end{proof}

\begin{lemma} \label{lem:isolated}
  If $\alpha$ is an isolated vertex of $M$, then a vertex subset of $\alpha$ is a positive tail reducer.
\end{lemma}

\begin{proof}
  Without loss of generality, we may assume that $\alpha = ab$ is red. Since $ab$ is isolated in $M$, the other neighbors of $a$ and $b$ are bluish. Depending on how these bluish neighbors overlap, there are three possibilities. In each case, $\sset{a,b}$ is a positive tail reducer by \cref{fig:ptr-b,fig:ptr-h,fig:ptr-k}.
  \begin{ctikzq}
    \begin{scope}
      \defc{a/0/1,b/1/1,c/0/0,d/1/0}
      \drawe{a/b}{}{a/c,a/d,b/c,b/d}
      \drawv{a,b}{}{}{}
      \drawvnl{}{}{}{c,d}
    \end{scope}
    \begin{scope}[shift={(3,0)}]
      \defc{a/0/1,b/1/1,c/-.5/0,d/.5/0,e/1.5/0}
      \drawe{a/b}{}{a/c,a/d,b/d,b/e}
      \drawv{a,b}{}{}{}
      \drawvnl{}{}{}{c,d,e}
    \end{scope}
    \begin{scope}[shift={(7,0)}]
      \defc{a/0/1,b/1/1}
      \defc{c/-1/0,d/0/0,e/1/0,f/2/0}
      \drawe{a/b}{}{a/c,a/d,b/e,b/f}
      \drawv{a,b}{}{}{}
      \drawvnl{}{}{}{c,d,e,f}
    \end{scope}
  \end{ctikzq}
\end{proof}

In the remaining proofs of this section, we shall obtain subsets $C$ and $D$ of $U$, for a specific vertex subset $U$ of $G$, such that $C$ is a cut preserver, and $D$ is a positive tail reducer with respect to $(A \oplus C, B \oplus C)$. We always automatically apply \cref{lem:mod-cut-preserver} to obtain a subset of $U$ that is a positive tail reducer with respect to $(A, B)$.

\begin{lemma} \label{lem:path-5-8}
  If $M(\alpha, 4)$ is a path on at least $5$ vertices, or a cycle on at least $6$ vertices, then a subset of $G(\alpha, 8)$ is a positive tail reducer.
\end{lemma}

\begin{proof}
  Observe that $M(\alpha, 4)$ contains an induced path on $5$ vertices, say $ab, cd, ef, gh, ij$, such that $\alpha \in \sset{ab,cd,ef,gh,ij}$. In view of \cref{lem:induced-path}, we may assume that $abcdefghij$ is an induced path in $G$, and moreover we may assume that $ab, ef, ij$ are red, and $cd, gh$ are blue. Let $k$ be the third neighbor of $e$. \cref{lem:blue-edge} implies that $k$ is bluish.
  \begin{ctikz}
    \defc{a/0/1,b/1/1,c/1/0,d/2/0,e/2/1,f/3/1,g/3/0,h/4/0,i/4/1,j/5/1,k/0/0}
    \drawe{a/b,e/f,i/j}{c/d,g/h}{b/c,d/e,e/k,f/g,h/i}
    \drawv{a,b,e,f,i,j}{}{c,d,g,h}{k}
  \end{ctikz}

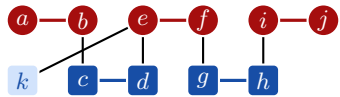
  \begin{kase}
    \item If $k \sim b$, then $\sset{b,e}$ is a positive tail reducer by \cref{fig:ptr-c}.
    \begin{ctikz}
      \defc{a/0/1,b/1/1,c/1/0,d/2/0,e/2/1,f/3/1,g/3/0,h/4/0,i/4/1,j/5/1,k/0/0}
      \drawe{a/b,e/f,i/j}{c/d,g/h}{b/c,d/e,e/k,f/g,h/i,k/b}
      \drawv{a,b,e,f,i,j}{}{c,d,g,h}{k}
    \end{ctikz}
    \item Suppose $k \nsim b$. \ce{b} implies that $k \nsim i$ because otherwise $\sset{d,e,h,i,k}$ is a cut enhancer. Flip the cut preserver $\sset{b,c}$, and then flip the cut preserver $\sset{h,i}$. The two neighbors of $e$ in $B$ become bluish. By a symmetric argument, we may assume that the two neighbors of $f$ become bluish as well. Then the red edge $ef$ becomes isolated in $M$, and so a subset of $\sset{e,f}$ becomes a positive tail reducer by \cref{lem:isolated}.
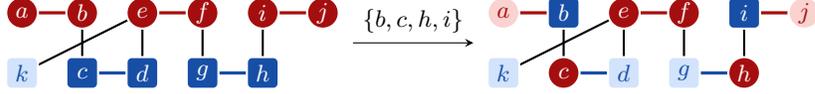
    \begin{ctikz}
      \defc{a/0/1,b/1/1,c/1/0,d/2/0,e/2/1,f/3/1,g/3/0,h/4/0,i/4/1,j/5/1,k/0/0}
      \drawe{a/b,e/f,i/j}{c/d,g/h}{b/c,d/e,e/k,f/g,h/i}
      \drawv{a,b,e,f,i,j}{}{c,d,g,h}{k}
      \draw[-stealth, thin] (5.5,0.5) -- (6.5,0.5) node[midway, above, font=\footnotesize] {$\sset{b,c}$};
      \node at (7,0.5) {$\dots$};
      \draw[-stealth, thin] (7.5,0.5) -- (8.5,0.5) node[midway, above, font=\footnotesize] {$\sset{h,i}$};
      \begin{scope}[shift={(9,0)}]
        \defc{a/0/1,b/1/1,c/1/0,d/2/0,e/2/1,f/3/1,g/3/0,h/4/0,i/4/1,j/5/1,k/0/0}
        \drawe{a/b,e/f,i/j}{c/d,g/h}{b/c,d/e,e/k,f/g,h/i}
        \drawv{c,e,f,h}{a,j}{b,i}{d,g,k}
      \end{scope}
    \end{ctikz}
  \end{kase}

  Finally, notice that the distance between $a$ and $ij$ is $8$, and it is the furthest distance between any vertex in a cut preserver or a positive tail reducer and $\alpha \in \sset{ab, cd, ef, gh, ij}$.
\end{proof}

\begin{lemma} \label{lem:path-5-subgraph-9}
  If $M$ contains a path $\alpha\beta\gamma\delta\epsilon$ as a subgraph, then a subset of $G(\gamma,9)$ is a positive tail reducer.
\end{lemma}

\begin{proof}
  \begin{kase}
    \item If $M(\gamma,4)$ contains multiple edges, then a subset of $G(\gamma,8)$ is a positive tail reducer by \cref{lem:no-multi}.
    \item If $M(\gamma,4)$ is a simple graph of maximum degree at least $3$, then a subset of $G(\gamma,9)$ is a positive tail reducer by \cref{lem:m-deg-2}.
    \item If $M(\gamma,4)$ is a simple graph of maximum degree at most $2$, then $M(\gamma,4)$ must be a path on at least $5$ vertices, or any even cycle on at least $6$ vertices, and so a subset of $G(\gamma,8)$ by \cref{lem:path-5-8}. \qedhere
  \end{kase}
\end{proof}

\begin{lemma} \label{lem:cycle-4}
  If $M(\alpha, 4)$ is a cycle on $4$ vertices, then a subset of $G(\alpha, 3)$ is a positive tail reducer.
\end{lemma}

\begin{proof}
  Suppose that $M(\alpha, 4)$ is the cycle $ab,cd,ef,gh$ such that $\alpha \in \sset{ab,cd,ef,gh}$. In view of \cref{lem:induced-cycle}, we may assume that $abcdefgh$ is an induced cycle in $G$, and moreover we may assume that $ab, ef$ are red, and $cd, gh$ are blue.
  \begin{ctikz}
    \defc{a/0/1,b/1/1,c/1/0,d/2/0,e/2/1,f/3/1,g/3/0,h/4/0}
    \drawe{a/b,e/f}{c/d,g/h}{b/c,d/e,f/g}
    \draw (a) to [bend right] (2,0.5) to [bend left] (h);
    \drawv{a,b,e,f}{}{c,d,g,h}{}
  \end{ctikz}
  \cref{lem:blue-edge} implies that the third neighbors of $a,b,e,f$ are all bluish. \ce{b} implies that $a$ and $e$ do not share a bluish neighbor because otherwise $\sset{a,d,e,h,i}$ is a cut enhancer, where $i$ is the shared bluish neighbor of $a$ and $e$. Similarly $b$ and $f$ do not share a bluish neighbor.
  \begin{kase}
    \item If $b$ and $e$ share a bluish neighbor, then $\sset{b,e}$ is a positive tail reducer by \cref{fig:ptr-c}.
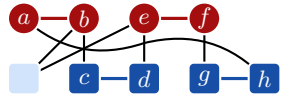
    \begin{ctikz}
      \defc{a/0/1,b/1/1,c/1/0,d/2/0,e/2/1,f/3/1,g/3/0,h/4/0,i/0/0}
      \drawe{a/b,e/f}{c/d,g/h}{b/c,d/e,f/g,e/i,b/i}
      \draw (a) to [bend right] (2,0.59) to [bend left] (h);
      \drawv{a,b,e,f}{}{c,d,g,h}{}
      \drawvnl{}{}{}{i}
    \end{ctikz}
    \item Suppose $b$ and $e$ do not share a bluish neighbor. By a symmetric argument, we may assume that neither do $a$ and $f$. Depending on how the bluish neighbors of $a,b,e,f$ overlap, there are three possibilities up to symmetry. In each case, $\sset{a,b,e,f}$ is a positive tail reducer by \cref{fig:ptr-v,fig:ptr-w,fig:ptr-y}.
    \begin{ctikz}
      \begin{scope}
        \defc{a/0/1,b/1/1,c/1/0,d/2/0,e/2/1,f/3/1,g/3/0,h/4/0,i/0/0,k/5/0}
        \drawe{a/b,e/f}{c/d,g/h}{b/c,d/e,f/g,a/i,b/i,f/k}
        \draw (a) to [bend right] (2,0.5) to [bend left] (h);
        \draw (e) .. controls (3.14,1.6) .. (k);
        \drawv{a,b,e,f}{}{c,d,g,h}{}
        \drawvnl{}{}{}{i,k}
      \end{scope}
      \begin{scope}[shift={(7.5,0)}]
        \defc{a/0/1,b/1/1,c/1/0,d/2/0,e/2/1,f/3/1,g/3/0,h/4/0,i/-1/0,j/0/0,k/5/0}
        \drawe{a/b,e/f}{c/d,g/h}{b/c,d/e,f/g,a/i,b/j,f/k}
        \draw (a) to [bend right] (2,0.5) to [bend left] (h);
        \draw (e) .. controls (3.14,1.6) .. (k);
        \drawv{a,b,e,f}{}{c,d,g,h}{}
        \drawvnl{}{}{}{i,j,k}
      \end{scope}
      \begin{scope}[shift={(3.75,-2.5)}]
        \defc{a/0/1,b/1/1,c/1/0,d/2/0,e/2/1,f/3/1,g/3/0,h/4/0,i/-1/0,j/0/0,k/6/0,l/5/0}
        \drawe{a/b,e/f}{c/d,g/h}{b/c,d/e,f/g,a/i,b/j,f/l}
        \draw (a) to [bend right] (2,0.5) to [bend left] (h);
        \draw (e) .. controls (3,1.6) .. (k);
        \drawv{a,b,e,f}{}{c,d,g,h}{}
        \drawvnl{}{}{}{i,j,k,l}
      \end{scope}
    \end{ctikz}
  \end{kase}

  Finally, note that the distance of $a$ and $ef$ is $3$, and it is the furthest distance between any vertex in a positive tail reducer and $\alpha \in \sset{ab,cd,ef,gh}$.
\end{proof}

Behold, our first foray into the depths: a 7-level case analysis. Let's dive in!

\begin{lemma} \label{lem:path-4-12}
  If $M(\alpha, 4)$ is a path on $4$ vertices, then a subset of $G(\alpha, 12)$ is a positive tail reducer.
\end{lemma}

\begin{proof}
  Suppose that $M(\alpha, 4)$ is the path $ab,cd,ef,gh$ such that $\alpha \in \sset{ab,cd,ef,gh}$. In view of \cref{lem:induced-path}, we may assume that $abcdefgh$ is an induced path in $G$, and moreover we may assume that $ab$ and $ef$ are red, and $cd$ and $gh$ are blue. Let $i$ be the third neighbor of $d$. \cref{lem:blue-edge} implies that $i$ is reddish.

  \begin{ctikz}
    \defc{a/0/1,b/1/1,c/1/0,d/2/0,e/2/1,f/3/1,g/3/0,h/4/0,i/4/1}
    \drawe{a/b,e/f}{c/d,g/h}{b/c,d/e,f/g,d/i}
    \drawv{a,b,e,f}{i}{c,d,g,h}{}
  \end{ctikz}

  \begin{kase}
    \item If $i \sim g$, then $\sset{d,g}$ is a positive tail reducer by \cref{fig:ptr-c}.
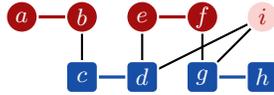
    \begin{ctikz}
      \defc{a/0/1,b/1/1,c/1/0,d/2/0,e/2/1,f/3/1,g/3/0,h/4/0,i/4/1}
      \drawe{a/b,e/f}{c/d,g/h}{b/c,d/e,f/g,d/i,i/g}
      \drawv{a,b,e,f}{i}{c,d,g,h}{}
    \end{ctikz}
    \item Suppose $i \nsim g$ and $i \nsim h$. Flip the cut preserver $\sset{d, e}$. Then the blue edge $gh$ becomes isolated in $M$, and so a subset of $\sset{g, h}$ becomes a positive tail reducer by \cref{lem:isolated}.
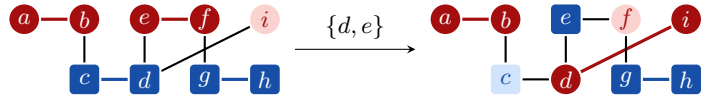
    \begin{ctikz}
      \begin{scope}
        \defc{a/0/1,b/1/1,c/1/0,d/2/0,e/2/1,f/3/1,g/3/0,h/4/0,i/4/1}
        \drawe{a/b,e/f}{c/d,g/h}{b/c,d/e,f/g,d/i}
        \drawv{a,b,e,f}{i}{c,d,g,h}{}
      \end{scope}
      \draw[-stealth, thin] (4.5,0.5) -- (6.5,0.5) node[midway, above, font=\footnotesize] {$\sset{d,e}$};
      \begin{scope}[shift={(7,0)}]
        \defc{a/0/1,b/1/1,c/1/0,d/2/0,e/2/1,f/3/1,g/3/0,h/4/0,i/4/1}
        \drawe{a/b,d/i}{g/h}{b/c,c/d,d/e,e/f,f/g}
        \drawv{a,b,d,i}{f}{e,g,h}{c}
      \end{scope}
    \end{ctikz}
    \item Suppose $i \nsim g$ and $i \sim h$. Let $j$ be the third neighbor of $e$. \cref{lem:blue-edge} implies that $j$ is bluish. By a symmetric argument, we may assume that $j \sim a$ but $j \nsim b$. Clearly, the third neighbor of $a$ is bluish. \cref{lem:blue-edge} implies that the third neighbor of $b$ is bluish.
    \begin{ctikz}
      \defc{a/0/1,b/1/1,c/1/0,d/2/0,e/2/1,f/3/1,g/3/0,h/4/0,i/4/1,j/0/0}
      \drawe{a/b,e/f}{c/d,g/h}{b/c,d/e,f/g,d/i,i/h,j/a,j/e}
      \drawv{a,b,e,f}{i}{c,d,g,h}{j}
    \end{ctikz}
    \begin{kase}
      \item Suppose $i \sim j$. Flip the cut preserver $\sset{d, e}$. Then $M$ contains multiple edges between $di$ and $ej$, and so a subset of $\sset{d, e, i, j}$ becomes a positive tail reducer by \cref{lem:no-multi}.
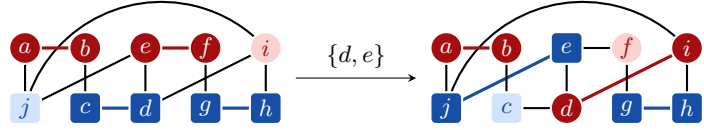
      \begin{ctikz}
        \begin{scope}
          \defc{a/0/1,b/1/1,c/1/0,d/2/0,e/2/1,f/3/1,g/3/0,h/4/0,i/4/1,j/0/0}
          \arcfourone{j}
          \drawe{a/b,e/f}{c/d,g/h}{b/c,d/e,f/g,d/i,i/h,j/a,j/e}
          \drawv{a,b,e,f}{i}{c,d,g,h}{j}
        \end{scope}
        \draw[-stealth, thin] (4.5,0.5) -- (6.5,0.5) node[midway, above, font=\footnotesize] {$\sset{d,e}$};
        \begin{scope}[shift={(7,0)}]
          \defc{a/0/1,b/1/1,c/1/0,d/2/0,e/2/1,f/3/1,g/3/0,h/4/0,i/4/1,j/0/0}
          \arcfourone{j}
          \drawe{a/b,i/d}{g/h,j/e}{b/c,c/d,d/e,e/f,f/g,a/j,i/h}
          \drawv{a,b,d,i}{f}{e,g,h,j}{c}
        \end{scope}
      \end{ctikz}
      \item If $a$ and $b$ share a bluish neighbor, then $\sset{a,b,e}$ is a positive tail reducer by \cref{fig:ptr-n}.
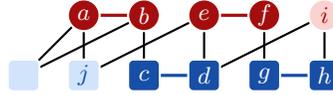
      \begin{ctikz}
        \defc{a/0/1,b/1/1,c/1/0,d/2/0,e/2/1,f/3/1,g/3/0,h/4/0,i/4/1,j/0/0,k/-1/0}
        \drawe{a/b,e/f}{c/d,g/h}{b/c,d/e,f/g,d/i,i/h,j/a,j/e,a/k,b/k}
        \drawv{a,b,e,f}{i}{c,d,g,h}{j}
        \drawvnl{}{}{}{k}
      \end{ctikz}
      \item If $a$ and $b$ share no bluish neighbors, $i \sim c$, and $j \sim f$, then $\sset{a,b,f,i}$ is a positive tail reducer by \cref{fig:ptr-x}.
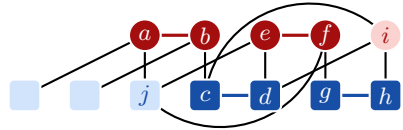
      \begin{ctikz}
        \defc{a/0/1,b/1/1,c/1/0,d/2/0,e/2/1,f/3/1,g/3/0,h/4/0,i/4/1,j/0/0,k/-1/0,l/-2/0}
        \arcthreeoneneg{j}
        \arcthreeone{c}
        \drawe{a/b,e/f}{c/d,g/h}{b/c,d/e,f/g,d/i,i/h,j/a,j/e,a/l,k/b}
        \drawv{a,b,e,f}{i}{c,d,g,h}{j}
        \drawvnl{}{}{}{k,l}
      \end{ctikz}
      \item Suppose $a$ and $b$ share no bluish neighbors, $i \nsim j$, and $i \nsim c$ or $j \nsim f$. We may assume that $i \nsim c$ up to symmetry. Let $k$ be the third neighbor of $c$. \cref{lem:blue-edge} implies that $k$ is reddish, and \ce{b} implies that $k \nsim g$ because otherwise $\sset{b,c,f,g,k}$ is a cut enhancer. We may assume that $k$ is of degree $3$ because otherwise $\sset{b,k}$ is an absolute tail reducer by \cref{lem:deg-3-flip}. 
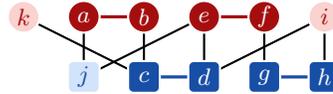
\label{4-path-a-b-no-share}
      \begin{ctikz}
        \defc{a/0/1,b/1/1,c/1/0,d/2/0,e/2/1,f/3/1,g/3/0,h/4/0,i/4/1,j/0/0,k/-1/1}
        \drawe{a/b,e/f}{c/d,g/h}{b/c,d/e,f/g,d/i,i/h,j/a,j/e,k/c}
        \drawv{a,b,e,f}{i,k}{c,d,g,h}{j}
      \end{ctikz}
      \begin{kase}
        \item Suppose $k \sim j$. Flip the cut preserver $\sset{f, g}$. Then $\sset{c, j}$ becomes a positive tail reducer by \cref{fig:ptr-m}.
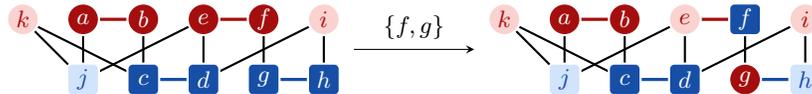
        \begin{ctikz}
          \begin{scope}
            \defc{a/0/1,b/1/1,c/1/0,d/2/0,e/2/1,f/3/1,g/3/0,h/4/0,i/4/1,j/0/0,k/-1/1}
            \drawe{a/b,e/f}{c/d,g/h}{b/c,d/e,f/g,d/i,i/h,j/a,j/e,k/c,k/j}
            \drawv{a,b,e,f}{i,k}{c,d,g,h}{j}
          \end{scope}
          \draw[-stealth, thin] (4.5,0.5) -- (6.5,0.5) node[midway, above, font=\footnotesize] {$\sset{f,g}$};
          \begin{scope}[shift={(8,0)}]
            \defc{a/0/1,b/1/1,c/1/0,d/2/0,e/2/1,f/3/1,g/3/0,h/4/0,i/4/1,j/0/0,k/-1/1}
            \drawe{a/b,e/f}{c/d,g/h}{b/c,c/k,k/j,d/e,f/g,d/i,i/h,e/j,j/a}
            \drawv{a,b,g}{k,e,i}{c,d,f}{h,j}
          \end{scope}
        \end{ctikz}
        \item Suppose $k \nsim j$ and $k \sim h$. Flip the cut preserver $\sset{b,c}$, and then flip the cut preserver $\sset{h,k}$.
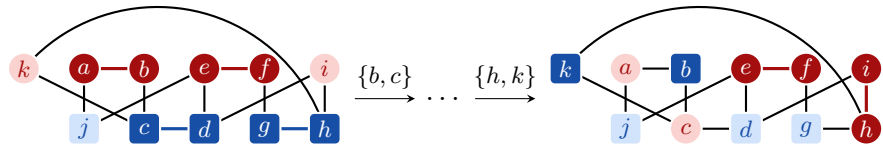
        \begin{ctikz}
          \begin{scope}
            \defc{a/0/1,b/1/1,c/1/0,d/2/0,e/2/1,f/3/1,g/3/0,h/4/0,i/4/1,j/0/0,k/-1/1}
            \arcfiveminusone{k}
            \drawe{a/b,e/f}{c/d,g/h}{b/c,d/e,f/g,d/i,i/h,j/a,j/e,k/c}
            \drawv{a,b,e,f}{i,k}{c,d,g,h}{j}
          \end{scope}
          \draw[-stealth, thin] (4.5,0.5) -- (5.5,0.5) node[midway, above, font=\footnotesize] {$\sset{b,c}$};
          \draw[-stealth, thin] (6.5,0.5) -- (7.5,0.5) node[midway, above, font=\footnotesize] {$\sset{h,k}$};
          \node at (6,0.5) {$\dots$};
          \begin{scope}[shift={(9,0)}]
            \defc{a/0/1,b/1/1,c/1/0,d/2/0,e/2/1,f/3/1,g/3/0,h/4/0,i/4/1,j/0/0,k/-1/1}
            \arcfiveminusone{k}
            \drawe{e/f,i/h}{}{a/b,b/c,c/k,c/d,d/e,f/g,g/h,d/i,e/j,j/a}
            \drawv{e,f,i,h}{a,c}{b,k}{d,g,j}
          \end{scope}
        \end{ctikz}
        \begin{kase}
          \item If $f$ has no blue neighbor, then $ef$ becomes isolated in $M$, and so a subset of $\sset{e, f}$ becomes a positive tail reducer by \cref{lem:isolated}.
          \item Suppose $f$ has a blue neighbor, say $l$. \ce{b} implies that $l \sim k$ because otherwise $\sset{f,g,h,k,l}$ is a cut enhancer. Then $\sset{f, h}$ becomes a positive tail reducer by \cref{fig:ptr-c}.
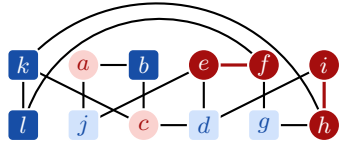
          \begin{ctikz}
            \defc{a/0/1,b/1/1,c/1/0,d/2/0,e/2/1,f/3/1,g/3/0,h/4/0,i/4/1,j/0/0,k/-1/1,l/-1/0}
            \arcfourone{l}
            \arcfiveminusone{k}
            \drawe{e/f,i/h}{}{a/b,b/c,c/k,c/d,d/e,f/g,g/h,d/i,e/j,j/a,k/l}
            \drawv{e,f,i,h}{a,c}{b,k,l}{d,g,j}
          \end{ctikz}
        \end{kase}
        \item Suppose $k \nsim j$ and $k \nsim h$. Let $l$ and $m$ be the other neighbors of $k$.
        \begin{ctikz}
          \defc{a/0/1,b/1/1,c/1/0,d/2/0,e/2/1,f/3/1,g/3/0,h/4/0,i/4/1,j/0/0,k/-1/1}
          \drawe{a/b,e/f}{c/d,g/h}{b/c,d/e,f/g,d/i,i/h,j/a,j/e,k/c}
          \drawv{a,b,e,f}{i,k}{c,d,g,h}{j}
        \end{ctikz}
        \begin{kase}
          \item Suppose $l$ and $m$ are bluish. 
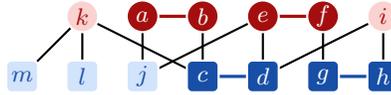
\label{3.3.3.1}
          \begin{ctikz}
            \defc{a/0/1,b/1/1,c/1/0,d/2/0,e/2/1,f/3/1,g/3/0,h/4/0,i/4/1,j/0/0,k/-1/1,l/-1/0,m/-2/0}
            \drawe{a/b,e/f}{c/d,g/h}{b/c,d/e,f/g,d/i,i/h,j/a,j/e,k/c,k/l,k/m}
            \drawv{a,b,e,f}{i,k}{c,d,g,h}{j,l,m}
          \end{ctikz}
          \begin{kase}
            \item If either $l \sim b$ or $m \sim b$, then $\sset{b, k}$ is a positive tail reducer by \cref{fig:ptr-l}.
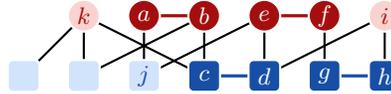
            \begin{ctikz}
              \defc{a/0/1,b/1/1,c/1/0,d/2/0,e/2/1,f/3/1,g/3/0,h/4/0,i/4/1,j/0/0,k/-1/1,l/-1/0,m/-2/0}
              \drawe{a/b,e/f}{c/d,g/h}{b/c,d/e,f/g,d/i,i/h,j/a,j/e,k/c,b/l,l/k,k/m}
              \drawv{a,b,e,f}{i,k}{c,d,g,h}{j}
              \drawvnl{}{}{}{l,m}
            \end{ctikz}
            \item If $l \nsim b$, $m \nsim b$, $l \nsim a$ and $m \nsim a$, then $\sset{a,b,k}$ is a positive tail reducer by \cref{fig:ptr-u}.
            \begin{ctikz}
              \defc{a/0/1,b/1/1,c/1/0,d/2/0,e/2/1,f/3/1,g/3/0,h/4/0,i/4/1,j/0/0,k/-1/1,l/-1/0,m/-2/0,n/-3/0,o/-4/0}
              \drawe{a/b,e/f}{c/d,g/h}{b/c,d/e,f/g,d/i,i/h,j/a,j/e,k/c,k/l,k/m}
              \draw (a) .. controls (-1.14,1.6) .. (n);
              \draw (b) .. controls (-1.3,1.78) .. (o);
              \drawv{a,b,e,f}{i,k}{c,d,g,h}{j,l,m}
              \drawvnl{}{}{}{n,o}
            \end{ctikz}
            \item Suppose $l \nsim b$, $m \nsim b$, and $l \sim a$ or $m \sim a$. We may assume that $l \sim a$ up to symmetry.
            \begin{kase}
              \item If $l$ has no other red neighbors, then $\sset{c,l}$ is a positive tail reducer by \cref{fig:degree-ptr-a} when $l$ has degree $2$, and by \cref{fig:ptr-m} when $l$ has degree $3$.
              \begin{ctikz}
                \defc{a/0/1,b/1/1,c/1/0,d/2/0,e/2/1,f/3/1,g/3/0,h/4/0,i/4/1,j/0/0,k/-1/1,l/-1/0,m/-2/0,o/-2/1}
                \drawe{a/b,e/f}{c/d,g/h}{b/c,d/e,f/g,d/i,i/h,j/a,j/e,k/c,k/l,k/m,l/a,l/o}
                \drawv{a,b,e,f}{i,k}{c,d,g,h}{j,l,m}
                \drawvnl{}{o}{}{}
              \end{ctikz}
              \item \label{case:path-4-l-sim-f} Suppose $l \sim f$. \Cref{lem:blue-edge} implies that the third neighbor of $g$ is reddish. Then $\sset{c,g,l}$ is a positive tail reducer by \cref{fig:ptr-r}.
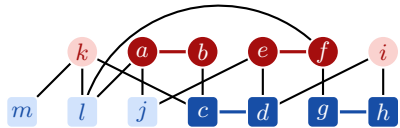
              \begin{ctikz}
                \defc{a/0/1,b/1/1,c/1/0,d/2/0,e/2/1,f/3/1,g/3/0,h/4/0,i/4/1,j/0/0,k/-1/1,l/-1/0,m/-2/0,n/5/1}
                \drawe{a/b,e/f}{c/d,g/h}{b/c,d/e,f/g,d/i,i/h,j/a,j/e,k/c,k/l,k/m,l/a,g/n}
                \arcfourone{l}
                \drawv{a,b,e,f}{i,k}{c,d,g,h}{j,l,m}
                \drawvnl{}{n}{}{}
              \end{ctikz}
              \item Suppose $l \nsim f$, and $l$ has another red neighbor, say $n$.
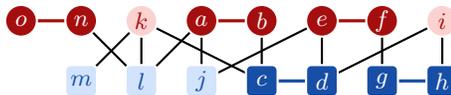
              \begin{ctikz}
                \defc{a/0/1,b/1/1,c/1/0,d/2/0,e/2/1,f/3/1,g/3/0,h/4/0,i/4/1,j/0/0,k/-1/1,l/-1/0,m/-2/0,n/-2/1,o/-3/1}
                \drawe{a/b,e/f,n/o}{c/d,g/h}{b/c,d/e,f/g,d/i,i/h,j/a,j/e,k/c,k/l,k/m,l/a,l/n}
                \drawv{a,b,e,f,n,o}{i,k}{c,d,g,h}{j,l,m}
              \end{ctikz}
              \begin{kase}
                \item Suppose $n$ has a blue neighbor. \ce{a} implies $n \nsim g$ because otherwise $\sset{f,g,n}$ is a cut enhancer. \ce{b} implies that $n \nsim h$ because otherwise $\sset{d,e,i,h,n}$ is a cut enhancer. Let $o$ be the blue neighbor of $n$. Flip the cut preserver $\sset{n,o}$. Then $M$ contains the path $gh, ef, cd, ab, ln$ as a subgraph, and so a subset of $G(cd,9)$ is a positive tail reducer by \cref{lem:path-5-subgraph-9}. \label{3.3.3.2.3.2.1}
                \begin{ctikz}
                  \begin{scope}
                    \defc{a/0/1,b/1/1,c/1/0,d/2/0,e/2/1,f/3/1,g/3/0,h/4/0,i/4/1,j/0/0,k/-1/1,l/-1/0,m/-2/0,n/-2/1,o/-3/0}
                    \drawe{a/b,e/f}{c/d,g/h}{b/c,d/e,f/g,d/i,i/h,j/a,j/e,k/c,k/l,k/m,l/a,l/n,n/o}
                    \drawv{a,b,e,f,n}{i,k}{c,d,g,h,o}{j,l,m}
                  \end{scope}
                  \begin{scope}[shift={(6,-2)}]
                    \draw[-stealth, thin] (-5.5,0.5) -- (-3.5,.5) node[midway,above,font=\footnotesize] {$\sset{n,o}$};
                    \defc{a/0/1,b/1/1,c/1/0,d/2/0,e/2/1,f/3/1,g/3/0,h/4/0,i/4/1,j/0/0,k/-1/1,l/-1/0,m/-2/0,n/-2/1,o/-3/0}
                    \drawe{a/b,e/f}{c/d,g/h,n/l}{b/c,d/e,f/g,d/i,i/h,j/a,j/e,k/c,k/l,k/m,l/a,n/o}
                    \drawv{a,b,e,f,o}{k}{c,d,g,h,n,l}{}
                    \node[a-vertex] at (i) {$i$};
                    \node[b-vertex] at (j) {$j$};
                    \node[b-vertex] at (m) {$m$};
                  \end{scope}
                \end{ctikz}
                \item Suppose $n$ has no blue neighbors. Let $o$ be the other bluish neighbor of $n$. \ce{c} implies $n \nsim j$ because otherwise $\sset{a,d,e,j,n}$ is a cut enhancer. Then $\sset{a,n}$ is a positive tail reducer by \cref{fig:ptr-g} regardless of whether $o$ is identical to $m$.
                \begin{ctikz}
                  \begin{scope}
                    \defc{a/0/1,b/1/1,c/1/0,d/2/0,e/2/1,f/3/1,g/3/0,h/4/0,i/4/1,j/0/0,k/-1/1,l/-1/0,m/-2/0,n/-2/1,o/-3/0}
                    \drawe{a/b,e/f}{c/d,g/h}{b/c,d/e,f/g,d/i,i/h,j/a,j/e,k/c,k/l,k/m,l/a,l/n,n/o}
                    \drawv{a,b,e,f,n}{i,k}{c,d,g,h}{j,l,m,o}
                  \end{scope}
                \end{ctikz}
              \end{kase}
            \end{kase}
          \end{kase}
          \item Suppose $l$ and $m$ are blue. Flip the cut preserver $\sset{b,c}$. \ce{a} implies that $l \sim m$ because otherwise $\sset{k,l,m}$ is a cut enhancer. Then $M$ contains multiple edges between $ck$ and $lm$, and so a subset of $\sset{c, k, l, m}$ becomes a positive tail reducer by \cref{lem:no-multi}.
          \begin{ctikz}
            \defc{a/0/1,b/1/1,c/1/0,d/2/0,e/2/1,f/3/1,g/3/0,h/4/0,i/4/1,j/0/0,k/-1/1,l/-1/0,m/-2/0}
            \drawe{c/k,e/f}{g/h,m/l}{e/j,j/a,a/b,b/c,c/d,d/e,d/i,i/h,f/g,m/k,k/l}
            \drawv{k,c,e,f}{a,i}{b,g,h,m,l}{j,d}
          \end{ctikz}
          \item Suppose exactly one of $l$ and $m$ is bluish, and the other is blue. We may assume that $l$ is bluish and $m$ is blue. Let $n$ be the blue neighbor of $m$. \cref{lem:blue-edge} implies that $m \nsim b$, $m \nsim f$, $n \nsim b$ and $n \nsim f$. \ce{b} implies that $m \nsim a$ and $n \nsim a$ because otherwise $\sset{a,d,e,j,m}$ or $\sset{a,d,e,j,n}$ is a cut enhancer. \ce{c} implies that $m \nsim i$ and $n \nsim i$ because otherwise $\sset{d,e,h,i,m}$ or $\sset{d,e,h,i,n}$ is a cut enhancer. Let $o$ be the third neighbor of $m$. \ce{b} implies that $o$ is reddish because otherwise $\sset{b,c,k,m,o}$ is a cut enhancer.
          \begin{ctikz}
            \defc{a/0/1,b/1/1,c/1/0,d/2/0,e/2/1,f/3/1,g/3/0,h/4/0,i/4/1,j/0/0,k/-1/1,l/-1/0,m/-2/0,n/-3/0,o/-2/1}
            \drawe{a/b,e/f}{c/d,g/h,m/n}{b/c,d/e,f/g,d/i,i/h,j/a,j/e,k/c,k/l,k/m,o/m}
            \drawv{a,b,e,f}{i,k,o}{c,d,g,h,m,n}{j,l}
          \end{ctikz}
          \begin{kase}
            \item If $n$ has no red neighbor, then $mn$ is isolated in $M$, and so a subset of $\sset{m, n}$ is a positive tail reducer by \cref{lem:isolated}.
            \item Suppose that $n$ has a red neighbor, say $p$. \ce{a} implies that $p \nsim g$ because otherwise $\sset{g,f,p}$ is a cut enhancer. \ce{b} implies that $p \nsim h$ and $p \nsim j$ because otherwise $\sset{d,e,h,i,p}$ or $\sset{d,e,j,n,p}$ is a cut enhancer.
            \begin{ctikz}
              \defc{a/0/1,b/1/1,c/1/0,d/2/0,e/2/1,f/3/1,g/3/0,h/4/0,i/4/1,j/0/0,k/-1/1,l/-1/0,m/-2/0,n/-3/0,o/-2/1,p/-3/1}
              \drawe{a/b,e/f}{c/d,g/h,m/n}{b/c,d/e,f/g,d/i,i/h,j/a,j/e,k/c,k/l,k/m,o/m,p/n}
              \drawv{a,b,e,f,p}{i,k,o}{c,d,g,h,m,n}{j,l}
            \end{ctikz}
            \begin{kase}
              \item If $p \sim l$, then $\sset{e,k,p}$ is a positive tail reducer by \cref{fig:ptr-s}.
              \begin{ctikz}
                \defc{a/0/1,b/1/1,c/1/0,d/2/0,e/2/1,f/3/1,g/3/0,h/4/0,i/4/1,j/0/0,k/-1/1,l/-1/0,m/-2/0,n/-3/0,o/-2/1,p/-3/1}
                \drawe{a/b,e/f}{c/d,g/h,m/n}{b/c,d/e,f/g,d/i,i/h,j/a,j/e,k/c,k/l,k/m,o/m,p/n,p/l}
                \drawv{a,b,e,f,p}{i,k,o}{c,d,g,h,m,n}{j,l}
              \end{ctikz}
              \item Suppose $p \nsim l$. Flip the cut preserver $\sset{n,p}$. Let $q$ be the third neighbor of $a$. \ce{b} implies that $q$ is bluish because otherwise $\sset{a,e,d,j,q}$ is a cut enhancer. \Cref{lem:blue-edge} implies that the third neighbor of $b$ is bluish. Therefore, the current case reduces to \ref{3.3.3.1}.
              \begin{ctikz}
                \node at (0,0.5) {$\dots$};
                \draw[-stealth, thin] (0.5,0.5) -- (2.5,0.5) node[midway, above, font=\footnotesize] {$\sset{n,p}$};
                \begin{scope}[shift={(7,0)}]
                  \defc{a/0/1,b/1/1,c/1/0,d/2/0,e/2/1,f/3/1,g/3/0,h/4/0,i/4/1,j/0/0,k/-1/1,l/-1/0,m/-2/0,n/-3/0,o/-2/1,p/-3/1,q/-4/0}
                  \drawe{a/b,e/f}{c/d,g/h}{b/c,d/e,f/g,d/i,i/h,j/a,j/e,k/c,k/l,k/m,o/m,m/n,p/n}
                  \arcfourone{q}
                  \drawv{a,b,e,f,n}{i,k,o}{c,d,g,h,p}{j,l,m,q}
                \end{scope}
              \end{ctikz}
            \end{kase}
          \end{kase}
        \end{kase}
      \end{kase}
    \end{kase}
  \end{kase}

  Finally, notice that the distance between a vertex in $G(cd,9)$ and $gh$ is at most $12$ in \ref{3.3.3.2.3.2.1}, and it is the furthest distance possible between any vertex in a cut preserver or a positive tail reducer and $\alpha \in \sset{ab, cd, ef, gh}$.
\end{proof}

\begin{lemma} \label{lem:path-4-subgraph-12}
  If $M$ contains a path $\alpha\beta\gamma\delta$ as a subgraph, then a subset of $G(\beta, 12)$ is a positive tail reducer.
\end{lemma}

\begin{proof}
  \begin{kase}
    \item If $M(\beta,4)$ contains multiple edges, then a subset of $G(\beta,8)$ is a positive tail reducer by \cref{lem:no-multi}.
    \item If $M(\beta,4)$ is a simple graph of maximum degree at least $3$, then a subset of $G(\beta,9)$ is a positive tail reducer by \cref{lem:m-deg-2}.
    \item If $M(\beta,4)$ is a simple graph of maximum degree at most $2$, then $M(\beta,4)$ must be a path on at least $4$ vertices, or an even cycle (on at least $4$ vertices), and so a subset of $G(\beta, 12)$ is a positive tail reducer by \cref{lem:path-5-8,lem:cycle-4,lem:path-4-12}. \qedhere
  \end{kase}
\end{proof}

Lo and behold, we now descend into yet another 7-level case analysis. May your coffee be strong and your patience stronger.

\begin{lemma} \label{lem:path-3-14}
  If $M(\alpha, 4)$ is a path on $3$ vertices, then a subset of $G(\alpha, 14)$ is a positive tail reducer.
\end{lemma}

\begin{proof}
  Suppose that $M(\alpha, 4)$ is the path $ab,cd,ef$ such that $\alpha \in \sset{ab,cd,ef}$. In view of \cref{lem:induced-path}, we may assume that $abcdef$ is an induced path in $G$, and moreover we may assume that $ab,ef$ are red, and $cd$ is blue. Let $g$ be the third neighbor of $b$. \cref{lem:blue-edge} implies that $g$ is bluish.
  \begin{ctikz}
    \defc{a/0/1,b/1/1,c/1/0,d/2/0,e/2/1,f/3/1,g/3/0}
    \drawe{a/b,e/f}{c/d}{b/c,d/e,b/g}
    \drawv{a,b,e,f}{}{c,d}{g}
  \end{ctikz}
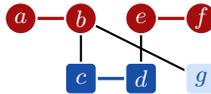
  \begin{kase}
    \item If $g \sim e$, then $\sset{b,e}$ is a positive tail reducer by \cref{fig:ptr-c}.
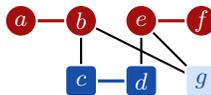
    \begin{ctikz}
      \defc{a/0/1,b/1/1,c/1/0,d/2/0,e/2/1,f/3/1,g/3/0}
      \drawe{a/b,e/f}{c/d}{b/c,d/e,b/g,g/e}
      \drawv{a,b,e,f}{}{c,d}{g}
    \end{ctikz}
    \item Suppose $g \nsim e$ and $g \nsim f$. Flip the cut preserver $\sset{b,c}$. Then the red edge $ef$ becomes isolated in $M$, and so a subset of $\sset{e,f}$ becomes a positive tail reducer by \cref{lem:isolated}.
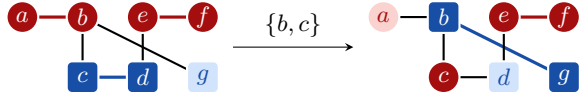
    \begin{ctikz}
      \begin{scope}
        \defc{a/0/1,b/1/1,c/1/0,d/2/0,e/2/1,f/3/1,g/3/0}
        \drawe{a/b,e/f}{c/d}{b/c,d/e,b/g}
        \drawv{a,b,e,f}{}{c,d}{g}
      \end{scope}
      \draw[-stealth, thin] (3.5,0.5) -- (5.5,0.5) node[midway, above, font=\footnotesize] {$\sset{b,c}$};
      \begin{scope}[shift={(6,0)}]
        \defc{a/0/1,b/1/1,c/1/0,d/2/0,e/2/1,f/3/1,g/3/0}
        \drawe{e/f}{b/g}{a/b,c/d,b/c,d/e}
        \drawv{c,e,f}{a}{b,g}{d}
      \end{scope}
    \end{ctikz}
    \item Suppose $g \nsim e$ and $g \sim f$. Let $h$ be the third neighbor of $e$. \cref{lem:blue-edge} implies that $h$ is bluish. By a symmetric argument, we may assume that $h \nsim b$ and $h \sim a$.
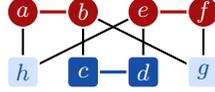
    \begin{ctikz}
      \defc{a/0/1,b/1/1,c/1/0,d/2/0,e/2/1,f/3/1,g/3/0,h/0/0}
      \drawe{a/b,e/f}{c/d}{b/c,d/e,b/g,a/h,e/h,f/g}
      \drawv{a,b,e,f}{}{c,d}{g,h}
    \end{ctikz}
    \begin{kase}
      \item If $g \sim a$ or $h \sim f$, then $\sset{a,b,e}$ or $\sset{b,e,f}$ is a positive tail reducer by \cref{fig:ptr-n}.
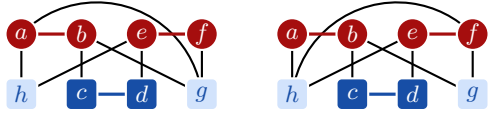
      \begin{ctikz}
        \begin{scope}
          \defc{a/0/1,b/1/1,c/1/0,d/2/0,e/2/1,f/3/1,g/3/0,h/0/0}
          \arcthreeminusone{a}
          \drawe{a/b,e/f}{c/d}{b/c,d/e,b/g,a/h,e/h,f/g}
          \drawv{a,b,e,f}{}{c,d}{g,h}
        \end{scope}
        \begin{scope}[shift={(4.5,0)}]
          \defc{a/0/1,b/1/1,c/1/0,d/2/0,e/2/1,f/3/1,g/3/0,h/0/0}
          \arcthreeone{h}
          \drawe{a/b,e/f}{c/d}{b/c,d/e,b/g,a/h,e/h,f/g}
          \drawv{a,b,e,f}{}{c,d}{g,h}
        \end{scope}
      \end{ctikz}
      \item Suppose $g \nsim a$ and $h \nsim f$. Let $i$ be the third neighbor of $c$. \cref{lem:blue-edge} implies that $i$ is reddish. We may assume that $i$ is of degree $3$ because otherwise $\sset{b,i}$ is an absolute tail reducer by \cref{lem:deg-3-flip}. \ce{d} implies that $i \nsim g$ because otherwise $\sset{b,c,f,g,i}$ is a cut enhancer.
      \begin{kase}
        \item Suppose some neighbor, say $j$, of $i$ outside $\sset{c,d}$ is blue. Let $k$ be the blue neighbor of $j$. Flip the cut preserver $\sset{b,c}$. Then $M$ contains the path $ef,bg,ci,jk$ as a subgraph, and so a subset of $G(bg,12)$ becomes a positive tail reducer by \cref{lem:path-4-subgraph-12}. 
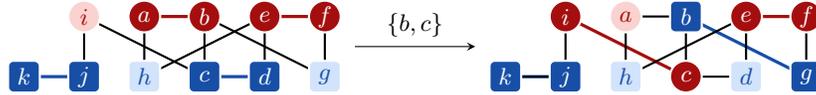
\label{3.2.1}
        \begin{ctikz}
          \begin{scope}
            \defc{a/0/1,b/1/1,c/1/0,d/2/0,e/2/1,f/3/1,g/3/0,h/0/0,i/-1/1,j/-1/0,k/-2/0}
            \drawe{a/b,e/f}{c/d,j/k}{b/c,d/e,b/g,a/h,e/h,f/g,i/c,i/j}
            \drawv{a,b,e,f}{i}{c,d,j,k}{g,h}
          \end{scope}
          \draw[-stealth, thin] (3.5,0.5) -- (5.5,0.5) node[midway, above, font=\footnotesize] {$\sset{b,c}$};
          \begin{scope}[shift={(8,0)}]
            \defc{a/0/1,b/1/1,c/1/0,d/2/0,e/2/1,f/3/1,g/3/0,h/0/0,i/-1/1,j/-1/0,k/-2/0}
            \drawe{i/c,e/f}{b/g,j/k}{a/b,c/d,b/c,d/e,a/h,e/h,f/g,j/k,i/j}
            \drawv{c,e,f,i}{a}{b,g,j,k}{d,h}
          \end{scope}
        \end{ctikz}
        \item Suppose every neighbor of $i$ outside $\sset{c,d}$ is bluish.
        \begin{kase}
          \item Suppose $i \sim d$. \ce{d} implies that $i \nsim h$ because otherwise $\sset{a,d,e,h,i}$ is a cut enhancer. Let $j$ be the third neighbor of $i$. Then $\sset{b,e,i}$ is a positive tail reducer by \cref{fig:ptr-p}.
          \begin{ctikz}
            \defc{a/0/1,b/1/1,c/1/0,d/2/0,e/2/1,f/3/1,g/3/0,h/0/0,i/-1/1,j/-1/0}
            \arcthreeminusoneneg{i}
            \drawe{a/b,e/f}{c/d}{b/c,d/e,b/g,a/h,e/h,f/g,i/c,i/j}
            \drawv{a,b,e,f}{i}{c,d}{g,h,j}
          \end{ctikz}
          \item Suppose $i \nsim d$ and $i \sim h$. Let $j$ be the third neighbor of $i$. Then $\sset{e,i}$ is a positive tail reducer by \cref{fig:ptr-m}.
          \begin{ctikz}
            \defc{a/0/1,b/1/1,c/1/0,d/2/0,e/2/1,f/3/1,g/3/0,h/0/0,i/-1/1,j/-1/0}
            \drawe{a/b,e/f}{c/d}{b/c,d/e,b/g,a/h,e/h,f/g,i/c,i/h,i/j}
            \drawv{a,b,e,f}{i}{c,d}{g,h,j}
          \end{ctikz}
          \item Suppose $i \nsim d$ and $i \nsim h$. Let $j$ and $k$ be the bluish neighbors of $i$.
          \begin{kase}
            \item Suppose $j \nsim a$ and $k \nsim a$. Let $l$ be the third neighbor of $a$. \ce{b} implies that $l$ is bluish because otherwise $\sset{a,d,e,h,l}$ is a cut enhancer. Then $\sset{a,b,i}$ is a positive tail reducer by \cref{fig:ptr-u}.
            \begin{ctikz}
              \defc{a/0/1,b/1/1,c/1/0,d/2/0,e/2/1,f/3/1,g/3/0,h/0/0,i/-1/1,j/-1/0,k/-2/0,l/-3/0}
              \draw (a) .. controls (-1.14,1.6) .. (l);
              \drawe{a/b,e/f}{c/d}{b/c,d/e,b/g,a/h,e/h,f/g,i/c,i/k,i/j}
              \drawv{a,b,e,f}{i}{c,d}{g,h,j,k,l}
            \end{ctikz}
            \item Suppose $j \sim a$ or $k \sim a$. We may assume that $j \sim a$ up to symmetry.
            \begin{kase}
              \item If $j$ has no other red neighbor, $\sset{c,j}$ is a positive tail reducer by \cref{fig:degree-ptr-a} when $j$ has degree $2$, and by \cref{fig:ptr-m} when $j$ has degree $3$.
              \label{3.2.4.2.1}
              \begin{ctikz}
                \defc{a/0/1,b/1/1,c/1/0,d/2/0,e/2/1,f/3/1,g/3/0,h/0/0,i/-1/1,j/-1/0,k/-2/0,l/-2/1}
                \drawe{a/b,e/f}{c/d}{b/c,d/e,b/g,a/h,e/h,f/g,i/c,i/j,i/k,j/a,j/l}
                \drawv{a,b,e,f}{i}{c,d}{g,h,k,j}
                \drawvnl{}{l}{}{}
              \end{ctikz}
              \item If $j \sim f$, then $\sset{a,f}$ is a positive tail reducer by \cref{fig:ptr-g}.
              \begin{ctikz}
                \defc{a/0/1,b/1/1,c/1/0,d/2/0,e/2/1,f/3/1,g/3/0,h/0/0,i/-1/1,j/-1/0,k/-2/0}
                \drawe{a/b,e/f}{c/d}{b/c,d/e,b/g,a/h,e/h,f/g,i/c,i/j,i/k,j/a}
                \arcfourone{j}
                \drawv{a,b,e,f}{i}{c,d}{g,h,k,j}
              \end{ctikz}
              \item Suppose $j$ has another red neighbor, and $j \nsim f$. Let $l$ be the other red neighbor of $j$, and let $m$ be the red neighbor of $l$. \ce{a} implies that $l \nsim d$ and $m \nsim d$ because otherwise $\sset{d,e,l}$ or $\sset{d,e,m}$ is a cut enhancer.
              \begin{ctikz}
                \defc{a/0/1,b/1/1,c/1/0,d/2/0,e/2/1,f/3/1,g/3/0,h/0/0,i/-1/1,j/-1/0,k/-2/0,l/-2/1,m/-3/1}
                \drawe{a/b,e/f,l/m}{c/d}{b/c,d/e,b/g,a/h,e/h,f/g,i/c,i/j,i/k,j/a,j/l}
                \drawv{a,b,e,f,l,m}{i}{c,d}{g,h,k,j}
              \end{ctikz}
              \begin{kase}
                \item If neither $l$ nor $m$ has blue neighbors, then $lm$ is isolated in $M$, and so a subset of $\sset{l,m}$ is a positive tail reducer by \cref{lem:isolated}.
                \item Suppose $l$ has a blue neighbor. Let $n$ be a blue neighbor of $l$. Flip the cut preserver $\sset{l,n}$. Then $M$ contains the path $ef, cd, ab, jl$, and so a subset of $G(\sset{cd},12)$ is a positive tail reducer by \cref{lem:path-4-subgraph-12}. \label{3.2.2.3.2.3.2}
                \begin{ctikz}
                  \begin{scope}
                    \defc{a/0/1,b/1/1,c/1/0,d/2/0,e/2/1,f/3/1,g/3/0,h/0/0,i/-1/1,j/-1/0,l/-2/1,m/-3/1,n/-3/0,k/-2/0}
                    \drawe{a/b,e/f,l/m}{c/d}{b/c,d/e,b/g,a/h,e/h,f/g,i/c,i/j,i/k,j/a,j/l,l/n}
                    \drawv{a,b,e,f,l,m}{i}{c,d,n}{g,h,k,j}
                  \end{scope}
                  \draw[-stealth, thin] (3.5,0.5) -- (4.5,.5) node[midway, above, font=\footnotesize] {$\sset{l,n}$};
                  \begin{scope}[shift={(8,0)}]
                    \defc{a/0/1,b/1/1,c/1/0,d/2/0,e/2/1,f/3/1,g/3/0,h/0/0,i/-1/1,j/-1/0,l/-2/1,m/-3/1,n/-3/0,k/-2/0}
                    \drawe{a/b,e/f}{c/d,j/l}{b/c,d/e,b/g,a/h,e/h,f/g,i/c,i/j,i/k,j/a,l/m,l/n}
                    \drawv{a,b,e,f,n}{i}{c,d,l,j}{g,h,k}
                    \node[a-vertex] at (m) {$m$};
                  \end{scope}
                \end{ctikz}
              \item Suppose $l$ has no blue neighbors, and $m$ has a blue neighbor. Let $n$ be a blue neighbor of $m$. \ce{b} implies that $m \nsim h$ and $m \nsim g$ because otherwise $\sset{d,e,h,m,n}$ or $\sset{b,c,g,m,n}$ is a cut enhancer. Flip the cut preserver $\sset{m,n}$. Then, depending on whether $m \sim k$, the current case reduces to \ref{3.2.4.2.1}, where $j$ has no red neighbors other than $a$, or \ref{3.2.1}, where some neighbor of $i$ outside $\sset{c,d}$ is blue.
              
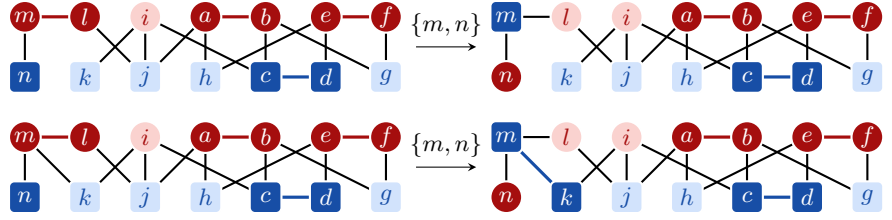
\begin{ctikz}
                \begin{scope}
                  \defc{a/0/1,b/1/1,c/1/0,d/2/0,e/2/1,f/3/1,g/3/0,h/0/0,i/-1/1,j/-1/0,l/-2/1,m/-3/1,n/-3/0,k/-2/0}
                  \drawe{a/b,e/f,l/m}{c/d}{b/c,d/e,b/g,a/h,e/h,f/g,i/c,i/j,i/k,j/a,j/l,m/n}
                  \drawv{a,b,e,f,l,m}{i}{c,d,n}{g,h,k,j}
                  \draw[-stealth, thin] (3.5,0.5) -- (4.5,0.5) node[midway, above, font=\footnotesize] {$\sset{m,n}$};
                \end{scope}
                \begin{scope}[shift={(8,0)}]
                  \defc{a/0/1,b/1/1,c/1/0,d/2/0,e/2/1,f/3/1,g/3/0,h/0/0,i/-1/1,j/-1/0,l/-2/1,m/-3/1,n/-3/0,k/-2/0}
                  \drawe{a/b,e/f}{c/d}{b/c,d/e,b/g,a/h,e/h,f/g,i/c,i/j,i/k,j/a,j/l,m/n,l/m}
                  \drawv{a,b,e,f,n}{i,l}{c,d,m}{g,h,k,j}
                \end{scope}
                \begin{scope}[shift={(0,-2)}]
                  \defc{a/0/1,b/1/1,c/1/0,d/2/0,e/2/1,f/3/1,g/3/0,h/0/0,i/-1/1,j/-1/0,l/-2/1,m/-3/1,n/-3/0,k/-2/0}
                  \drawe{a/b,e/f,l/m}{c/d}{b/c,d/e,b/g,a/h,e/h,f/g,i/c,i/j,i/k,j/a,j/l,m/n,m/k}
                  \drawv{a,b,e,f,l,m}{i}{c,d,n}{g,h,k,j}
                  \draw[-stealth, thin] (3.5,0.5) -- (4.5,0.5) node[midway, above, font=\footnotesize] {$\sset{m,n}$};
                \end{scope}
                \begin{scope}[shift={(8,-2)}]
                  \defc{a/0/1,b/1/1,c/1/0,d/2/0,e/2/1,f/3/1,g/3/0,h/0/0,i/-1/1,j/-1/0,l/-2/1,m/-3/1,n/-3/0,k/-2/0}
                  \drawe{a/b,e/f}{c/d,m/k}{b/c,d/e,b/g,a/h,e/h,f/g,i/c,i/j,i/k,j/a,j/l,m/n,l/m}
                  \drawv{a,b,e,f,n}{i,l}{c,d,m,k}{g,h,j}
                \end{scope}
              \end{ctikz}
              \end{kase}
            \end{kase}
          \end{kase}
        \end{kase}
      \end{kase}
    \end{kase}
  \end{kase}

  Finally notice that the distance between a vertex of $G(cd,12)$ and $ef$ is at most $14$ in \ref{3.2.2.3.2.3.2}, and it is the furthest distance possible between any vertex in a cut preserver or a positive tail reducer and $\alpha \in \sset{ab,cd,ef}$.
\end{proof}

\begin{lemma} \label{lem:path-3-subgraph-14}
  If $M$ contains a path $\alpha\beta\gamma$ as a subgraph, then a subset of $G(\beta,14)$ is a positive tail reducer.
\end{lemma}

\begin{proof}
  \begin{kase}
    \item If $M(\beta, 4)$ contains multiple edges, then a subset of $G(\beta, 8)$ is a positive tail reducer by \cref{lem:no-multi}.
    \item If $M(\beta, 4)$ is a simple graph of maximum degree at least $3$, then a subset of $G(\beta, 9)$ is a positive tail reducer by \cref{lem:m-deg-2}.
    \item If $M(\beta, 4)$ is a simple graph of maximum degree at most $2$, then $M(\beta, 4)$ must be a path on at least $3$ vertices, or an even cycle (on at least $4$ vertices), and so a subset of $G(\beta,14)$ is a positive tail reducer by \cref{lem:path-5-8,lem:cycle-4,lem:path-4-12,lem:path-3-14}. \qedhere
  \end{kase}
\end{proof}

\begin{lemma} \label{lem:path-2-16}
  If $M(\alpha, 4)$ is a path on $2$ vertices, then a subset of $G(\alpha, 16)$ is a positive tail reducer.
\end{lemma}

\begin{proof}
  Suppose that $M(\alpha, 4)$ is the path $ab, cd$ such that $\alpha \in \sset{ab,cd}$. In view of \cref{lem:induced-path}, we may assume that $abcd$ is an induced path in $G$, and moreover we may assume that $ab$ is red, and $cd$ is blue. Let $e$ and $f$ be the third neighbors of $b$ and $c$ respectively. \cref{lem:blue-edge} implies that $e$ is bluish and $f$ is reddish. We may assume that both $e$ and $f$ are of degree $3$ because otherwise $\sset{c,e}$ or $\sset{b,f}$ is an absolute tail reducer by \cref{lem:deg-3-flip}.
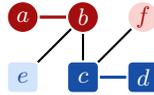
  \begin{ctikz}
    \defc{a/0/1,b/1/1,c/1/0,d/2/0,e/0/0,f/2/1}
    \drawe{a/b}{c/d}{b/c,b/e,c/f}
    \drawv{a,b}{f}{c,d}{e}
  \end{ctikz}
  \begin{kase}
    \item Suppose $e \sim f$. Flip the cut preserver $\sset{b,c}$. Then $M$ contains multiple edges between $be$ and $cf$, and so a subset of $\sset{b,c,e,f}$ becomes a positive tail reducer by \cref{lem:no-multi}.
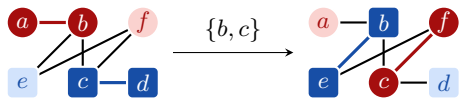
    \begin{ctikz}
      \begin{scope}
        \defc{a/0/1,b/1/1,c/1/0,d/2/0,e/0/0,f/2/1}
        \drawe{a/b}{c/d}{b/c,b/e,c/f,e/f}
        \drawv{a,b}{f}{c,d}{e}
      \end{scope}
      \draw[-stealth, thin] (2.5,0.5) -- (4.5,0.5) node[midway, above, font=\footnotesize] {$\sset{b,c}$};
      \begin{scope}[shift={(5,0)}]
        \defc{a/0/1,b/1/1,c/1/0,d/2/0,e/0/0,f/2/1}
        \drawe{c/f}{b/e}{a/b,c/d,b/c,e/f}
        \drawv{c,f}{a}{b,e}{d}
      \end{scope}
    \end{ctikz}
    \item Suppose $e \nsim f$, $e \nsim a$, and $e$ has another red neighbor. Let $g$ be a red neighbor of $e$, and let $h$ be the red neighbor of $g$. Flip the cut preserver $\sset{b,c}$. Then $M$ contains the path $cf, be, hg$ on $3$ as a subgraph, and so a subset of $G(be, 14)$ is a positive tail reducer by \cref{lem:path-3-subgraph-14}.
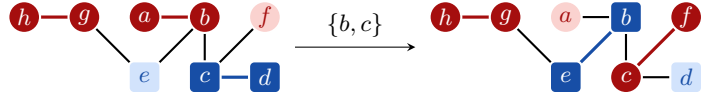
    \begin{ctikz}
      \begin{scope}
        \defc{a/0/1,b/1/1,c/1/0,d/2/0,e/0/0,f/2/1,g/-1/1,h/-2/1}
        \drawe{a/b,g/h}{c/d}{b/c,b/e,c/f,e/g}
        \drawv{a,b,g,h}{f}{c,d}{e}
      \end{scope}
      \draw[-stealth, thin] (2.5,0.5) -- (4.5,0.5) node[midway, above, font=\footnotesize] {$\sset{b,c}$};
      \begin{scope}[shift={(7,0)}]
        \defc{a/0/1,b/1/1,c/1/0,d/2/0,e/0/0,f/2/1,g/-1/1,h/-2/1}
        \drawe{c/f,g/h}{b/e}{b/c,a/b,c/d,e/g}
        \drawv{c,g,h,f}{a}{b,e}{d}
      \end{scope}
    \end{ctikz}
    \item Suppose $e \nsim f$, $e \nsim a$, and $e$ has no other red neighbors. Clearly, both $d$ and $e$ have two reddish neighbors. \label{3-recall}
    \begin{kase}
      \item Suppose $d$ and $e$ share a reddish neighbor. Let $g$ be a shared reddish neighbor of $d$ and $e$.

      \begin{kase}
        \item Suppose $g$ is not adjacent to any blue vertex other than $d$. Then either $g$ is of degree 2 or has a third neighbor $h$ which is bluish. In each case, $\sset{b,g}$ is a positive tail reducer by \cref{fig:degree-ptr-a} or \cref{fig:ptr-m}, respectively. 
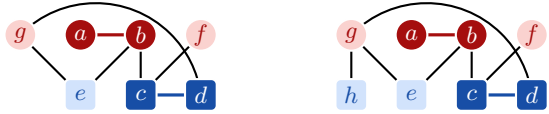
\label{3.1.1}
        \begin{ctikz}
        \begin{scope}
            \defc{a/0/1,b/1/1,c/1/0,d/2/0,e/0/0,f/2/1,g/-1/1}
        \arcthreeminusone{g}
        \drawe{a/b}{c/d}{b/c,b/e,c/f,e/g}
        \drawv{a,b}{f,g}{c,d}{e}
        \end{scope}
          \begin{scope}[shift={(5.5,0)}]
              \defc{a/0/1,b/1/1,c/1/0,d/2/0,e/0/0,f/2/1,g/-1/1,h/-1/0}
          \arcthreeminusone{g}
          \drawe{a/b}{c/d}{b/c,b/e,c/f,e/g,g/h}
          \drawv{a,b}{f,g}{c,d}{e,h}
          \end{scope}
        \end{ctikz}
        \item Suppose $g$ is adjacent to a blue vertex, say $h$, other than $d$. Let $i$ be the blue neighbor of $h$. Clearly, $h \nsim a$ and $i \nsim a$.
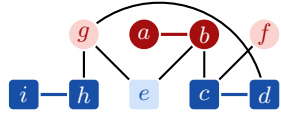
        \begin{ctikz}
          \defc{a/0/1,b/1/1,c/1/0,d/2/0,e/0/0,f/2/1,g/-1/1,h/-1/0,i/-2/0}
          \arcthreeminusone{g}
          \drawe{a/b}{c/d,h/i}{b/c,b/e,c/f,e/g,g/h}
          \drawv{a,b}{f,g}{c,d,h,i}{e}
        \end{ctikz}
        \begin{kase}
          \item If neither $h$ nor $i$ has red neighbors, then $hi$ is isolated in $M$, and so a subset of $\sset{h, i}$ is a positive tail reducer by \cref{lem:isolated}.
          \item Suppose $h$ has a red neighbor. Let $j$ be the red neighbor of $h$. Flip the cut preserver $\sset{h,j}$. Then $M$ contains the path $ab,cd,gh$ as a subgraph, and so a subset of $G(cd, 14)$ is a positive tail reducer by \cref{lem:path-3-subgraph-14}.
          \begin{ctikz}
            \begin{scope}
              \defc{a/0/1,b/1/1,c/1/0,d/2/0,e/0/0,f/2/1,g/-1/1,h/-1/0,i/-2/0,j/-2/1}
              \arcthreeminusone{g}
              \drawe{a/b}{c/d,h/i}{b/c,b/e,c/f,e/g,g/h,h/j}
              \drawv{a,b,j}{f,g}{c,d,h,i}{e}
            \end{scope}
            \draw[-stealth, thin] (2.5,0.5) -- (4.5,0.5) node[midway, above, font=\footnotesize] {$\sset{h,j}$};
            \begin{scope}[shift={(7,0)}]
              \defc{a/0/1,b/1/1,c/1/0,d/2/0,e/0/0,f/2/1,g/-1/1,h/-1/0,i/-2/0,j/-2/1}
              \arcthreeminusone{g}
              \drawe{a/b}{c/d,h/i}{b/c,b/e,c/f,e/g,g/h,h/j}
              \drawv{a,b,h,g}{f}{c,d,j}{}
              \node[b-vertex] at (e) {$e$};
              \node[b-vertex] at (i) {$i$};
            \end{scope}
          \end{ctikz}
          \item Suppose $h$ has no red neighbor, and $i$ has a red neighbor. Let $j$ be a red neighbor of $i$. Clearly $j \nsim h$. Recall from \ref{3-recall} that $b$ is the only red neighbor of $e$, and so $j \nsim e$. \ce{b} implies that $i \nsim f$ because otherwise $\sset{b,c,f,i,j}$ is a cut enhancer. Flip the cut preserver $\sset{i,j}$. Then the current case reduces to \ref{3.1.1}, where $h$ is bluish.
          
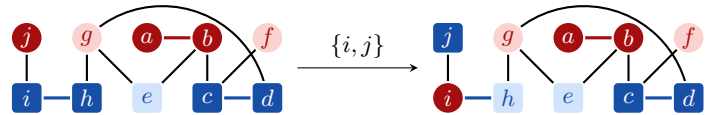
\begin{ctikz}
            \begin{scope}
              \defc{a/0/1,b/1/1,c/1/0,d/2/0,e/0/0,f/2/1,g/-1/1,h/-1/0,i/-2/0,j/-2/1}
              \arcthreeminusone{g}
              \drawe{a/b}{c/d,h/i}{b/c,b/e,c/f,e/g,g/h,i/j}
              \drawv{a,b,j}{f,g}{c,d,h,i}{e}
            \end{scope}
            \draw[-stealth, thin] (2.5,0.5) -- (4.5,0.5) node[midway, above, font=\footnotesize] {$\sset{i,j}$};
            \begin{scope}[shift={(7,0)}]
              \defc{a/0/1,b/1/1,c/1/0,d/2/0,e/0/0,f/2/1,g/-1/1,h/-1/0,i/-2/0,j/-2/1}
              \arcthreeminusone{g}
              \drawe{a/b}{c/d,h/i}{b/c,b/e,c/f,e/g,g/h,i/j}
              \drawv{a,b,i}{f,g}{c,d,j}{e,h}
            \end{scope}
          \end{ctikz}
        \end{kase}
      \end{kase}
      \item Suppose $d$ and $e$ share no reddish neighbors. Depending on whether $f \sim d$, there are two possibilities. In each case, $\sset{c,d,e}$ is a positive tail reducer by \cref{fig:ptr-t,fig:ptr-u}.
      \begin{ctikz}
        \begin{scope}
          \defc{a/0/1,b/1/1,c/1/0,d/2/0,e/0/0,f/2/1,g/-1/1,h/-2/1,i/3/1,j/4/1}
          \drawe{a/b}{c/d}{b/c,b/e,c/f,e/g,j/d,d/i,e/h}
          \drawv{a,b}{f,g,h,i,j}{c,d}{e}
        \end{scope}
        \begin{scope}[shift={(7.5,0)}]
          \defc{a/0/1,b/1/1,c/1/0,d/2/0,e/0/0,f/2/1,g/-1/1,h/-2/1,i/3/1}
          \drawe{a/b}{c/d}{b/c,b/e,c/f,e/g,f/d,d/i,e/h}
          \drawv{a,b}{f,g,h,i}{c,d}{e}
        \end{scope}
      \end{ctikz}
    \end{kase}
    \item Suppose $e \nsim f$ and $e \sim a$. By a symmetric argument, we may assume that $f \sim d$. Let $g$ be the third neighbor of $e$.
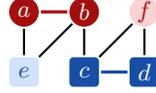
    \begin{ctikz}
      \defc{a/0/1,b/1/1,c/1/0,d/2/0,e/0/0,f/2/1}
      \drawe{a/b}{c/d}{b/c,b/e,c/f,a/e,d/f}
      \drawv{a,b}{f}{c,d}{e}
    \end{ctikz}
    \begin{kase}
      \item Suppose $g$ is red. Let $h$ be the red neighbor of $g$. Flip the cut preserver $\sset{b,c}$. Then $M$ contains the path $cf, be, hg$ as a subgraph, and so a subset of $G(be, 14)$ is a positive tail reducer by \cref{lem:path-3-subgraph-14}. 
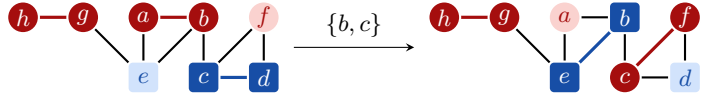
\label{4.1}
      \begin{ctikz}
        \begin{scope}
          \defc{a/0/1,b/1/1,c/1/0,d/2/0,e/0/0,f/2/1,g/-1/1,h/-2/1}
          \drawe{a/b,g/h}{c/d}{b/c,b/e,c/f,a/e,d/f,e/g}
          \drawv{a,b,g,h}{f}{c,d}{e}
        \end{scope}
        \draw[-stealth, thin] (2.5,0.5) -- (4.5,0.5) node[midway, above, font=\footnotesize] {$\sset{b,c}$};
        \begin{scope}[shift={(7,0)}]
          \defc{a/0/1,b/1/1,c/1/0,d/2/0,e/0/0,f/2/1,g/-1/1,h/-2/1}
          \drawe{c/f,g/h}{b/e}{b/c,a/b,c/d,e/g,a/e,f/d}
          \drawv{c,g,h,f}{a}{b,e}{d}
        \end{scope}
      \end{ctikz}
      \item Suppose $g$ is reddish. Depending on whether $g \sim d$, there are two possibilities. In each case, $\sset{c,d,e}$ is a positive tail reducer by \cref{fig:ptr-o,fig:ptr-q}.
      \begin{ctikz}
        \begin{scope}
          \defc{a/0/1,b/1/1,c/1/0,d/2/0,e/0/0,f/2/1,g/-1/1}
          \arcthreeminusoneneg{g}
          \drawe{a/b}{c/d}{b/c,b/e,c/f,a/e,d/f,e/g}
          \drawv{a,b}{f,g}{c,d}{e}
        \end{scope}
        \begin{scope}[shift={(4.5,0)}]
          \defc{a/0/1,b/1/1,c/1/0,d/2/0,e/0/0,f/2/1,g/-1/1,h/3/1}
          \drawe{a/b}{c/d}{b/c,b/e,c/f,a/e,d/f,e/g,h/d}
          \drawv{a,b}{f,g,h}{c,d}{e}
        \end{scope}
      \end{ctikz}
    \end{kase}
  \end{kase}

  Finally notice that the distance between a vertex of $G(be,14)$ and $cd$ is at most $16$ in \ref{4.1}, and it is the furthest distance possible between any vertex in a cut preserver or a positive tail reducer and $\alpha \in \sset{ab,cd}$.
\end{proof}

\begin{proof}[Proof of \cref{lem:pos-tail-16}]
  Let $\alpha$ be a vertex of $M$.
  \begin{kase}
    \item If $M(\alpha,4)$ contains multiple edges, then a subset of $G(\alpha,8)$ is a positive tail reducer by \cref{lem:no-multi}.
    \item If $M(\alpha,4)$ is a simple graph of maximum degree at least $3$, then a subset of $G(\alpha,9)$ is a positive tail reducer by \cref{lem:m-deg-2}.
    \item If $M(\alpha,4)$ is a simple graph of maximum degree at most $2$, then $M(\alpha,4)$ must be an isolated vertex, a path (on at least $2$ vertices), or an even cycle (on at least $4$ vertices), and so a subset of $G(\alpha,16)$ is a positive tail reducer by \cref{lem:isolated,lem:path-5-8,lem:cycle-4,lem:path-4-12,lem:path-3-14,lem:path-2-16}. \qedhere
  \end{kase}
\end{proof}

\begin{figure}
  \centering
  \begin{csubfig}{2.0cm}{ntr-a}
    \defc{a/0/1,b/0/0,c/1/0}
    \drawe{}{c/b}{a/b,a/c}
    \drawv{a}{}{b,c}{}
  \end{csubfig}
  \begin{csubfig}{2.0cm}{ntr-b}
    \defc{a/0/1,b/1/1,c/0/0,d/1/0}
    \drawe{}{}{a/c,a/d,b/c,b/d}
    \drawv{a,b}{}{}{c,d}
  \end{csubfig}
  \begin{csubfig}{2.0cm}{ntr-c}
    \defc{a/0/1,b/1/1,c/0/0,d/1/0}
    \drawe{a/b}{}{a/c,a/d,b/c,b/d}
    \drawv{a,b}{}{}{c,d}
  \end{csubfig}
  \begin{csubfig}{2.8cm}{ntr-d}
    \defc{a/0/1,b/1/1,c/-1/0,d/0/0,e/1/0}
    \drawe{}{}{a/c,a/d,b/d,b/e}
    \drawv{a,b}{}{}{c,d,e}
  \end{csubfig}
  \begin{csubfig}{2.8cm}{ntr-e}
    \defc{a/0/1,b/1/1,c/0/0,d/1/0,e/2/0}
    \drawe{}{}{a/c,a/d,b/c,b/d,b/e}
    \drawv{a}{b}{}{c,d,e}
  \end{csubfig}
  \begin{csubfig}{3.6cm}{ntr-f}
    \defc{a/0/1,b/1/1,d/0/0,e/1/0,f/2/0,c/-1/0}
    \drawe{a/b}{d/e}{a/c,a/d,b/f,b/e}
    \drawv{a,b}{}{e,d}{c,f}
  \end{csubfig}

  \par\bigskip

  \begin{csubfig}{3.6cm}{ntr-g}
    \defc{a/0/1,b/1/1,c/-1/0,d/0/0,e/1/0,f/2/0}
    \drawe{a/b}{}{a/c,a/d,b/e,b/f}
    \drawv{a,b}{}{}{c,d,e,f}
  \end{csubfig}
  \begin{csubfig}{3.6cm}{ntr-h}
    \defc{a/0/1,b/1/1,d/0/0,e/1/0,f/2/0,c/-1/0}
    \drawe{}{e/f}{a/c,a/d,b/d,b/e,b/f}
    \drawv{a}{b}{e,f}{c,d}
  \end{csubfig}
  \begin{csubfig}{3.6cm}{ntr-i}
    \defc{a/0/1,b/1/1,c/-1/0,d/0/0,e/1/0,f/2/0}
    \drawe{}{c/d}{a/d,a/e,b/d,b/e,b/f}
    \drawv{a}{b}{c,d}{e,f}
  \end{csubfig}
  \begin{csubfig}{4.4cm}{ntr-j}
    \defc{a/0/1,b/1/1,e/0/0,f/1/0,g/2/0,d/-1/0,c/-2/0}
    \drawe{}{d/e,f/g}{a/c,a/e,b/e,b/f,b/g}
    \drawv{a}{b}{d,e,f,g}{c}
  \end{csubfig}

  \par\bigskip

  \begin{csubfig}{2.8cm}{ntr-k}
    \defc{a/0/1,b/1/1,c/2/1,d/0/0,e/1/0,f/2/0}
    \drawe{a/b}{d/e}{a/d,a/f,b/f,b/e,c/d,c/e,c/f}
    \drawv{a,b}{c}{e,d}{f}
  \end{csubfig}
  \begin{csubfig}{3.6cm}{ntr-l}
    \defc{a/0/1,b/1/1,c/2/1,d/-1/0,e/0/0,f/1/0,g/2/0}
    \drawe{a/b}{f/g}{a/d,a/e,b/d,b/f,c/e,c/g}
    \drawv{a,b,c}{}{f,g}{d,e}
  \end{csubfig}
  \begin{csubfig}{3.6cm}{ntr-m}
    \defc{a/0/1,b/1/1,c/2/1,d/0/0,e/1/0,f/2/0,g/3/0}
    \drawe{a/b}{d/e}{a/d,a/f,b/f,b/e,c/f,c/g}
    \drawv{a,b,c}{}{e,d}{f,g}
  \end{csubfig}
  \begin{csubfig}{3.6cm}{ntr-n}
    \defc{a/1/1,b/2/1,c/3/1,d/0/0,e/1/0,f/2/0,g/3/0}
    \drawe{a/b}{}{a/f,a/d,b/f,b/e,c/f,c/g}
    \drawv{a,b,c}{}{}{d,e,f,g}
  \end{csubfig}

  \par\medskip

  \begin{csubfig}{3.6cm}{ntr-o}
    \defc{a/0/1,b/1/1,c/2/1,d/0/0,e/1/0,f/2/0,g/3/0}
    \drawe{a/b}{d/e}{a/d,a/f,b/f,b/e,c/d,c/e,c/g}
    \drawv{a,b}{c}{e,d}{f,g}
  \end{csubfig}
  \begin{csubfig}{3.6cm}{ntr-p}
    \defc{a/0/1,b/1/1,c/2/1,d/0/0,e/1/0,f/2/0,g/3/0}
    \drawe{}{e/f}{a/d,a/e,a/f,b/d,b/e,c/f,c/g}
    \drawv{b,c}{a}{e,f}{d,g}
  \end{csubfig}
  \begin{csubfig}{4.4cm}{ntr-q}
    \defc{a/0/1,b/1/1,c/2/1,e/0/0,f/1/0,g/2/0,h/3/0,d/-1/0}
    \drawe{}{f/g}{a/d,a/e,a/f,b/e,b/g,c/h,c/f}
    \drawv{b,c}{a}{g,f}{d,e,h}
  \end{csubfig}
  \begin{csubfig}{4.4cm}{ntr-r}
    \defc{a/0/1,b/1/1,c/2/1,e/0/0,f/1/0,g/2/0,h/3/0,d/-1/0}
    \drawe{}{e/f}{a/d,a/e,a/f,b/e,b/g,c/h,c/f}
    \drawv{b,c}{a}{e,f}{d,g,h}
  \end{csubfig}

  \par\medskip

  \begin{csubfig}{4.4cm}{ntr-s}
    \defc{a/0/1,b/1/1,c/2/1,e/0/0,f/1/0,g/2/0,d/-1/0,h/3/0}
    \arcthreeone{d}
    \drawe{a/b}{f/g}{a/e,b/e,b/f,c/f,h/c,a/d}
    \drawv{a,b}{c}{f,g}{d,e,h}
  \end{csubfig}
  \begin{csubfig}{4.4cm}{ntr-t}
    \defc{a/1/1,b/2/1,c/3/1,d/0/0,f/2/0,e/1/0,g/3/0,h/4/0}
    \drawe{a/b}{}{a/f,a/d,b/f,b/e,c/f,c/g,c/h}
    \drawv{a,b}{c}{}{d,e,f,g,h}
  \end{csubfig}
  \begin{csubfig}{4.4cm}{ntr-u}
    \defc{a/0/1,b/1/1,f/1/0,g/2/0,e/0/0,c/2/1,d/-1/0,h/3/0}
    \arcthreeone{d}
    \arcthreeminusone{a}
    \drawe{a/b}{f/g}{b/f,b/e,f/c,a/d,c/h}
    \drawv{a,b}{c}{d,f,g}{d,e,h}
  \end{csubfig}

  \par\medskip

  \begin{csubfig}{4.4cm}{ntr-v}
    \defc{a/0/1,b/1/1,c/2/1,d/-1/0,e/0/0,f/1/0,g/2/0,h/3/0}
    \drawe{}{f/g}{a/d,a/e,a/f,b/e,b/f,b/g,c/g,c/h}
    \drawv{c}{a,b}{f,g}{d,e,h}
  \end{csubfig}
  \begin{csubfig}{4.4cm}{ntr-w}
    \defc{a/0/1,b/1/1,c/2/1,d/-1/0,e/0/0,f/1/0,g/2/0,h/3/0}
    \drawe{}{g/h}{a/d,a/e,a/f,b/d,b/e,b/g,c/f,c/h}
    \drawv{c}{a,b}{g,h}{d,e,f}
  \end{csubfig}
  \begin{csubfig}{4.4cm}{ntr-x}
    \defc{a/0/1,b/1/1,c/2/1,d/-1/0,e/0/0,f/1/0,g/2/0,h/3/0}
    \drawe{}{f/g}{a/d,a/f,b/e,b/f,b/g,c/g,c/h}
    \drawv{a,c}{b}{f,g}{d,e,h}
  \end{csubfig}

  \par\bigskip

  \begin{csubfig}{5.2cm}{ntr-y}
    \defc{a/0/1,b/1/1,c/2/1,f/0/0,g/1/0,h/2/0,i/3/0,e/-1/0,d/-2/0}
    \drawe{a/b}{e/f,g/h}{a/d,a/f,b/g,c/f,c/g,c/i}
    \draw (b) .. controls (-0.14,1.6) .. (d);
    \drawv{a,b}{c}{e,f,g,h}{d,i}
  \end{csubfig}
  \begin{csubfig}{5.2cm}{ntr-z}
    \defc{a/0/1,b/1/1,c/2/1,e/0/0,f/1/0,g/2/0,h/3/0,i/4/0,d/-1/0}
    \drawe{}{g/h}{a/d,a/e,b/e,b/f,b/g,c/g,c/i}
    \drawv{a,c}{b}{g,h}{d,e,f,i}
  \end{csubfig}
  \begin{csubfig}{5.2cm}{ntr-aa}
    \defc{a/0/1,b/1/1,c/2/1,e/0/0,f/1/0,g/2/0,h/3/0,i/4/0,d/-1/0}
    \drawe{b/c}{h/g}{a/d,a/e,a/g,b/f,b/e,c/g,c/i}
    \drawv{b,c}{a}{h,g}{d,e,f,i}
  \end{csubfig}
  \caption{Negative tail reducers.} \label{fig:ntr}
\end{figure}

\begin{figure}
  \ContinuedFloat
  \centering
  \begin{csubfig}{5.2cm}{ntr-ab}
    \defc{a/0/1,b/1/1,c/2/1,e/0/0,f/1/0,g/2/0,h/3/0,i/4/0,d/-1/0}
    \drawe{}{f/g,h/i}{a/d,a/e,b/e,b/f,c/g,c/h,c/i}
    \drawv{a,b}{c}{f,g,h,i}{d,e}
  \end{csubfig}
  \begin{csubfig}{5.2cm}{ntr-ac}
    \defc{a/0/1,b/1/1,c/2/1,d/-2/0,e/-1/0,f/0/0,g/1/0,h/2/0,i/3/0}
    \drawe{}{g/h}{a/d,a/e,a/f,b/e,b/f,b/g,c/h,c/i}
    \drawv{c}{a,b}{g,h}{d,e,f,i}
  \end{csubfig}
  \begin{csubfig}{5.2cm}{ntr-ad}
    \defc{a/0/1,b/1/1,c/2/1,d/-2/0,e/-1/0,f/0/0,g/1/0,h/2/0,i/3/0}
    \drawe{}{h/i}{a/d,a/f,a/g,b/e,b/f,b/h,c/g,c/i}
    \drawv{c}{a,b}{h,i}{d,e,f,g}
  \end{csubfig}

  \par

  \begin{csubfig}{5.2cm}{ntr-ae}
    \defc{a/0/1,b/1/1,c/2/1,d/-1/0,e/0/0,f/1/0,g/2/0,h/3/0,i/4/0}
    \drawe{}{g/h}{a/d,a/e,a/f,b/f,b/g,b/h,c/e,c/i}
    \drawv{c}{a,b}{g,h}{d,e,f,i}
  \end{csubfig}
  \begin{csubfig}{6.0cm}{ntr-af}
    \defc{a/0/1,b/1/1,c/2/1,g/0/0,h/1/0,i/2/0,j/3/0,f/-1/0,e/-2/0,d/-3/0}
    \drawe{a/b}{f/g,h/i}{a/e,a/g,b/h,c/h,c/g,c/j}
    \draw (b) .. controls (0,1.6) .. (d);
    \drawv{a,b}{c}{f,g,h,i}{d,j,e}
  \end{csubfig}

  \par

  \begin{csubfig}{6.0cm}{ntr-ag}
    \defc{a/0/1,b/1/1, c/2/1,f/0/0,g/1/0,h/2/0,i/3/0,j/4/0,e/-1/0,d/-2/0}
    \drawe{b/c}{h/i}{a/d,a/e,a/h,b/f,b/g,c/h,c/j}
    \drawv{b,c}{a}{h,i}{d,e,f,g,j}
  \end{csubfig}
  \begin{csubfig}{5.2cm}{ntr-ah}
    \defc{g/0/0,h/1/0,i/2/0,b/0/1,c/1/1,d/2/1,a/-1/1,f/-1/0,j/3/0,e/-2/0}
    \arcfourminusone{a}
    \drawe{c/d}{g/h}{b/g,b/h,h/c,c/i,d/i,a/g,b/f,a/e,d/j}
    \drawv{c,d}{a,b}{g,h}{e,f,i,j}
  \end{csubfig}
  \begin{csubfig}{5.2cm}{ntr-ai}
    \defc{g/0/0,h/1/0,i/2/0,b/0/1,c/1/1,d/2/1,a/-1/1,f/-1/0,e/-2/0,j/3/0}
    \arcfourone{e}
    \arcfourminusone{a}
    \drawe{c/d}{g/h}{b/g,b/h,c/h,c/i,a/g,b/f,a/e,d/j}
    \drawv{c,d}{a,b}{g,h}{e,f,i,j}
  \end{csubfig}

  \par

  \begin{csubfig}{6.0cm}{ntr-aj}
    \defc{a/-2/1,b/-1/1,c/0/1,d/1/1,e/-3/0,f/-2/0,g/-1/0,h/0/0,i/1/0,j/2/0,k/3/0}
    \arcfourminusone{b}
    \drawe{c/d}{e/f,j/k}{a/e,a/f,a/g,b/g,b/h,c/h,c/i,d/i,d/j}
    \drawv{c,d}{a,b}{e,f,j,k}{g,h,i}
  \end{csubfig}
  \begin{csubfig}{6.0cm}{ntr-ak}
    \defc{a/0/1,b/1/1,c/2/1,d/3/1,f/0/0,g/1/0,h/2/0,i/3/0,j/4/0,k/5/0,e/-1/0}
    \drawe{a/b,c/d}{i/j,g/h}{a/e,b/f,b/g,c/h,d/i,d/k}
    \draw (a) to [bend right] (2,0.5) to [bend left] (j);
    \draw (c) .. controls (3.14,1.6) .. (k);
  \drawv{a,b,c,d}{}{g,h,i,j}{e,f,k}
  \end{csubfig}

  \par

  \begin{csubfig}{6.8cm}{ntr-al}
    \defc{a/0/1,b/1/1,c/2/1,d/3/1,f/0/0,g/1/0,h/2/0,i/3/0,j/4/0,k/5/0,l/6/0,e/-1/0}
    \drawe{a/b,c/d}{i/j,g/h}{a/e,b/f,b/g,c/h,d/i,d/k}
    \draw (a) to [bend right] (2,0.5) to [bend left] (j);
    \draw (c) .. controls (3,1.6) .. (l);
    \drawv{a,b,c,d}{}{g,h,i,j}{e,f,k,l}
  \end{csubfig}
  \begin{csubfig}{6.8cm}{ntr-am}
    \defc{a/0/1,b/1/1,c/2/1,d/3/1,e/-2/0,f/-1/0,g/0/0,h/1/0,i/2/0,j/3/0,k/4/0,l/5/0,m/6/0}
    \drawe{c/d}{h/i}{a/e,a/f,a/h,b/g,b/h,b/i,c/i,c/j,d/k,d/l}
  \drawv{c,d}{a,b}{h,i}{e,f,g,j,k,l}
  \end{csubfig}

  \par

  \begin{csubfig}{6.0cm}{ntr-an}
    \defc{a/0/1,b/1/1,c/2/1,d/3/1,e/4/1,f/0/0,g/1/0,h/2/0,i/3/0,j/4/0,k/5/0,l/6/0}
    \drawe{b/c,d/e}{h/i,j/k}{a/f,a/h,b/g,c/h,c/g,d/i,e/j,e/l}
    \arcthreeminusone{a}
    \draw (b) to [bend right] (3,0.5) to [bend left] (k);
    \arcthreeminusone{d}
    \drawv{b,c,d,e}{a}{h,i,j,k}{f,g,l}
  \end{csubfig}
  \begin{csubfig}{7.6cm}{ntr-ao}
    \defc{a/-1/1,b/0/1,c/1/1,d/2/1,e/3/1,f/-3/0,g/-2/0,h/-1/0,i/0/0,j/1/0,k/2/0,l/3/0,m/4/0,n/5/0}
    \drawe{c/d}{i/j}{a/f,a/g,a/i,b/h,b/i,b/j,c/j,c/k,d/k,d/l,e/k,e/m,e/n}
  \drawv{c,d}{a,b,e}{i,j}{f,g,h,k,l,m,n}
  \end{csubfig}

  \par

  \begin{csubfig}{7.6cm}{ntr-ap}
    \defc{a/0/1,b/1/1,c/2/1,d/3/1,e/4/1,f/-3/0,g/-2/0,h/-1/0,i/0/0,j/1/0,k/2/0,l/3/0,m/4/0,n/5/0}
    \arcfourone{i}
    \drawe{d/e}{k/l}{a/f,a/g,a/i,b/h,b/i,b/k,c/j,c/k,c/l,d/l,d/m,e/n}
    \drawv{d,e}{a,b,c}{k,l}{f,g,h,i,j,m,n}
  \end{csubfig}
  \caption{Negative tail reducers (continued).}
\end{figure}

\begin{figure}
  \centering
  \begin{csubfig}{2.0cm}{degree-ntr-a}
    \defc{a/0/1,b/0/0,c/1/0}
    \drawe{}{b/c}{a/b,a/c}
    \node[reddish-degree-two-vertex] at (a) {$a$};
    \drawv{}{}{b,c}{}
  \end{csubfig}%
  \begin{csubfig}{2.0cm}{degree-ntr-e}
    \defc{a/0/1,b/1/1,c/0/0,d/1/0}
    \drawe{}{}{a/c,a/d,b/c,b/d}
    \node[reddish-degree-two-vertex] at (b) {$b$};
    \drawv{a}{}{}{c,d}
  \end{csubfig}%
  \begin{csubfig}{3.6cm}{degree-ntr-l}
    \defc{a/0/1,b/1/1,c/2/1,d/-1/0,e/0/0,f/1/0,g/2/0}
    \drawe{a/b}{f/g}{a/d,a/e,b/d,b/f,c/e,c/g}
    \node[reddish-degree-two-vertex] at (c) {$c$};
    \drawv{a,b}{}{f,g}{d,e}
  \end{csubfig}%
  \begin{csubfig}{3.6cm}{degree-ntr-t-two}
    \defc{a/1/1,b/2/1,c/3/1,d/0/0,e/1/0,f/2/0,g/3/0}
    \drawe{a/b}{}{a/f,a/d,b/f,b/e,c/f,c/g}
    \node[reddish-degree-two-vertex] at (c) {$c$};
    \drawv{a,b}{}{}{d,e,f,g}
  \end{csubfig}

  \par\medskip

  \begin{csubfig}{3.6cm}{degree-ntr-t-one}
    \defc{a/1/1,b/2/1,c/3/1,d/0/0,e/1/0,f/2/0}
    \drawe{a/b}{}{a/f,a/d,b/f,b/e,c/f}
    \node[reddish-degree-one-vertex] at (c) {$c$};
    \drawv{a,b}{}{}{d,e,f}
  \end{csubfig}%
  \begin{csubfig}{3.6cm}{degree-ntr-v}
    \defc{a/0/1,b/1/1,c/2/1,d/0/0,e/1/0,f/2/0,g/3/0}
    \drawe{}{e/f}{a/d,a/e,b/d,b/e,b/f,c/f,c/g}
    \node[reddish-degree-two-vertex] at (a) {$a$};
    \drawv{c}{b}{e,f}{d,g}
  \end{csubfig}%
  \begin{csubfig}{4.4cm}{degree-ntr-x-two}
    \defc{a/0/1,b/1/1,c/2/1,d/-1/0,e/0/0,f/1/0,g/2/0,h/3/0}
    \drawe{}{f/g}{a/d,a/f,b/e,b/f,b/g,c/g,c/h}
    \node[reddish-degree-two-vertex] at (a) {$a$};
    \drawv{c}{b}{f,g}{d,e,h}
  \end{csubfig}%
  \begin{csubfig}{3.6cm}{degree-ntr-x-one}
    \defc{a/0/1,b/1/1,c/2/1,d/0/0,e/1/0,f/2/0,g/3/0}
    \drawe{}{e/f}{a/e,b/d,b/e,b/f,c/f,c/g}
    \node[reddish-degree-one-vertex] at (a) {$a$};
    \drawv{c}{b}{e,f}{d,g}
  \end{csubfig}

  \par\medskip

  \begin{csubfig}{4.4cm}{degree-ntr-ac}
    \defc{a/0/1,b/1/1,c/2/1,d/-1/0,e/0/0,f/1/0,g/2/0,h/3/0}
    \drawe{}{f/g}{a/d,a/e,b/d,b/e,b/f,c/g,c/h}
    \node[reddish-degree-two-vertex] at (a) {$a$};
    \drawv{c}{b}{f,g}{d,e,h}
  \end{csubfig}%
  \begin{csubfig}{4.4cm}{degree-ntr-ad}
    \defc{a/0/1,b/1/1,c/2/1,d/-1/0,e/0/0,f/1/0,g/2/0,h/3/0}
    \drawe{}{g/h}{a/e,a/f,b/d,b/e,b/g,c/f,c/h}
    \node[reddish-degree-two-vertex] at (a) {$a$};
    \drawv{c}{b}{g,h}{d,e,f}
  \end{csubfig}%
  \begin{csubfig}{4.4cm}{degree-ntr-ae}
    \defc{a/0/1,b/1/1,c/2/1,d/-1/0,e/0/0,f/1/0,g/2/0,h/3/0}
    \drawe{}{f/g}{a/d,a/e,b/e,b/f,b/g,c/d,c/h}
    \node[reddish-degree-two-vertex] at (a) {$a$};
    \drawv{c}{b}{f,g}{d,e,h}
  \end{csubfig}

  \par\medskip

  \begin{csubfig}{6.0cm}{degree-ntr-ap-one}
    \defc{a/0/1,b/1/1,c/2/1,d/3/1,e/4/1,f/-1/0,g/0/0,h/1/0,i/2/0,j/3/0,k/4/0,l/5/0}
    \arcfourone{g}
    \drawe{d/e}{i/j}{a/g,b/f,b/g,b/i,c/h,c/i,c/j,d/j,d/k,e/l}
    \node[reddish-degree-one-vertex] at (a) {$a$};
    \drawv{d,e}{b,c}{i,j}{f,g,h,k,l}
  \end{csubfig}%
  \begin{csubfig}{6.8cm}{degree-ntr-ap-two}
    \defc{a/0/1,b/1/1,c/2/1,d/3/1,e/4/1,f/-2/0,g/-1/0,h/0/0,i/1/0,j/2/0,k/3/0,l/4/0,m/5/0}
    \arcfourone{h}
    \drawe{d/e}{j/k}{a/f,a/h,b/g,b/h,b/j,c/i,c/j,c/k,d/k,d/l,e/m}
    \node[reddish-degree-two-vertex] at (a) {$a$};
    \drawv{d,e}{b,c}{j,k}{f,g,h,i,l,m}
  \end{csubfig}

  \caption{Negative tail reducers with degree constraints.} \label{fig:degree-ntr}
\end{figure}

\section{Negative tail reducers} \label{sec:neg}

We collect all the negative tail reducers ever needed in the following proposition.

\begin{proposition} \label{prop:ntr}
  For every maximum cut $(A, B)$ of a subcubic graph $G$, if either one of the graphs in \cref{fig:ntr} is an induced subgraph $F$ of $G$, or one of the graphs in \cref{fig:degree-ntr} is an induced subgraph $F$ of $G$ satisfying the indicated degree constraints, then $A \cap V(F)$ is a negative tail reducer.
\end{proposition}

\begin{proof}
  Let $D = A \cap V(F)$. Suppose that $F$ is one of the graphs in \cref{fig:ntr}. Each reddish vertex has degree $3$ in $F$, each red vertex has two neighbors in $B \cap V(F)$, and each blue vertex has its blue neighbor in $F$. By subcubicity and \cref{lem:max-cut}, no edge of $G[D \cup B]$ joins $V(F)$ to its complement. Thus $F$ is the connected component of $G[D \cup B]$ containing $D$. The inequality $t_F^- < \abs{D}$ is verified in \cref{sec:app}.

  Finally suppose that $F$ is one of the graphs in \cref{fig:degree-ntr}. The displayed colors and degree constraints imply, by the same argument, that $F$ is the connected component of $G[D \cup B]$ containing $D$. Ignoring colors and degree constraints, $F$ is obtained from the corresponding graph $F'$ in \cref{fig:ntr} by deleting vertices, possibly none, in $B$. By the computation in \cref{sec:app} and the Cauchy interlacing theorem, $t_F^- \le t_{F'}^- < \abs{A \cap V(F')} = \abs{D}$.
\end{proof}

Just like in \cref{sec:pos}, throughout the rest of this section, we assume that $(A, B)$ is a maximum cut of a subcubic graph $G$, and $M$ is the underlying multigraph of $G$ with respect to $(A, B)$. When the conclusion of a result in this section is of the form ``a subset of $U$ is a negative tail reducer'' for a specific vertex subset $U$ of $G$, in view of \cref{lem:deg-3}, we always assume from the start of the proof that every red or blue vertex in $U$ is of degree $3$.

We establish a series of lemmas in parallel to those in \cref{sec:pos}. In their proofs, we often obtain subsets $C$ and $D$ of $U$, for a specific vertex subset $U$ of $G$, such that $C$ is a cut preserver, and $D$ is a negative tail reducer with respect to $(A \oplus C, B \oplus C)$. We always automatically apply \cref{lem:mod-cut-preserver} to obtain a subset of $U$ that is a negative tail reducer with respect to $(A, B)$.

\begin{lemma} \label{lem:no-multi-neg-2}
  If $M$ contains multiple edges between $\alpha$ and $\beta$, then a subset of $G(\alpha,2)$ is a negative tail reducer.
\end{lemma}

\begin{proof}
  Let $\alpha = ab$ and $\beta = cd$. In view of \cref{lem:multigraph-bipartite}, we may assume that $ab$ is red, and $cd$ is blue up to symmetry.
  \begin{kase}
    \item Suppose that the edges between $\sset{a,b}$ and $\sset{c,d}$ do not form a matching in $G$. We may assume that $a \sim c$ and $a \sim d$ up to symmetry. Then $\sset{a}$ is a negative tail reducer by \cref{fig:ntr-a}.
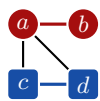
    \begin{ctikz}
      \defc{a/0/1,b/1/1,c/0/0,d/1/0}
      \drawe{a/b}{c/d}{a/c,a/d}
      \drawv{a,b}{}{c,d}{}
    \end{ctikz}
    \item Suppose that the edges between $\sset{a,b}$ and $\sset{c,d}$ form a matching in $G$. We may assume that $a \sim c$ and $b \sim d$ up to symmetry. Let $e$ be the third neighbor of $b$. \cref{lem:blue-edge} implies that $e$ is bluish.
    \begin{kase}
      \item Suppose $e \nsim a$. Let $f$ be the third neighbor of $a$. \cref{lem:blue-edge} implies that $f$ is bluish. Then $\sset{a,b}$ is a negative tail reducer by \cref{fig:ntr-f}.
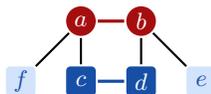
      \begin{ctikz}
        \defc{a/0/1,b/1/1,c/0/0,d/1/0,e/2/0,f/-1/0}
        \drawe{a/b}{c/d}{a/c,b/d,a/f,b/e}
        \drawv{a,b}{}{c,d}{e,f}
      \end{ctikz}
      \item Suppose $e \sim a$. Let $f$ be the third neighbor of $d$. \cref{lem:blue-edge} implies that $f$ is reddish. By a symmetric argument, we may assume that $f \sim c$. We may assume that $e$ is of degree $3$ because otherwise $\sset{c,e}$ is an absolute tail reducer by \cref{lem:deg-3-flip}.
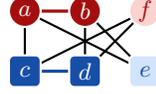
      \begin{ctikz}
        \defc{a/0/1,b/1/1,c/0/0,d/1/0,e/2/0,f/2/1}
        \drawe{a/b}{c/d}{a/c,b/d,c/f,d/f,b/e,a/e}
        \drawv{a,b}{f}{c,d}{e}
      \end{ctikz}
      \begin{kase}
        \item If $e \sim f$, then $\sset{a,b,f}$ is a negative tail reducer by \cref{fig:ntr-k}.
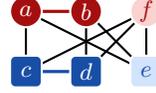
        \begin{ctikz}
          \defc{a/0/1,b/1/1,c/0/0,d/1/0,e/2/0,f/2/1}
          \drawe{a/b}{c/d}{a/c,b/d,c/f,d/f,b/e,a/e,e/f}
          \drawv{a,b}{f}{c,d}{e}
        \end{ctikz}
        \item Suppose $e \nsim f$. Let $g$ be the third neighbor of $e$.
        \begin{kase}
          \item If $g$ is reddish, then $\sset{c,d,e}$ is a negative tail reducer by \cref{fig:ntr-o}.
          \begin{ctikz}
            \defc{a/0/1,b/1/1,c/0/0,d/1/0,e/2/0,f/2/1,g/3/1}
            \drawe{a/b}{c/d}{a/c,b/d,c/f,d/f,b/e,a/e,e/g}
            \drawv{a,b}{f,g}{c,d}{e}
          \end{ctikz}
          \item Suppose $g$ is red. Let $h$ be the second neighbor of $g$ in $B$. \ce{b} implies that $h$ is bluish because otherwise $\sset{b,d,e,g,h}$ is a cut enhancer. Then $\sset{a,b,g}$ is a negative tail reducer by \cref{fig:ntr-m}.
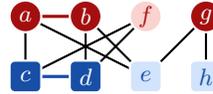
          \begin{ctikz}
            \defc{a/0/1,b/1/1,c/0/0,d/1/0,e/2/0,f/2/1,g/3/1,h/3/0}
            \drawe{a/b}{c/d}{a/c,b/d,c/f,d/f,b/e,a/e,e/g,g/h}
            \drawv{a,b,g}{f}{c,d}{e,h}
          \end{ctikz}
        \end{kase}
      \end{kase}
    \end{kase}
  \end{kase}

  Finally, notice that the distance between $g$ and $ab$ is $2$, and it is the furthest distance from a vertex in a negative tail reducer to $ab$.
\end{proof}

\begin{lemma} \label{lem:m-deg-2-neg-2}
  If $\alpha$ is a vertex of $M$, then $\alpha$ is of degree at most $2$ in $M$, or a subset of $G(\alpha,2)$ is a negative tail reducer.
\end{lemma}

\begin{proof}
  Without loss of generality, we may assume that $\alpha = ab$ is red. \cref{lem:blue-edge} implies that neither $a$ nor $b$ can be adjacent to two blue edges in $G$. If $\alpha$ is of degree at least $3$ in $M$, then $M$ contains multiple edges between $\alpha$ and its neighbor $\beta$, and so a subset of $G(\alpha, 2)$ is a negative tail reducer by \cref{lem:no-multi-neg-2}. Finally, notice that the distance between a vertex of $G(\alpha, 2)$ and $\alpha$ is at most $2$.
\end{proof}

\begin{lemma} \label{lem:isolated-neg-4}
  If $\alpha$ is an isolated vertex of $M$, then a subset of $G(\alpha, 4)$ is a negative tail reducer.
\end{lemma}

\begin{proof}
  Without loss of generality, we may assume that $\alpha = ab$ is red. Since $ab$ is isolated in $M$, the other neighbors of $a$ and $b$ are bluish.
  \begin{kase}
    \item If $a$ and $b$ share zero or two bluish neighbors, then, in either case, $\sset{a,b}$ is a negative tail reducer by \cref{fig:ntr-g,fig:ntr-c}.
    \begin{ctikz}
      \begin{scope}
        \defc{a/0/1,b/1/1}
        \defc{c/-1/0,d/0/0,e/1/0,f/2/0}
        \drawe{a/b}{}{a/c,a/d,b/e,b/f}
        \drawv{a,b}{}{}{}
        \drawvnl{}{}{}{c,d,e,f}
      \end{scope}
      \begin{scope}[shift={(3.5,0)}]
        \defc{a/0/1,b/1/1,c/0/0,d/1/0}
        \drawe{a/b}{}{a/c,a/d,b/c,b/d}
        \drawv{a,b}{}{}{}
        \drawvnl{}{}{}{c,d}
      \end{scope}
    \end{ctikz}
    \item If $a$ and $b$ share exactly one bluish vertex. Let $c, d$ and $d, e$ be the bluish neighbors of $a$ and $b$ respectively.
    \begin{ctikz}
      \defc{a/0/1,b/1/1,d/1/0,c/0/0,e/2/0}
      \drawe{a/b}{}{a/c,a/d,b/d,b/e}
      \drawv{a,b}{}{}{c,d,e}
    \end{ctikz}
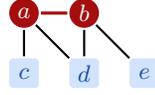
    \begin{kase}
      \item If $d$ is of degree $2$, then $\sset{d}$ is a negative tail reducer by \cref{fig:degree-ntr-a}.
      \item Suppose $d$ is of degree $3$. Let $f$ be the third neighbor of $d$.
      \begin{kase}
        \item Suppose $f$ has no blue neighbor, and $f \sim c$ or $f \sim e$. We may assume that $f \sim e$ up to symmetry. Depending on whether $f\sim c$ and the number of bluish neighbors of $f$, there are three possibilities. In each case, $\sset{b,f}$ is a negative tail reducer by \cref{fig:ntr-b,fig:ntr-e,fig:degree-ntr-e}.
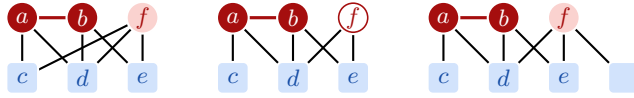
        \begin{ctikz}
          \begin{scope}
            \defc{a/0/1,b/1/1,c/0/0,d/1/0,e/2/0,f/2/1}
            \drawe{a/b}{}{a/c,a/d,b/d,b/e,f/c,f/d,f/e}
            \drawv{a,b}{f}{}{c,d,e}
          \end{scope}
          \begin{scope}[shift={(3.5,0)}]
            \defc{a/0/1,b/1/1,c/0/0,d/1/0,e/2/0,f/2/1}
            \drawe{a/b}{}{a/c,a/d,b/d,b/e,f/d,f/e}
            \drawv{a,b}{}{}{c,d,e}
            \node[a-vertex] at (f) {$f$};
          \end{scope}
          \begin{scope}[shift={(7,0)}]
            \defc{a/0/1,b/1/1,c/0/0,d/1/0,e/2/0,f/2/1,g/3/0}
            \drawe{a/b}{}{a/c,a/d,b/d,b/e,f/d,f/e,f/g}
            \drawv{a,b}{f}{}{c,d,e}
            \drawvnl{}{}{}{g}
          \end{scope}
        \end{ctikz}
        \item Suppose $f$ has no blue neighbor, and $f \nsim c$ and $f \nsim e$. Depending on the number of bluish neighbors of $f$, there are three possibilities. In each case, $\sset{a,b,f}$ is a negative tail reducer by \cref{fig:ntr-t,fig:ntr-n,fig:degree-ntr-t-two,fig:degree-ntr-t-one}. 
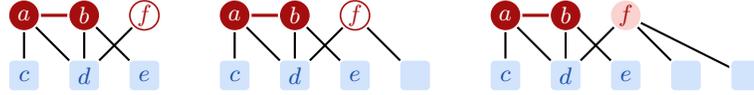
\label{2.2.2}
        \begin{ctikz}
          \begin{scope}
            \defc{a/0/1,b/1/1,c/0/0,d/1/0,e/2/0,f/2/1}
            \drawe{a/b}{}{a/c,a/d,b/d,b/e,f/d}
            \drawv{a,b}{}{}{c,d,e}
            \node[a-vertex] at (f) {$f$};
          \end{scope}
          \begin{scope}[shift={(3.5,0)}]
            \defc{a/0/1,b/1/1,c/0/0,d/1/0,e/2/0,f/2/1,g/3/0}
            \drawe{a/b}{}{a/c,a/d,b/d,b/e,f/d,f/g}
            \drawv{a,b}{}{}{c,d,e}
            \drawvnl{}{}{}{g}
            \node[a-vertex] at (f) {$f$};
          \end{scope}
          \begin{scope}[shift={(8,0)}]
            \defc{a/0/1,b/1/1,c/0/0,d/1/0,e/2/0,f/2/1,g/3/0,h/4/0}
            \drawe{a/b}{}{a/c,a/d,b/d,b/e,f/d,f/g,f/h}
            \drawv{a,b}{f}{}{c,d,e}
            \drawvnl{}{}{}{g,h}
          \end{scope}
        \end{ctikz}
        \item Suppose $f$ has a blue neighbor, and either $f \sim c$ or $f \sim e$. Let $g$ be the blue neighbor of $f$. We may assume that $f \sim e$ up to symmetry.
        \begin{ctikz}
          \defc{a/0/1,b/1/1,c/0/0,d/1/0,e/2/0,f/2/1, g/3/0}
          \drawe{a/b}{}{a/c,a/d,b/d,b/e,f/d,f/e,f/g}
          \drawv{a,b}{f}{g}{c,d,e}
        \end{ctikz}
        \begin{kase}
          \item Suppose $g$ has no red neighbor. Let $h$ be the second reddish neighbor of $g$. Then $\sset{d,g}$ is a negative tail reducer by \cref{fig:ntr-h}.
          \begin{ctikz}
            \defc{a/0/1,b/1/1,c/0/0,d/1/0,e/2/0,f/2/1,g/3/0,h/3/1}
            \drawe{a/b}{}{a/c,a/d,b/d,b/e,f/d,f/e,f/g,g/h}
            \drawv{a,b}{f,h}{g}{c,d,e}
          \end{ctikz}
          \item Suppose $g$ has a red neighbor, say $h$. \ce{d} implies that $h \nsim e$ because otherwise $\sset{b,e,f,g,h}$ is a cut enhancer. Flip the cut preserver $\sset{g,h}$. Then $\sset{b,f}$ becomes a negative tail reducer by \cref{fig:ntr-b}.
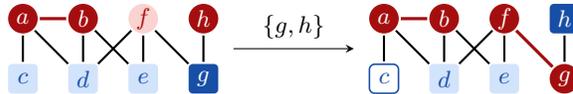
          \begin{ctikz}
            \begin{scope}
              \defc{a/0/1,b/1/1,c/0/0,d/1/0,e/2/0,f/2/1,g/3/0,h/3/1}
              \drawe{a/b}{}{a/c,a/d,b/d,b/e,f/d,f/e,f/g,g/h}
              \drawv{a,b,h}{f}{g}{c,d,e}
            \end{scope}
            \draw[-stealth, thin] (3.5,.5) -- (5.5,.5) node[midway, above, font=\footnotesize] {$\sset{g,h}$};
            \begin{scope}[shift={(6,0)}]
              \defc{a/0/1,b/1/1,c/0/0,d/1/0,e/2/0,f/2/1,g/3/0,h/3/1}
              \drawe{a/b,f/g}{}{a/c,a/d,b/d,b/e,f/d,f/e,g/h}
              \drawv{a,b,f,g}{}{h}{d,e}
              \node[b-vertex] at (c) {$c$};
            \end{scope}
          \end{ctikz}
        \end{kase}
        \item Suppose $f$ has a blue neighbor, $f \nsim c$, $f \nsim e$, and $f$ is red. Let $g$ be the blue neighbor of $f$. Flip the cut preserver $\sset{f,g}$. Then $\sset{d}$ become a negative tail reducer by \cref{fig:ntr-a}. 
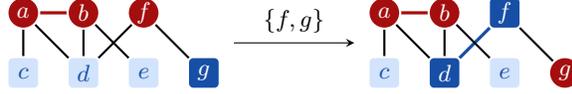
\label{2.2.4}
        \begin{ctikz}
          \begin{scope}
            \defc{a/0/1,b/1/1,c/0/0,d/1/0,e/2/0,f/2/1,g/3/0}
            \drawe{a/b}{}{a/c,a/d,b/d,b/e,f/d,f/g}
            \drawv{a,b,f}{}{g}{c,d,e}
          \end{scope}
          \draw[-stealth, thin] (3.5,0.5) -- (5.5,0.5) node[midway, above, font=\footnotesize] {$\sset{f,g}$};
          \begin{scope}[shift={(6,0)}]
            \defc{a/0/1,b/1/1,c/0/0,d/1/0,e/2/0,f/2/1, g/3/0}
            \drawe{a/b}{f/d}{a/c,a/d,b/d,b/e,f/g}
            \drawv{a,b,g}{}{d,f}{c,e}
          \end{scope}
        \end{ctikz}
        \item Suppose $f$ has a blue neighbor, $f \nsim c$, $f \nsim e$, and $f$ is reddish. Let $g$ be a blue neighbor of $f$, and let $h$ be the second neighbor of $g$ in $A$. \label{2.2.5}
        \begin{kase}
          \item If $h$ is reddish, then $\sset{d,g}$ is a negative tail reducer by \cref{fig:ntr-h}.
          \begin{ctikz}
            \defc{a/0/1,b/1/1,c/0/0,d/1/0,e/2/0,f/2/1,g/3/0,h/3/1}
            \drawe{a/b}{}{a/c,a/d,b/d,b/e,f/d,f/g,g/h}
            \drawv{a,b}{f,h}{g}{c,d,e}
          \end{ctikz}
          \item Suppose $h$ is red, $h \nsim c$ and $h \nsim e$. Flipping the cut preserver $\sset{g,h}$ reduces to \ref{2.2.2} or \ref{2.2.4}.
          \begin{ctikz}
            \begin{scope}
              \defc{a/0/1,b/1/1,c/0/0,d/1/0,e/2/0,f/2/1,g/3/0,h/3/1}
              \drawe{a/b}{}{a/c,a/d,b/d,b/e,f/d,f/g,g/h}
              \drawv{a,b,h}{f}{g}{c,d,e}
            \end{scope}
            \draw[-stealth, thin] (3.5,0.5) -- (5.5,0.5) node[midway, above, font=\footnotesize] {$\sset{g,h}$};
            \begin{scope}[shift={(6,0)}]
              \defc{a/0/1,b/1/1,c/0/0,d/1/0,e/2/0,f/2/1,g/3/0,h/3/1}
              \drawe{a/b}{}{a/c,a/d,b/d,b/e,f/d,f/g,g/h}
              \drawv{a,b,f,g}{}{h}{c,d,e}
            \end{scope}
          \end{ctikz}
          \item Suppose $h$ is red, $h \sim c$ or $h \sim e$. We may assume that $h \sim e$ up to symmetry. Let $i$ be the blue neighbor of $g$. We may assume that $f$ is of degree $3$ because otherwise $\sset{f,h}$ is an absolute tail reducer by \cref{lem:deg-3-flip}. 
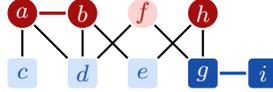
\label{case:isolated-neg-2.2.5.3}
          \begin{ctikz}
            \defc{a/0/1,b/1/1,c/0/0,d/1/0,e/2/0,f/2/1,g/3/0,h/3/1,i/4/0}
            \drawe{a/b}{g/i}{a/c,a/d,b/d,b/e,f/d,f/g,g/h,h/e}
            \drawv{a,b,h}{f}{g,i}{c,d,e}
          \end{ctikz}
          \begin{kase}
            \item If $f \sim i$, then $\sset{b,f}$ is a negative tail reducer by \cref{fig:ntr-h}.
            \begin{ctikz}
              \defc{a/0/1,b/1/1,c/0/0,d/1/0,e/2/0,f/2/1,g/3/0,h/3/1,i/4/0}
              \drawe{a/b}{g/i}{a/c,a/d,b/d,b/e,f/d,f/g,g/h,h/e,f/i}
              \drawv{a,b,h}{f}{g,i}{c,d,e}
            \end{ctikz}
            \item Suppose $f \nsim i$. Let $j$ be the third neighbor of $f$. Suppose for contradiction that $j$ is blue. Then \ce{a} implies that $j \nsim h$ because otherwise $\sset{g,h,j}$ is a cut enhancer, and \ce{e} implies that $\sset{b,d,e,f,g,h,j}$ is a cut enhancer, a contradiction. Therefore $j$ is bluish, and $\sset{a,f,h}$ is a negative tail reducer by \cref{fig:ntr-z}. 
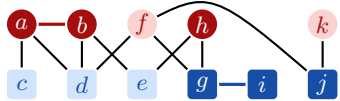
\label{case:isolated-neg-2.2.5.3.2}
            \begin{ctikz}
              \defc{a/0/1,b/1/1,c/0/0,d/1/0,e/2/0,f/2/1,g/3/0,h/3/1,i/4/0,j/5/0}
              \drawe{a/b}{g/i}{a/c,a/d,b/d,b/e,f/d,f/g,g/h,h/e}
              \draw (f) .. controls (3.14,1.6) .. (j);
              \drawv{a,b,h}{f}{g,i}{c,d,e,j}
            \end{ctikz}
          \end{kase}
        \end{kase}
      \end{kase}
    \end{kase}
  \end{kase}

  Finally, notice that the distance between $h$ and $\sset{a,b}$ is $4$ in \ref{2.2.5}, and it is the furthest distance between any vertex in a cut preserver or a negative tail reducer and $\sset{a,b}$.
\end{proof}

The last three lemmas establish that every connected component of $M$ is a path or a cycle of even length if $G$ lacks negative tail reducers. We introduce the following versatile ingredient on pentagons.

\begin{lemma} \label{lem:pentagon-neg-versatile}
  If the pentagon on $\sset{a,b,c,d,e}$ below is an induced subgraph of $G$, and $a$, $b$ and $c$ have no blue neighbors outside $\sset{d,e}$, then a subset of $\sset{a,b,c}$ is a negative tail reducer.
  \begin{ctikz}
    \defc{a/0/1,b/1/1,c/2/1,d/0/0,e/1/0}
    \drawe{a/b}{}{a/d,b/e,c/d,c/e}
    \drawv{a,b}{c}{d}{}
    \node[b-vertex] at (e) {$e$};
  \end{ctikz}
\end{lemma}

\begin{proof}
  We may assume that $c$ is of degree $3$ because otherwise $\sset{a,c}$ is an absolute tail reducer by \cref{lem:deg-3-flip}.
  \begin{kase}
    \item Suppose $e$ is bluish.
    \begin{kase}
      \item If $a$ and $c$ share a bluish neighbor, then $\sset{a,c}$ is a negative tail reducer by \cref{fig:ntr-i}.
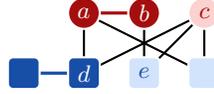
      \begin{ctikz}
        \defc{a/0/1,b/1/1,c/2/1,d/0/0,e/1/0,y/-1/0,z/2/0}
        \drawe{a/b}{d/y}{a/d,b/e,c/d,c/e,a/z,c/z}
        \drawv{a,b}{c}{d}{e}
        \drawvnl{}{}{y}{z}
      \end{ctikz}
      \item Suppose $a$ and $c$ share no bluish neighbors. Depending on whether $b$ shares a bluish neighbor with $a$ and $c$, there are three possibilities. In each case, $\sset{a,b,c}$ is a negative tail reducer by \cref{fig:ntr-s,fig:ntr-u,fig:ntr-aa}.
      \begin{ctikz}
        \begin{scope}
          \defc{a/0/1,b/1/1,c/2/1,d/0/0,e/1/0,y/-1/0,u/2/0,v/3/0}
          \drawe{a/b}{d/y}{a/d,b/e,c/d,c/e,a/u,b/u,c/v}
          \drawv{a,b}{c}{d}{e}
          \drawvnl{}{}{y}{u,v}
        \end{scope}
        \begin{scope}[shift={(5.5,0)}]
          \defc{a/0/1,b/1/1,c/2/1,d/0/0,e/1/0,y/-1/0,u/2/0,v/3/0}
          \drawe{a/b}{d/y}{a/d,b/e,c/d,c/e,a/u,b/v,c/v}
          \drawv{a,b}{c}{d}{e}
          \drawvnl{}{}{y}{u,v}
        \end{scope}
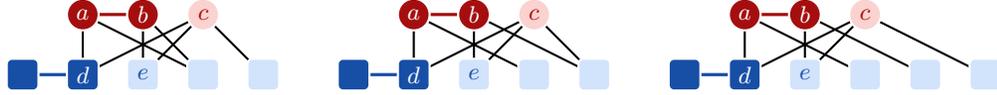
        \begin{scope}[shift={(11,0)}]
          \defc{a/0/1,b/1/1,c/2/1,d/0/0,e/1/0,y/-1/0,u/2/0,v/3/0,w/4/0}
          \drawe{a/b}{d/y}{a/d,b/e,c/d,c/e,a/u,b/v,c/w}
          \drawv{a,b}{c}{d}{e}
          \drawvnl{}{}{y}{u,v,w}
        \end{scope}
      \end{ctikz}
    \end{kase}
    \item Suppose $e$ is blue. Let $f$ be the bluish neighbor of $c$. \ce{d} implies that $f \nsim a$ and $f \nsim b$ because otherwise $\sset{a,c,d,e,f}$ or $\sset{b,c,d,e,f}$ is a cut enhancer. Depending on whether $a$ and $b$ share a bluish neighbor, there are two possibilities. In each case, $\sset{a,b,c}$ is a negative tail reducer by \cref{fig:ntr-y,fig:ntr-af}.
    \begin{ctikz}
      \begin{scope}
        \defc{a/0/1,b/1/1,c/2/1,d/0/0,e/1/0,y/-1/0,z/2/0,f/3/0,u/-2/0}
        \arcthreeone{u}
        \drawe{a/b}{d/y,e/z}{a/d,b/e,c/d,c/e,c/f,a/u}
        \drawv{a,b}{c}{d,e}{f}
        \drawvnl{}{}{y,z}{u}
      \end{scope}
      \begin{scope}[shift={(6.5,0)}]
        \defc{a/0/1,b/1/1,c/2/1,d/0/0,e/1/0,y/-1/0,z/2/0,x/3/0,u/-2/0,f/4/0}
        \drawe{a/b}{d/y,e/z}{a/d,b/e,c/d,c/e,c/f,a/u,b/x}
        \drawv{a,b}{c}{d,e}{f}
        \drawvnl{}{}{y,z}{u,x}
      \end{scope}
    \end{ctikz}
  \end{kase}

  Finally, notice that every negative tail reducer above is a subset of $\sset{a,b,c}$.
\end{proof}

As an application of \cref{lem:pentagon-neg-versatile}, we build the following tool to deal with paths and cycles.

\begin{lemma} \label{lem:pentagon-neg-2}
  If $abcdef$ is an induced path in $G$ such that $ab$ and $ef$ are red, and $cd$ is blue, then $b$ and $e$ do not share a neighbor, or a subset of $G(cd, 2)$ is a negative tail reducer.
\end{lemma}

\begin{proof}
  Suppose $b$ and $e$ share a neighbor, say $g$. \cref{lem:blue-edge} implies that $g$ is bluish. We may assume that $g$ is of degree $3$ because otherwise $\sset{d,g}$ is an absolute tail reducer by \cref{lem:deg-3-flip}.
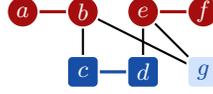
  \begin{ctikz}
    \defc{a/0/1,b/1/1,c/1/0,d/2/0,e/2/1,f/3/1,g/3/0}
    \drawe{a/b,e/f}{c/d}{b/c,d/e,g/b,g/e}
    \drawv{a,b,e,f}{}{c,d}{g}
  \end{ctikz}
  \begin{kase}
    \item Suppose $g \sim f$. Let $h$ be the third neighbor of $c$. \cref{lem:blue-edge} implies that $h$ is reddish. Then $\sset{c,g}$ is a negative tail reducer by \cref{fig:ntr-j}.
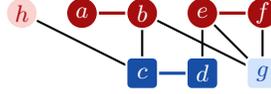
    \begin{ctikz}
      \defc{a/0/1,b/1/1,c/1/0,d/2/0,e/2/1,f/3/1,g/3/0,h/-1/1}
      \drawe{a/b,e/f}{c/d}{b/c,d/e,g/b,g/e,g/f,c/h}
      \drawv{a,b,e,f}{h}{c,d}{g}
    \end{ctikz}
    \item Suppose $g \nsim f$. By a symmetric argument, we may assume that $g \nsim a$. Let $h$ be the third neighbor of $g$. \cref{lem:blue-edge} implies that the third neighbors of $c$ and $d$ are both reddish. \ce{c} implies that $h$ is reddish because otherwise $\sset{b,c,e,g,h}$ is a cut enhancer. Then a subset of $\sset{c,d,g}$ is a negative tail reducer by \cref{lem:pentagon-neg-versatile}.
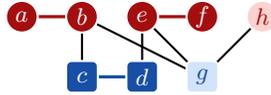
    \begin{ctikz}
      \defc{a/0/1,b/1/1,c/1/0,d/2/0,e/2/1,f/3/1,g/3/0,h/4/1}
      \drawe{a/b,e/f}{c/d}{b/c,d/e,g/b,g/e,g/h}
      \drawv{a,b,e,f}{h}{c,d}{g}
    \end{ctikz}
  \end{kase}

  Finally, notice that the distance between $g$ and $cd$ is $2$, and it is the furthest distance between any vertex in a negative tail reducer and $cd$.
\end{proof}

\begin{remark}
  We do not need the corresponding version of \cref{lem:pentagon-neg-2} for positive tail reducers in \cref{sec:pos} because of \cref{fig:ptr-c}.
\end{remark}

\begin{lemma} \label{lem:path-5-neg-8}
  If the graph $M(\alpha, 4)$ is a path on at least $5$ vertices, or a cycle on at least $6$ vertices, then a subset of $G(\alpha, 8)$ is a negative tail reducer.
\end{lemma}

\begin{proof}
  Observe that $M(\alpha, 4)$ contains an induced path on $5$ vertices, say $ab,cd,ef,gh,ij$, such that $\alpha \in \sset{ab,cd,ef,gh,ij}$. In view of \cref{lem:induced-path}, we may assume that $abcdefghij$ is an induced path in $G$, and moreover we may assume that $ab,ef,ij$ are red, and $cd, gh$ are blue.

  Let $k$ be the third neighbor of $e$. \cref{lem:blue-edge} implies that $k$ is bluish. \ce{b} implies that $k \nsim i$ because otherwise $\sset{d,e,h,i,k}$ is a cut enhancer.
  \begin{kase}
    \item If $k \sim b$, then a subset of $G(cd,2)$ is a negative tail reducer by \cref{lem:pentagon-neg-2}.
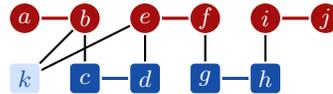
    \begin{ctikz}
      \defc{a/0/1,b/1/1,c/1/0,d/2/0,e/2/1,f/3/1,g/3/0,h/4/0,i/4/1,j/5/1,k/0/0}
      \drawe{a/b,e/f,i/j}{c/d,g/h}{b/c,d/e,e/k,f/g,h/i,k/b}
      \drawv{a,b,e,f,i,j}{}{c,d,g,h}{k}
    \end{ctikz}
    \item Suppose $k \nsim b$. Flip the cut preserver $\sset{b,c}$, and then flip the cut preserver $\sset{h,i}$. Two neighbors of $e$ become bluish. By a symmetric argument, we may assume that two neighbors of $f$ become bluish as well. Then the red edge $ef$ becomes isolated in $M$, and so a subset of $G(ef, 4)$ becomes a negative tail reducer by \cref{lem:isolated-neg-4}.
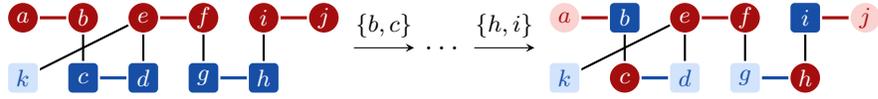
    \begin{ctikz}
      \begin{scope}
        \defc{a/0/1,b/1/1,c/1/0,d/2/0,e/2/1,f/3/1,g/3/0,h/4/0,i/4/1,j/5/1,k/0/0}
        \drawe{a/b,e/f,i/j}{c/d,g/h}{b/c,d/e,e/k,f/g,h/i}
        \drawv{a,b,e,f,i,j}{}{c,d,g,h}{k}
      \end{scope}
      \draw[-stealth, thin] (5.5,0.5) -- (6.5,0.5) node[midway, above, font=\footnotesize] {$\sset{b,c}$};
      \node at (7,0.5) {$\dots$};
      \draw[-stealth, thin] (7.5,0.5) -- (8.5,0.5) node[midway, above, font=\footnotesize] {$\sset{h,i}$};
      \begin{scope}[shift={(9,0)}]
        \defc{a/0/1,b/1/1,c/1/0,d/2/0,e/2/1,f/3/1,g/3/0,h/4/0,i/4/1,j/5/1,k/0/0}
        \drawe{a/b,e/f,i/j}{c/d,g/h}{b/c,d/e,e/k,f/g,h/i}
        \drawv{c,e,f,h}{a,j}{b,i}{d,g,k}
      \end{scope}
    \end{ctikz}
  \end{kase}

  Finally, notice that the distance between $a$ and $ij$ is $8$, and it is the furthest distance between any vertex in a cut preserver or a tail reducer and $\alpha \in \sset{ab,cd,ef,gh,ij}$.
\end{proof}

\begin{lemma} \label{lem:cycle-4-neg-10}
  If $M(\alpha, 4)$ is a cycle on $4$ vertices, then a subset of $G(\alpha, 10)$ is a negative tail reducer.
\end{lemma}

\begin{proof}
  Suppose that $M(\alpha, 4)$ is the cycle $ab,cd,ef,gh$ such that $\alpha \in \sset{ab,cd,ef,gh}$. In view of \cref{lem:induced-cycle}, we may assume that $abcdefgh$ is an induced cycle in $G$, and moreover we may assume that $ab$ and $ef$ are red, and $cd$ and $gh$ are blue.
  \begin{ctikz}
    \defc{a/0/1,b/1/1,c/1/0,d/2/0,e/2/1,f/3/1,g/3/0,h/4/0}
    \drawe{a/b,e/f}{c/d,g/h}{b/c,d/e,f/g}
    \draw (a) to [bend right] (2,0.5) to [bend left] (h);
    \drawv{a,b,e,f}{}{c,d,g,h}{}
  \end{ctikz}

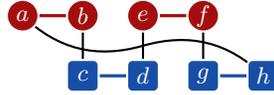
  \cref{lem:blue-edge} implies that the third neighbors of $a,b,e,f$ are all bluish. \ce{b} implies that $a$ and $e$ do not share a bluish neighbor because otherwise $\sset{a,d,e,h,\circ}$ is a cut enhancer, where $\circ$ is the shared bluish neighbor of $a$ and $e$. Similarly $b$ and $f$ do not share a bluish neighbor.
  \begin{kase}
    \item If $b$ and $e$ share a bluish neighbor, then a subset of $G(cd, 2)$ is a negative tail reducer by \cref{lem:pentagon-neg-2}.
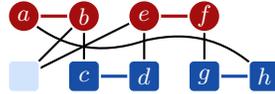
    \begin{ctikz}
      \defc{a/0/1,b/1/1,c/1/0,d/2/0,e/2/1,f/3/1,g/3/0,h/4/0,i/0/0}
      \drawe{a/b,e/f}{c/d,g/h}{b/c,d/e,f/g,e/i,b/i}
      \draw (a) to [bend right] (2,0.59) to [bend left] (h);
      \drawv{a,b,e,f}{}{c,d,g,h}{}
      \drawvnl{}{}{}{i}
    \end{ctikz}
    \item Suppose $b$ and $e$ do not share a bluish neighbor. By a symmetric argument, we may assume that neither do $a$ and $f$, and moreover, neither $c$ and $h$ nor $d$ and $g$ share a reddish neighbor. \label{item:cycle-4-no-sharing}
    \begin{kase}
      \item Suppose $a$ and $b$ do not share a bluish neighbor, or $e$ and $f$ do not share a bluish neighbor. There are two possibilities up to symmetry. In each case, $\sset{a,b,e,f}$ is a negative tail reducer by \cref{fig:ntr-ak,fig:ntr-al}.
      \begin{ctikz}
        \begin{scope}
          \defc{a/0/1,b/1/1,c/1/0,d/2/0,e/2/1,f/3/1,g/3/0,h/4/0,i/-1/0,j/0/0,k/5/0}
          \drawe{a/b,e/f}{c/d,g/h}{b/c,d/e,f/g,a/i,b/j,f/k}
          \draw (a) to [bend right] (2,0.5) to [bend left] (h);
          \draw (e) .. controls (3.14,1.6) .. (k);
          \drawv{a,b,e,f}{}{c,d,g,h}{}
          \drawvnl{}{}{}{i,j,k}
        \end{scope}
        \begin{scope}[shift={(7.5,0)}]
          \defc{a/0/1,b/1/1,c/1/0,d/2/0,e/2/1,f/3/1,g/3/0,h/4/0,i/-1/0,j/0/0,k/6/0,l/5/0}
          \drawe{a/b,e/f}{c/d,g/h}{b/c,d/e,f/g,a/i,b/j,f/l}
          \draw (a) to [bend right] (2,0.5) to [bend left] (h);
          \draw (e) .. controls (3,1.6) .. (k);
          \drawv{a,b,e,f}{}{c,d,g,h}{}
          \drawvnl{}{}{}{i,j,k,l}
        \end{scope}
      \end{ctikz}
      \item Suppose $a$ and $b$ share a bluish neighbor, and $e$ and $f$ share a bluish neighbor. Let $i$ be the shared bluish neighbor of $a$ and $b$, and let $j$ be that of $e$ and $f$. By a symmetric argument, we may assume that $c$ and $d$ share a reddish neighbor, say $k$. Recall from \ref{item:cycle-4-no-sharing} that $k \nsim g$ and $k \nsim h$. We may assume that $k$ is of degree $3$ because otherwise $\sset{b,k}$ is an absolute tail reducer by \cref{lem:deg-3-flip}.
      \begin{kase}
        \item Suppose $k \sim i$ or $k \sim j$. We may assume that $k \sim i$ up to symmetry. Then $\sset{b,e,k}$ is a negative tail reducer by \cref{fig:ntr-p}.
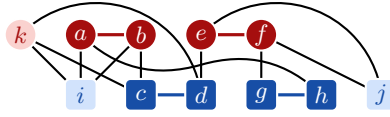
        \begin{ctikz}
            \defc{a/0/1,b/1/1,c/1/0,d/2/0,e/2/1,f/3/1,g/3/0,h/4/0,i/0/0,j/5/0,k/-1/1,l/-1/0}
            \drawe{a/b,e/f}{c/d,g/h}{b/c,d/e,f/g,a/i,b/i,f/j,c/k,k/i}
            \arcthreeminusone{k}
            \draw (a) to [bend right] (2,0.5) to [bend left] (h);
            \arcthreeminusone{e}
            \drawv{a,b,e,f}{k}{c,d,g,h}{i,j}
        \end{ctikz}
        \item Suppose $k \nsim i$ and $k \nsim j$. Let $l$ be the third neighbor of $k$.
        \begin{kase}
          \item If $l$ is bluish, then $\sset{a,b,e,f,k}$ is a negative tail reducer by \cref{fig:ntr-an}. \label{2.2.1}
          \begin{ctikz}
            \defc{a/0/1,b/1/1,c/1/0,d/2/0,e/2/1,f/3/1,g/3/0,h/4/0,i/0/0,j/5/0,k/-1/1,l/-1/0}
            \drawe{a/b,e/f}{c/d,g/h}{b/c,d/e,f/g,a/i,b/i,f/j,c/k,k/l}
            \arcthreeminusone{k}
            \draw (a) to [bend right] (2,0.5) to [bend left] (h);
            \arcthreeminusone{e}
            \drawv{a,b,e,f}{k}{c,d,g,h}{i,j,l}
          \end{ctikz}
          \item Suppose $l$ is blue. Let $m$ be the blue neighbor of $l$, and let $n$ be the third neighbor of $l$. \ce{b} implies that $n$ is reddish because otherwise $\sset{d,e,k,l,n}$ is a cut enhancer.
          \begin{ctikz}
            \defc{a/0/1,b/1/1,c/1/0,d/2/0,e/2/1,f/3/1,g/3/0,h/4/0,i/0/0,j/5/0,k/-1/1,l/-1/0,m/-2/0,n/-2/1}
            \drawe{a/b,e/f}{c/d,g/h}{b/c,d/e,f/g,a/i,b/i,f/j,c/k,k/l,l/m,l/n}
            \arcthreeminusone{k}
            \draw (a) to [bend right] (2,0.5) to [bend left] (h);
            \arcthreeminusone{e}
            \drawv{a,b,e,f}{k,n}{c,d,g,h,l,m}{i,j}
          \end{ctikz}
          \begin{kase}
            \item If $m$ has no red neighbor, then the blue edge $lm$ is isolated in $M$, and so a subset of $G(lm, 4)$ is a negative tail reducer by \cref{lem:isolated-neg-4}. \label{2.2.2.1}
            \item Suppose $m$ has a red neighbor, say $o$. \ce{b} implies that $o \nsim i$ and $o \nsim j$ because otherwise $\sset{b,c,i,m,o}$ or $\sset{f,g,j,m,o}$ is a cut enhancer. Flip the cut preserver $\sset{m,o}$. Then the current case reduces to \ref{2.2.1}, where $l$ is bluish.
            \begin{ctikz}
              \defc{a/0/1,b/1/1,c/1/0,d/2/0,e/2/1,f/3/1,g/3/0,h/4/0,i/0/0,j/5/0,k/-1/1,l/-1/0,m/-2/0,n/-2/1,o/-3/1}
              \drawe{a/b,e/f}{c/d,g/h}{b/c,d/e,f/g,a/i,b/i,f/j,c/k,k/l,l/m,m/o,l/n}
              \arcthreeminusone{k}
              \draw (a) to [bend right] (2,0.5) to [bend left] (h);
              \arcthreeminusone{e}
              \drawv{a,b,e,f,o}{k,n}{c,d,g,h,l,m}{i,j}
            \end{ctikz}
            \end{kase}
        \end{kase}
      \end{kase}
    \end{kase}
  \end{kase}

  Finally, note that the distance between a vertex in $G(lm,4)$ and $gh$ is at most $10$ in \ref{2.2.2.1}, and it is the furthest distance possible between any vertex in a cut preserver or a negative tail reducer and $\alpha \in \sset{ab,cd,ef,gh}$.
\end{proof}

We confront yet another 7-level case analysis. Fear not, for each step forward brings us closer to clarity, and perhaps a well-deserved coffee break.

\begin{lemma} \label{lem:path-4-neg}
  If $M(\alpha, 4)$ is a path on $4$ vertices, then a subset of $G(\alpha, \pathfourneg)$ is a negative tail reducer.
\end{lemma}

\begin{proof}
  Suppose that $M(\alpha, 4)$ is the path $ab,cd,ef,gh$ such that $\alpha \in \sset{ab,cd,ef,gh}$. In view of \cref{lem:induced-path}, we may assume that $abcdefgh$ is an induced path in $G$, and moreover we may assume that $ab$ and $ef$ are red, and $cd$ and $gh$ are blue. Let $i$ be the third neighbor of $d$. \cref{lem:blue-edge} implies that $i$ is reddish. We may assume that $i$ is of degree $3$ because otherwise $\sset{e, i}$ is an absolute tail reducer by \cref{lem:deg-3-flip}.
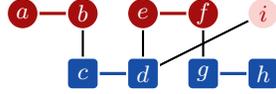
  \begin{ctikz}
    \defc{a/0/1,b/1/1,c/1/0,d/2/0,e/2/1,f/3/1,g/3/0,h/4/0,i/4/1}
    \drawe{a/b,e/f}{c/d,g/h}{b/c,d/e,f/g,d/i}
    \drawv{a,b,e,f}{i}{c,d,g,h}{}
  \end{ctikz}
  \begin{kase}
    \item If $i \sim g$, then a subset of $G(ef, 2)$ is a negative tail reducer by \cref{lem:pentagon-neg-2}.
    \item Suppose $i \nsim g$ and $i \nsim h$. Flip the cut preserver $\sset{d,e}$. Then the blue edge $gh$ becomes isolated in $M$, and so a subset of $G(gh, 4)$ becomes a negative tail reducer by \cref{lem:isolated-neg-4}. 
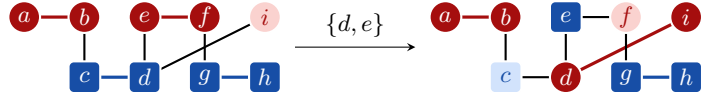
\label{3}
    \begin{ctikz}
      \begin{scope}
        \defc{a/0/1,b/1/1,c/1/0,d/2/0,e/2/1,f/3/1,g/3/0,h/4/0,i/4/1}
        \drawe{a/b,e/f}{c/d,g/h}{b/c,d/e,f/g,d/i}
        \drawv{a,b,e,f}{i}{c,d,g,h}{}
      \end{scope}
      \draw[-stealth, thin] (4.5,0.5) -- (6.5,0.5) node[midway, above, font=\footnotesize] {$\sset{d,e}$};
      \begin{scope}[shift={(7,0)}]
        \defc{a/0/1,b/1/1,c/1/0,d/2/0,e/2/1,f/3/1,g/3/0,h/4/0,i/4/1}
        \drawe{a/b,d/i}{g/h}{c/d,b/c,d/e,f/g,e/f}
        \drawv{a,b,d,i}{f}{e,g,h}{c}
      \end{scope}
    \end{ctikz}
    \item Suppose $i \nsim g$ and $i \sim h$. Let $j$ be the third neighbor of $e$. \cref{lem:blue-edge} implies that $j$ is bluish. By a symmetric argument, $j \sim a$ and $j \nsim b$. Clearly, the third neighbor of $a$ is bluish. \cref{lem:blue-edge} implies that the third neighbor of $b$ is bluish.
    \begin{ctikz}
      \defc{a/0/1,b/1/1,c/1/0,d/2/0,e/2/1,f/3/1,g/3/0,h/4/0,j/0/0,i/4/1}
      \drawe{a/b,e/f}{c/d,g/h}{b/c,d/e,f/g,e/j,a/j,d/i,h/i}
      \drawv{a,b,e,f}{i}{c,d,g,h}{j}
    \end{ctikz}
    \begin{kase}
      \item If $a$ and $b$ share a bluish neighbor, then $\sset{a,b,e}$ is a negative tail reducer by \cref{fig:ntr-l}.
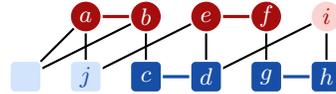
      \begin{ctikz}
        \defc{a/0/1,b/1/1,c/1/0,d/2/0,e/2/1,f/3/1,g/3/0,h/4/0,i/4/1,j/0/0,k/-1/0}
        \drawe{a/b,e/f}{c/d,g/h}{b/c,d/e,f/g,d/i,i/h,j/a,j/e,a/k,b/k}
        \drawv{a,b,e,f}{i}{c,d,g,h}{j}
        \drawvnl{}{}{}{k}
      \end{ctikz}
      \item Suppose $i \sim j$. Flip the cut preserver $\sset{d,e}$. Then $M$ contains multiple edges between $di$ and $ej$, and so a subset of $G(di,2)$ becomes a negative tail reducer by \cref{lem:no-multi-neg-2}.
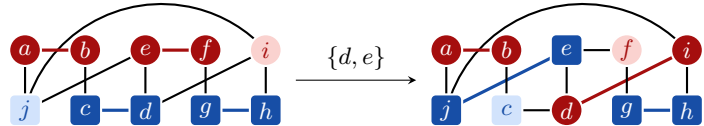
      \begin{ctikz}
        \begin{scope}
          \defc{a/0/1,b/1/1,c/1/0,d/2/0,e/2/1,f/3/1,g/3/0,h/4/0,i/4/1,j/0/0}
          \arcfourone{j}
          \drawe{a/b,e/f}{c/d,g/h}{b/c,d/e,f/g,d/i,i/h,j/a,j/e}
          \drawv{a,b,e,f}{i}{c,d,g,h}{j}
        \end{scope}
        \draw[-stealth, thin] (4.5,0.5) -- (6.5,0.5) node[midway, above, font=\footnotesize] {$\sset{d,e}$};
        \begin{scope}[shift={(7,0)}]
          \defc{a/0/1,b/1/1,c/1/0,d/2/0,e/2/1,f/3/1,g/3/0,h/4/0,i/4/1,j/0/0}
          \arcfourone{j}
          \drawe{a/b,i/d}{g/h,j/e}{b/c,c/d,d/e,e/f,f/g,a/j,i/h}
          \drawv{a,b,d,i}{f}{e,g,h,j}{c}
        \end{scope}
      \end{ctikz}
      \item Suppose $i \sim c$ and $j \sim f$. Clearly, the third neighbor of $a$ is bluish. Then $\sset{a,f,i}$ is a negative tail reducer by \cref{fig:ntr-ab}.
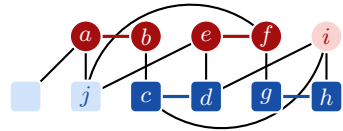
      \begin{ctikz}
        \defc{a/0/1,b/1/1,c/1/0,d/2/0,e/2/1,f/3/1,g/3/0,h/4/0,j/0/0,i/4/1,k/-1/0}
        \arcthreeone{j}
        \arcthreeoneneg{c}
        \drawe{a/b,e/f}{c/d,g/h}{b/c,d/e,f/g,e/j,a/j,d/i,h/i,a/k}
        \drawv{a,b,e,f}{i}{c,d,g,h}{j}
        \drawvnl{}{}{}{k}
      \end{ctikz}
      \item Suppose $a$ and $b$ share no bluish neighbors, $i \nsim j$, and $i \nsim c$ or $j \nsim f$. We may assume that $i \nsim c$ up to symmetry. \ce{c} implies that the third neighbor of $i$ is bluish because otherwise $\sset{d,e,h,i,\circ}$ is a cut enhancer, where $\circ$ is the third blue neighbor of $i$. Let $k$ be the third neighbor of $c$. \cref{lem:blue-edge} implies that $k$ is reddish, and \ce{b} implies that $k \nsim g$ because otherwise $\sset{b,c,f,g,k}$ is a cut enhancer. We may assume that $k$ is of degree $3$ because otherwise $\sset{b,k}$ is an absolute tail reducer by \cref{lem:deg-3-flip}. 
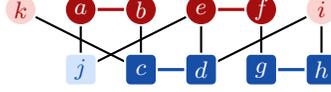
\label{4-path-neg-a-b-no-share}
      \begin{ctikz}
        \defc{a/0/1,b/1/1,c/1/0,d/2/0,e/2/1,f/3/1,g/3/0,h/4/0,i/4/1,j/0/0,k/-1/1}
        \drawe{a/b,e/f}{c/d,g/h}{b/c,d/e,f/g,d/i,i/h,j/a,j/e,k/c}
        \drawv{a,b,e,f}{i,k}{c,d,g,h}{j}
      \end{ctikz}
      \begin{kase}
        \item If $k \sim j$, and the third neighbor of $k$ is bluish, then a subset of $\sset{a,b,k}$ is a negative tail reducer by \cref{lem:pentagon-neg-versatile}. \label{path-4-neg:l-is-bluish}
        \item Suppose $k \sim j$, and the third neighbor of $k$ is blue. \ce{e} implies that $k \nsim h$ because otherwise $\sset{a,d,e,h,i,j,k}$ is a cut enhancer. Let $l$ be the other blue neighbor of $k$. Clearly $l \nsim a$. \ce{b} implies that $l$ has no red neighbors because otherwise $\sset{b,c,k,l,\circ}$ is a cut enhancer, where $\circ$ is the red neighbor of $l$. Let $m$ be the blue neighbor of $l$.
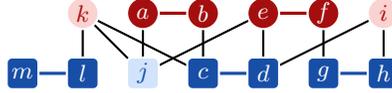
        \begin{ctikz}
          \defc{a/0/1,b/1/1,c/1/0,d/2/0,e/2/1,f/3/1,g/3/0,h/4/0,i/4/1,j/0/0,k/-1/1,l/-1/0,m/-2/0}
          \drawe{a/b,e/f}{c/d,g/h,l/m}{b/c,d/e,f/g,d/i,i/h,j/a,j/e,k/c,j/k,k/l}
          \drawv{a,b,e,f}{i,k}{c,d,g,h,l,m}{j}
        \end{ctikz}
        \begin{kase}
          \item If $m$ has no red neighbors, then $lm$ is isolated in $M$, and so a subset of $G(lm, 4)$ is a negative tail reducer by \cref{lem:isolated-neg-4}.
          \item Suppose $m$ has a red neighbor. Clearly, $m \nsim a$, $m \nsim b$. \cref{lem:blue-edge} implies that the third neighbor of $f$ is bluish, and so $m \nsim f$. Let $n$ be a red neighbor of $m$. \ce{c} implies that $m \nsim i$ because otherwise $\sset{d,e,h,i,m}$ is a cut enhancer. \ce{b} implies that $n \nsim l$ because otherwise $\sset{b,c,k,l,n}$ is a cut enhancer. Flip the cut preserver $\sset{m,n}$. \ce{b} implies that the third neighbor of $a$ is bluish because otherwise $\sset{a,d,e,j,\circ}$ is a cut enhancer, where $\circ$ is the third blue neighbor of $a$. \Cref{lem:blue-edge} implies that the third neighbor of $b$ is bluish. Then a subset of $\sset{a,b,k}$ is a negative tail reducer by \cref{lem:pentagon-neg-versatile}.
          \begin{ctikz}
            \begin{scope}
              \defc{a/0/1,b/1/1,c/1/0,d/2/0,e/2/1,f/3/1,g/3/0,h/4/0,i/4/1,j/0/0,k/-1/1,l/-1/0,m/-2/0,n/-2/1}
              \drawe{a/b,e/f}{c/d,g/h,l/m}{b/c,d/e,f/g,d/i,i/h,j/a,j/e,k/c,j/k,k/l,m/n}
              \drawv{a,b,e,f,n}{i,k}{c,d,g,h,l,m}{j}
            \end{scope}
            \begin{scope}[shift={(9,0)}]
              \draw[-stealth, thin] (-4.5,.5) -- (-2.5,.5) node[midway,above,font=\footnotesize] {$\sset{m,n}$};
              \defc{a/0/1,b/1/1,c/1/0,d/2/0,e/2/1,f/3/1,g/3/0,h/4/0,i/4/1,j/0/0,k/-1/1,l/-1/0,m/-2/0,n/-2/1}
              \drawe{a/b,e/f}{c/d,g/h}{b/c,d/e,f/g,d/i,i/h,j/a,j/e,k/c,j/k,k/l,m/n,l/m}
              \drawv{a,b,e,f,m}{i,k}{c,d,g,h,n}{j,l}
            \end{scope}
          \end{ctikz}
        \end{kase}
        \item Suppose $k \nsim j$ and $k \sim h$. Let $l$ be the third neighbor of $k$. \ce{c} implies that $l$ is bluish because otherwise $\sset{b,c,h,k,l}$ is a cut enhancer. \ce{b} implies that $b$ and $f$ do not share a bluish neighbor because otherwise $\sset{b,c,f,g,\circ}$ is a cut enhancer, where $\circ$ is the shared bluish neighbor of $b$ and $f$. Flip the cut preserver $\sset{b,c}$, and then the cut preserver $\sset{h,k}$.
        \begin{ctikz}
          \begin{scope}
            \defc{a/0/1,b/1/1,c/1/0,d/2/0,e/2/1,f/3/1,g/3/0,h/4/0,i/4/1,j/0/0,k/-1/1,l/-1/0}
            \arcfiveminusone{k}
            \drawe{a/b,e/f}{c/d,g/h}{b/c,d/e,f/g,d/i,i/h,j/a,j/e,k/c,k/l}
            \drawv{a,b,e,f}{i,k}{c,d,g,h}{j,l}
          \end{scope}
          \draw[-stealth, thin] (4.5,0.5) -- (5.5,0.5) node[midway, above, font=\footnotesize] {$\sset{b,c}$};
          \draw[-stealth, thin] (6.5,0.5) -- (7.5,0.5) node[midway, above, font=\footnotesize] {$\sset{h,k}$};
          \node at (6,0.5) {$\dots$};
          \begin{scope}[shift={(9,0)}]
            \defc{a/0/1,b/1/1,c/1/0,d/2/0,e/2/1,f/3/1,g/3/0,h/4/0,i/4/1,j/0/0,k/-1/1,l/-1/0}
            \arcfiveminusone{k}
            \drawe{e/f,i/h}{k/l}{a/b,b/c,c/k,c/d,d/e,f/g,g/h,d/i,e/j,j/a}
            \drawv{e,f,i,h}{a,c}{b,k,l}{d,g,j}
          \end{scope}
        \end{ctikz}
        \begin{kase}
          \item Suppose $l \nsim f$. Since $b$ and $f$ do not share a neighbor, $ef$ becomes isolated in $M$, and so a subset of $G(ef,4)$ becomes a negative tail reducer by \cref{lem:isolated-neg-4}.
          \item Suppose $l \sim f$. Let $m$ be the third neighbor of $l$. \ce{a} implies that $m$ is reddish because otherwise $\sset{f,l,m}$ is a cut enhancer. \ce{c} implies that the third neighbor of $g$ is reddish because otherwise $\sset{f,g,h,k,\circ}$ is a cut enhancer, where $\circ$ is the third red neighbor of $g$. Then a subset of $\sset{g,l,k}$ is a negative tail reducer by \cref{lem:pentagon-neg-versatile}.
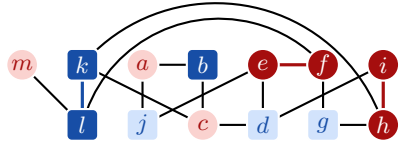
          \begin{ctikz}
            \defc{a/0/1,b/1/1,c/1/0,d/2/0,e/2/1,f/3/1,g/3/0,h/4/0,i/4/1,j/0/0,k/-1/1,l/-1/0,m/-2/1}
            \arcfiveminusone{k}
            \arcfourone{l}
            \drawe{e/f,i/h}{k/l}{a/b,b/c,c/k,c/d,d/e,f/g,g/h,d/i,e/j,j/a,l/m}
            \drawv{e,f,i,h}{a,c,m}{b,k,l}{d,g,j}
          \end{ctikz}
        \end{kase}
        \item Suppose $k \nsim j$ and $k \nsim h$. Let $l$ and $m$ be the other neighbors of $k$.
        \begin{ctikz}
          \defc{a/0/1,b/1/1,c/1/0,d/2/0,e/2/1,f/3/1,g/3/0,h/4/0,i/4/1,j/0/0,k/-1/1}
          \drawe{a/b,e/f}{c/d,g/h}{b/c,d/e,f/g,d/i,i/h,j/a,j/e,k/c}
          \drawv{a,b,e,f}{i,k}{c,d,g,h}{j}
        \end{ctikz}
        \begin{kase}
          \item Suppose $l$ and $m$ are bluish. 
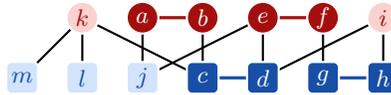
\label{4-path-neg-l-m-bluish}
          \begin{ctikz}
            \defc{a/0/1,b/1/1,c/1/0,d/2/0,e/2/1,f/3/1,g/3/0,h/4/0,i/4/1,j/0/0,k/-1/1,l/-1/0,m/-2/0}
            \drawe{a/b,e/f}{c/d,g/h}{b/c,d/e,f/g,d/i,i/h,j/a,j/e,k/c,k/l,k/m}
            \drawv{a,b,e,f}{i,k}{c,d,g,h}{j,l,m}
          \end{ctikz}
          \begin{kase}
            \item If either $l \sim b$ or $m \sim b$, then $\sset{b, k}$ is a negative tail reducer by \cref{fig:ntr-i}.
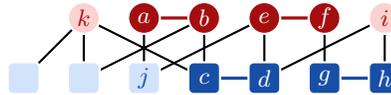
            \begin{ctikz}
              \defc{a/0/1,b/1/1,c/1/0,d/2/0,e/2/1,f/3/1,g/3/0,h/4/0,i/4/1,j/0/0,k/-1/1,l/-1/0,m/-2/0}
              \drawe{a/b,e/f}{c/d,g/h}{b/c,d/e,f/g,d/i,i/h,j/a,j/e,k/c,b/l,l/k,k/m}
              \drawv{a,b,e,f}{i,k}{c,d,g,h}{j}
              \drawvnl{}{}{}{l,m}
            \end{ctikz}
            \item If $l \nsim b$, $m \nsim b$, $l \nsim a$ and $m \nsim a$, then $\sset{a,b,k}$ is a negative tail reducer by \cref{fig:ntr-ag}.
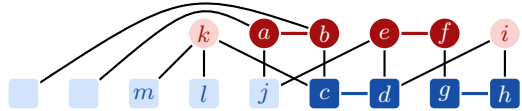
            \begin{ctikz}
              \defc{a/0/1,b/1/1,c/1/0,d/2/0,e/2/1,f/3/1,g/3/0,h/4/0,i/4/1,j/0/0,k/-1/1,l/-1/0,m/-2/0,n/-3/0,o/-4/0}
              \drawe{a/b,e/f}{c/d,g/h}{b/c,d/e,f/g,d/i,i/h,j/a,j/e,k/c,k/l,k/m}
              \draw (a) .. controls (-1.14,1.6) .. (n);
              \draw (b) .. controls (-1.3,1.78) .. (o);
              \drawv{a,b,e,f}{i,k}{c,d,g,h}{j,l,m}
              \drawvnl{}{}{}{n,o}
            \end{ctikz}
            \item Suppose $l \nsim b$, $m \nsim b$, and $l \sim a$ or $m \sim a$. We may assume that $l \sim a$ up to symmetry. Then a subset of $\sset{a,b,k}$ is a negative tail reducer by \cref{lem:pentagon-neg-versatile}.
            \begin{ctikz}
              \defc{a/0/1,b/1/1,c/1/0,d/2/0,e/2/1,f/3/1,g/3/0,h/4/0,i/4/1,j/0/0,k/-1/1,l/-1/0,m/-2/0}
              \drawe{a/b,e/f}{c/d,g/h}{b/c,d/e,f/g,d/i,i/h,j/a,j/e,k/c,k/l,k/m,l/a}
              \drawv{a,b,e,f}{i,k}{c,d,g,h}{j,l,m}
            \end{ctikz}
          \end{kase}
          \item Suppose $l$ and $m$ are blue. Flip the cut preserver $\sset{b,c}$. \ce{a} implies that $l \sim m$ because otherwise $\sset{k,l,m}$ is a cut enhancer. Then $M$ contains multiple edges between $ck$ and $lm$, and so a subset of $G(ck,2)$ becomes a negative tail reducer by \cref{lem:no-multi-neg-2}.
          \begin{ctikz}
            \defc{a/0/1,b/1/1,c/1/0,d/2/0,e/2/1,f/3/1,g/3/0,h/4/0,i/4/1,j/0/0,k/-1/1,l/-1/0,m/-2/0}
            \drawe{c/k,e/f}{g/h,m/l}{e/j,j/a,a/b,b/c,c/d,d/e,d/i,i/h,f/g,m/k,k/l}
            \drawv{k,c,e,f}{a,i}{b,g,h,m,l}{j,d}
          \end{ctikz}
          \item Suppose exactly one of $l$ and $m$ is bluish, and the other is blue. We may assume that $l$ is bluish and $m$ is blue. Let $n$ be the blue neighbor of $m$. \cref{lem:blue-edge} implies that $m \nsim b$, $m \nsim f$, $n \nsim b$ and $n \nsim f$. \ce{b} implies that $m \nsim a$ and $n \nsim a$ because otherwise $\sset{a,d,e,j,m}$ or $\sset{a,d,e,j,n}$ is a cut enhancer. \ce{c} implies that $m \nsim i$ and $n \nsim i$ because otherwise $\sset{d,e,h,i,m}$ or $\sset{d,e,h,i,n}$ is a cut enhancer. Let $o$ be the third neighbor of $m$. \ce{b} implies that $o$ is reddish because otherwise $\sset{b,c,k,m,o}$ is a cut enhancer.
          \begin{ctikz}
            \defc{a/0/1,b/1/1,c/1/0,d/2/0,e/2/1,f/3/1,g/3/0,h/4/0,i/4/1,j/0/0,k/-1/1,l/-1/0,m/-2/0,n/-3/0,o/-2/1}
            \drawe{a/b,e/f}{c/d,g/h,m/n}{b/c,d/e,f/g,d/i,i/h,j/a,j/e,k/c,k/l,k/m,o/m}
            \drawv{a,b,e,f}{i,k,o}{c,d,g,h,m,n}{j,l}
          \end{ctikz}
          \begin{kase}
            \item If $n$ has no red neighbor, then $mn$ is isolated in $M$, and so a subset of $G(mn,4)$ is a negative tail reducer by \cref{lem:isolated-neg-4}. \label{path-4-neg:longest}
            \item Suppose that $n$ has a red neighbor, say $p$. \cref{lem:blue-edge} implies that $p \nsim g$. Clearly, $p \nsim h$. \ce{b} implies that $p \nsim j$ because otherwise $\sset{d,e,j,n,p}$ is a cut enhancer.
            \begin{ctikz}
              \defc{a/0/1,b/1/1,c/1/0,d/2/0,e/2/1,f/3/1,g/3/0,h/4/0,i/4/1,j/0/0,k/-1/1,l/-1/0,m/-2/0,n/-3/0,o/-2/1,p/-3/1}
              \drawe{a/b,e/f}{c/d,g/h,m/n}{b/c,d/e,f/g,d/i,i/h,j/a,j/e,k/c,k/l,k/m,o/m,p/n}
              \drawv{a,b,e,f,p}{i,k,o}{c,d,g,h,m,n}{j,l}
            \end{ctikz}
            \begin{kase}
              \item Suppose $p \sim l$. \ce{f} implies $l \nsim a$ because otherwise $\sset{a,d,e,j,l,n,p}$ is a cut enhancer. \ce{b} implies $l \nsim b$ and $l \nsim f$ because otherwise $\sset{b,c,l,n,p}$ or $\sset{f,g,l,n,p}$ is a cut enhancer.
              \begin{ctikz}
                \defc{a/0/1,b/1/1,c/1/0,d/2/0,e/2/1,f/3/1,g/3/0,h/4/0,i/4/1,j/0/0,k/-1/1,l/-1/0,m/-2/0,n/-3/0,o/-2/1,p/-3/1}
                \drawe{a/b,e/f}{c/d,g/h,m/n}{b/c,d/e,f/g,d/i,i/h,j/a,j/e,k/c,k/l,k/m,o/m,p/n,p/l}
                \drawv{a,b,e,f,p}{i,k,o}{c,d,g,h,m,n}{j,l}
              \end{ctikz}
              \begin{kase}
                \item Suppose the third neighbor, say $q$, of $l$ is red. \ce{g} implies $q \sim p$ because otherwise $\sset{b,c,k,l,m,p,q}$ is a cut enhancer. Then $\sset{l,m}$ is a negative tail reducer by \cref{fig:ntr-h}.
                \begin{ctikz}
                  \defc{a/0/1,b/1/1,c/1/0,d/2/0,e/2/1,f/3/1,g/3/0,h/4/0,i/4/1,j/0/0,k/-1/1,l/-1/0,m/-2/0,n/-3/0,o/-2/1,p/-3/1,q/-4/1}
                  \drawe{a/b,e/f,p/q}{c/d,g/h,m/n}{b/c,d/e,f/g,d/i,i/h,j/a,j/e,k/c,k/l,k/m,o/m,p/n,p/l}
                  \arcthreeminusone{q}
                  \drawv{a,b,e,f,p,q}{i,k,o}{c,d,g,h,m,n}{j,l}
                \end{ctikz}
                \item Suppose the third neighbor of $l$ is reddish. Let $q$ be the third neighbor of $n$. Then $q$ is reddish, because otherwise $\sset{n}$ is a negative tail reducer by \cref{fig:ntr-a} when $p \sim q$, or $\sset{n,p,q}$ is a cut enhancer by \ce{a} when $p \nsim q$. Then a subset of $\sset{l,m,n}$ is a negative tail reducer by \cref{lem:pentagon-neg-versatile}.
                \begin{ctikz}
                  \defc{a/0/1,b/1/1,c/1/0,d/2/0,e/2/1,f/3/1,g/3/0,h/4/0,i/4/1,j/0/0,k/-1/1,l/-1/0,m/-2/0,n/-3/0,o/-2/1,p/-3/1,q/-4/1}
                  \drawe{a/b,e/f}{c/d,g/h,m/n}{b/c,d/e,f/g,d/i,i/h,j/a,j/e,k/c,k/l,k/m,o/m,p/n,p/l,n/q}
                  \drawv{a,b,e,f,p}{i,k,o,q}{c,d,g,h,m,n}{j,l}
                \end{ctikz}
              \end{kase}
              \item Suppose $p \nsim l$. Flip the cut preserver $\sset{n,p}$. \ce{b} implies that the third neighbor of $a$ is bluish because otherwise $\sset{a,d,e,j,\circ}$ is a cut enhancer, where $\circ$ is the third blue neighbor of $a$. \Cref{lem:blue-edge} implies that the third neighbor of $b$ is bluish. Then the current case reduces to \ref{4-path-neg-l-m-bluish}, where $l$ and $m$ are bluish.
            \end{kase}
          \end{kase}
        \end{kase}
      \end{kase}
    \end{kase}
  \end{kase}

  Finally, notice that the distance between a vertex in $G(mn,4)$ and $gh$ is at most $\pathfourneg$ in \ref{path-4-neg:longest}, and it is the furthest distance possible between any vertex in a cut preserver or a negative tail reducer and $\alpha \in \sset{ab, cd, ef, gh}$.
\end{proof}

\begin{lemma} \label{lem:path-4-subgraph-neg}
  If $M$ contains a path $\alpha\beta\gamma\delta$ as a subgraph, then a subset of $G(\beta, \pathfoursubgraphneg)$ is a negative tail reducer.
\end{lemma}

\begin{proof}
  \begin{kase}
    \item If $M(\beta,4)$ contains multiple edges, then a subset of $G(\beta, 10)$ is a negative tail reducer by \cref{lem:no-multi-neg-2}.
    \item If $M(\beta, 4)$ is a simple graph of maximum degree at least $3$, then a subset of $G(\beta, 10)$ is a negative tail reducer by \cref{lem:m-deg-2-neg-2}.
    \item If $M(\beta,4)$ is a simple graph of maximum degree at most $2$, then $M(\beta,4)$ must be a path on at least $4$ vertices, or an even cycle (on at least $4$ vertices), and so a subset of $G(\beta,\pathfourneg)$ is a negative tail reducer by \cref{lem:path-5-neg-8,lem:cycle-4-neg-10,lem:path-4-neg}.
  \end{kase}
\end{proof}

\begin{lemma} \label{lem:path-3-neg}
  If $M(\alpha, 4)$ is a path on $3$ vertices, then a subset of $G(\alpha, \paththreeneg)$ is a negative tail reducer.
\end{lemma}

\begin{proof}
  Suppose that $M(\alpha, 4)$ is the path $ab,cd,ef$ such that $\alpha \in \sset{ab,cd,ef}$. In view of \cref{lem:deg-3}, we may assume that $abcdef$ is an induced path in $G$, and moreover we may assume that $ab,ef$ are red, and $cd$ is blue. Let $g$ be the third neighbor of $b$. \cref{lem:blue-edge} implies that $g$ is bluish.
  \begin{ctikz}
    \defc{a/0/1,b/1/1,c/1/0,d/2/0,e/2/1,f/3/1,g/3/0}
    \drawe{a/b,e/f}{c/d}{b/c,d/e,b/g}
    \drawv{a,b,e,f}{}{c,d}{g}
  \end{ctikz}
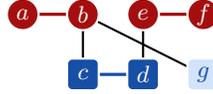
  \begin{kase}
    \item If $g \sim e$, then a subset of $G(cd,2)$ is a negative tail reducer by \cref{lem:pentagon-neg-2}.
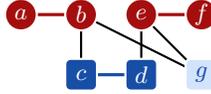
    \begin{ctikz}
      \defc{a/0/1,b/1/1,c/1/0,d/2/0,e/2/1,f/3/1,g/3/0}
      \drawe{a/b,e/f}{c/d}{b/c,d/e,b/g,g/e}
      \drawv{a,b,e,f}{}{c,d}{g}
    \end{ctikz}
    \item Suppose $g \nsim e$ and $g \nsim f$. Flip the cut preserver $\sset{b,c}$. Then the red edge $ef$ becomes isolated in $M$, and so a subset of $G(ef,4)$ is a negative tail reducer by \cref{lem:isolated-neg-4}.
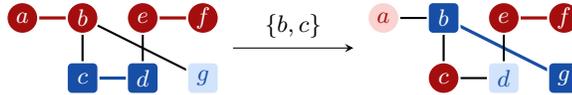
    \begin{ctikz}
      \begin{scope}
        \defc{a/0/1,b/1/1,c/1/0,d/2/0,e/2/1,f/3/1,g/3/0}
        \drawe{a/b,e/f}{c/d}{b/c,d/e,b/g}
        \drawv{a,b,e,f}{}{c,d}{g}
      \end{scope}
      \draw[-stealth, thin] (3.5,0.5) -- (5.5,0.5) node[midway, above, font=\footnotesize] {$\sset{b,c}$};
      \begin{scope}[shift={(6,0)}]
        \defc{a/0/1,b/1/1,c/1/0,d/2/0,e/2/1,f/3/1,g/3/0}
        \drawe{e/f}{b/g}{a/b,c/d,b/c,d/e}
        \drawv{c,e,f}{a}{b,g}{d}
      \end{scope}
    \end{ctikz}
    \item Suppose $g \nsim e$ and $g \sim f$. Let $h$ be the third neighbor of $e$. \cref{lem:blue-edge} implies that $h$ is bluish. By a symmetric argument, we may assume that $h \nsim b$ and $h \sim a$.
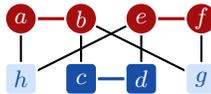
    \begin{ctikz}
      \defc{a/0/1,b/1/1,c/1/0,d/2/0,e/2/1,f/3/1,g/3/0,h/0/0}
      \drawe{a/b,e/f}{c/d}{b/c,d/e,b/g,a/h,e/h,f/g}
      \drawv{a,b,e,f}{}{c,d}{g,h}
    \end{ctikz}
    \begin{kase}
      \item If $g \sim a$ or $h \sim f$, then $\sset{a,b,e}$ or $\sset{b,e,f}$ is a negative tail reducer by \cref{fig:ntr-l}.
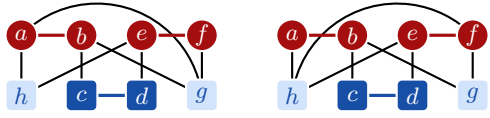
      \begin{ctikz}
        \begin{scope}
          \defc{a/0/1,b/1/1,c/1/0,d/2/0,e/2/1,f/3/1,g/3/0,h/0/0}
          \arcthreeminusone{a}
          \drawe{a/b,e/f}{c/d}{b/c,d/e,b/g,a/h,e/h,f/g}
          \drawv{a,b,e,f}{}{c,d}{g,h}
        \end{scope}
        \begin{scope}[shift={(4.5,0)}]
          \defc{a/0/1,b/1/1,c/1/0,d/2/0,e/2/1,f/3/1,g/3/0,h/0/0}
          \arcthreeone{h}
          \drawe{a/b,e/f}{c/d}{b/c,d/e,b/g,a/h,e/h,f/g}
          \drawv{a,b,e,f}{}{c,d}{g,h}
        \end{scope}
      \end{ctikz}
      \item Suppose $g \nsim a$ and $h \nsim f$. Let $i$ be the third neighbor of $c$. \cref{lem:blue-edge} implies that $i$ is reddish. We may assume that $i$ is of degree $3$ because otherwise $\sset{b,i}$ is an absolute tail reducer by \cref{lem:deg-3-flip}. \ce{d} implies that $i \nsim g$ because otherwise $\sset{b,c,f,g,i}$ is a cut enhancer.
      \begin{ctikz}
        \defc{a/0/1,b/1/1,c/1/0,d/2/0,e/2/1,f/3/1,g/3/0,h/0/0,i/-1/1}
        \drawe{a/b,e/f}{c/d}{b/c,d/e,b/g,a/h,e/h,f/g,i/c}
        \drawv{a,b,e,f}{i}{c,d}{g,h}
      \end{ctikz}
      \begin{kase}
        \item Suppose some neighbor, say $j$, of $i$ outside $\sset{c,d}$ is blue. Let $k$ be the blue neighbor of $j$. Flip the cut preserver $\sset{b,c}$. Then $M$ contains the path $ef,bg,ci,jk$ as a subgraph, and so a subset of $G(bg,10)$ becomes a negative tail reducer by \cref{lem:path-4-subgraph-neg}. 
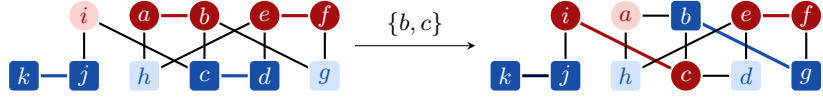
\label{3.2.1-neg}
        \begin{ctikz}
          \begin{scope}
            \defc{a/0/1,b/1/1,c/1/0,d/2/0,e/2/1,f/3/1,g/3/0,h/0/0,i/-1/1,j/-1/0,k/-2/0}
            \drawe{a/b,e/f}{c/d,j/k}{b/c,d/e,b/g,a/h,e/h,f/g,i/c,i/j}
            \drawv{a,b,e,f}{i}{c,d,j,k}{g,h}
          \end{scope}
          \draw[-stealth, thin] (3.5,0.5) -- (5.5,0.5) node[midway, above, font=\footnotesize] {$\sset{b,c}$};
          \begin{scope}[shift={(8,0)}]
            \defc{a/0/1,b/1/1,c/1/0,d/2/0,e/2/1,f/3/1,g/3/0,h/0/0,i/-1/1,j/-1/0,k/-2/0}
            \drawe{i/c,e/f}{b/g,j/k}{a/b,c/d,b/c,d/e,a/h,e/h,f/g,j/k,i/j}
            \drawv{c,e,f,i}{a}{b,g,j,k}{d,h}
          \end{scope}
        \end{ctikz}
        \item Suppose every neighbor of $i$ outside $\sset{c,d}$ is bluish.
        \begin{kase}
          \item Suppose $i \sim d$. \ce{d} implies that $i \nsim h$ because otherwise $\sset{a,d,e,h,i}$ is a cut enhancer. Let $j$ be the third neighbor of $i$. Then $\sset{b,e,i}$ is a negative tail reducer by \cref{fig:ntr-r}.
          \begin{ctikz}
            \defc{a/0/1,b/1/1,c/1/0,d/2/0,e/2/1,f/3/1,g/3/0,h/0/0,i/-1/1,j/-1/0}
            \arcthreeminusoneneg{i}
            \drawe{a/b,e/f}{c/d}{b/c,d/e,b/g,a/h,e/h,f/g,i/c,i/j}
            \drawv{a,b,e,f}{i}{c,d}{g,h,j}
          \end{ctikz}
          \item Suppose $i \nsim d$ and $i \sim h$. Let $j$ be the third neighbor of $i$. Then $\sset{b,e,i}$ is a negative tail reducer by \cref{fig:ntr-q}.
            \begin{ctikz}
              \defc{a/0/1,b/1/1,c/1/0,d/2/0,e/2/1,f/3/1,g/3/0,h/0/0,i/-1/1,j/-1/0}
              \drawe{a/b,e/f}{c/d}{b/c,d/e,b/g,a/h,e/h,f/g,i/c,i/h,i/j}
              \drawv{a,b,e,f}{i}{c,d}{g,h,j}
            \end{ctikz}
          \item Suppose $i \nsim d$ and $i \nsim h$. Let $j$ and $k$ be the bluish neighbors of $i$.
          \begin{kase}
            \item Suppose $j \nsim a$ and $k \nsim a$. Let $l$ be the third neighbor of $a$. Clearly $l$ is bluish. Then $\sset{a,b,i}$ is a negative tail reducer by \cref{fig:ntr-ag}.
            \begin{ctikz}
              \defc{a/0/1,b/1/1,c/1/0,d/2/0,e/2/1,f/3/1,g/3/0,h/0/0,i/-1/1,j/-1/0,k/-2/0,l/-3/0}
              \draw (a) .. controls (-1.14,1.6) .. (l);
              \drawe{a/b,e/f}{c/d}{b/c,d/e,b/g,a/h,e/h,f/g,i/c,i/k,i/j}
              \drawv{a,b,e,f}{i}{c,d}{g,h,j,k,l}
            \end{ctikz}
            \item Suppose $j \sim a$ or $k \sim a$. We may assume that $j \sim a$ up to symmetry. Then $\sset{a,b,i}$ is a negative tail reducer by \cref{fig:ntr-aa}.
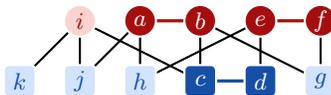
            \begin{ctikz}
              \defc{a/0/1,b/1/1,c/1/0,d/2/0,e/2/1,f/3/1,g/3/0,h/0/0,i/-1/1,j/-1/0,k/-2/0}
              \drawe{a/b,e/f}{c/d}{b/c,d/e,b/g,a/h,e/h,f/g,i/c,i/k,i/j,j/a}
              \drawv{a,b,e,f}{i}{c,d}{g,h,j,k}
            \end{ctikz}
          \end{kase}
        \end{kase}
      \end{kase}
    \end{kase}
  \end{kase}

  Finally notice that the distance between a vertex of $G(bg,\pathfoursubgraphneg)$ and $ef$ is at most $\paththreeneg$ in \ref{3.2.1-neg}, and it is the furthest distance possible between any vertex in a cut preserver or a negative tail reducer and $\alpha \in \sset{ab,cd,ef}$.
\end{proof}

\begin{lemma} \label{lem:path-3-subgraph-neg}
  If $M$ contains a path $\alpha\beta\gamma$ as a subgraph, then a subset of $G(\beta,\paththreesubgraphneg)$ is a negative tail reducer.
\end{lemma}

\begin{proof}
  \begin{kase}
    \item If $M(\beta, 4)$ contains multiple edges, then a subset of $G(\beta, 10)$ is a negative tail reducer by \cref{lem:no-multi-neg-2}.
    \item If $M(\beta, 4)$ is a simple graph of maximum degree at least $3$, then a subset of $G(\beta, 10)$ is a negative tail reducer by \cref{lem:m-deg-2-neg-2}.
    \item If $M(\beta, 4)$ is a simple graph of maximum degree at most $2$, then $M(\beta, 4)$ must be a path on at least $3$ vertices, or an even cycle (on at least $4$ vertices), and so a subset of $G(\beta,\paththreesubgraphneg)$ is a negative tail reducer by \cref{lem:path-5-neg-8,lem:cycle-4-neg-10,lem:path-4-neg,lem:path-3-neg}. \qedhere
  \end{kase}
\end{proof}

Congratulations, brave reader! You have made it to the final boss: a 7-level case analysis. This is the ultimate test of endurance, focus, and caffeine tolerance. Victory is near --- let's slay this proof together!

\begin{lemma} \label{lem:path-2-neg}
  If $M(\alpha, 4)$ is a path on $2$ vertices, then a subset of $G(\alpha, \pathtwoneg)$ is a negative tail reducer.
\end{lemma}

\begin{proof}
  Suppose that $M(\alpha, 4)$ is the path $ab, cd$ such that $\alpha \in \sset{ab,cd}$. In view of \cref{lem:induced-path}, we may assume that $abcd$ is an induced path in $G$, and moreover we may assume that $ab$ is red, and $cd$ is blue. Let $e$ and $f$ be the third neighbors of $b$ and $c$ respectively. \cref{lem:blue-edge} implies that $e$ is bluish, and $f$ is reddish. We may assume that both $e$ and $f$ are of degree $3$ because otherwise $\sset{c,e}$ or $\sset{b,f}$ is an absolute tail reducer by \cref{lem:deg-3-flip}.
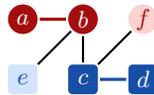
  \begin{ctikz}
    \defc{a/0/1,b/1/1,c/1/0,d/2/0,e/0/0,f/2/1}
    \drawe{a/b}{c/d}{b/c,b/e,c/f}
    \drawv{a,b}{f}{c,d}{e}
  \end{ctikz}
  \begin{kase}
    \item Suppose $e \sim f$. Flip the cut preserver $\sset{b,c}$. Then $M$ contains multiple edges between $be$ and $cf$, and so a subset of $G(be, 2)$ is a negative tail reducer by \cref{lem:no-multi-neg-2}.
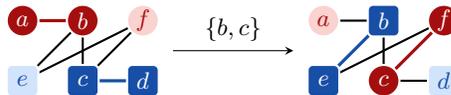
    \begin{ctikz}
      \begin{scope}
        \defc{a/0/1,b/1/1,c/1/0,d/2/0,e/0/0,f/2/1}
        \drawe{a/b}{c/d}{b/c,b/e,c/f,e/f}
        \drawv{a,b}{f}{c,d}{e}
      \end{scope}
      \draw[-stealth, thin] (2.5,0.5) -- (4.5,0.5) node[midway, above, font=\footnotesize] {$\sset{b,c}$};
      \begin{scope}[shift={(5,0)}]
        \defc{a/0/1,b/1/1,c/1/0,d/2/0,e/0/0,f/2/1}
        \drawe{c/f}{b/e}{a/b,c/d,b/c,e/f}
        \drawv{c,f}{a}{b,e}{d}
      \end{scope}
    \end{ctikz}
    \item Suppose some neighbor, say $g$, of $f$ outside $\sset{c,d}$ is blue. Let $h$ be the blue neighbor of $g$. Flip the cut preserver $\sset{b,c}$. Then $M$ contains the path $be, cf, gh$ as a subgraph, and so a subset of $G(cf, \paththreesubgraphneg)$ is a negative tail reducer by \cref{lem:path-3-subgraph-neg}. 
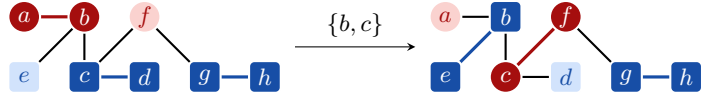
\label{1-neg}
    \begin{ctikz}
      \begin{scope}
        \defc{a/0/1,b/1/1,c/1/0,d/2/0,e/0/0,f/2/1,g/3/0,h/4/0}
        \drawe{a/b}{c/d,g/h}{b/c,b/e,c/f,f/g}
        \drawv{a,b}{f}{c,d,g,h}{e}
      \end{scope}
      \draw[-stealth, thin] (4.5,0.5) -- (6.5,0.5) node[midway, above, font=\footnotesize] {$\sset{b,c}$};
      \begin{scope}[shift={(7,0)}]
        \defc{a/0/1,b/1/1,c/1/0,d/2/0,e/0/0,f/2/1,g/3/0,h/4/0}
        \drawe{c/f}{b/e,g/h}{b/c,a/b,c/d,f/g}
        \drawv{c,f}{a}{b,e,g,h}{d}
      \end{scope}
    \end{ctikz}
    \item Suppose $e \nsim f$, every neighbor of $f$ outside $\sset{c,d}$ is bluish, $e \nsim a$ and $f \nsim d$. Clearly, both $a$ and $f$ have two bluish neighbors.
    \begin{kase}
      \item If $a$ and $f$ share a bluish neighbor, then a subset of $\sset{a,b,f}$ is a negative tail reducer by \cref{lem:pentagon-neg-versatile}.
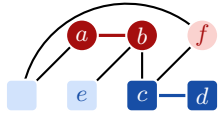
      \begin{ctikz}
        \defc{a/0/1,b/1/1,c/1/0,d/2/0,e/0/0,f/2/1,g/-1/0}
        \arcthreeone{g}
        \drawe{a/b}{c/d}{b/c,b/e,c/f,a/g}
        \drawv{a,b}{f}{c,d}{e}
        \drawvnl{}{}{}{g}
      \end{ctikz}
      \item If $a$ and $f$ share no bluish neighbors, then $\sset{a,b,f}$ is a negative tail reducer by \cref{fig:ntr-ag}.
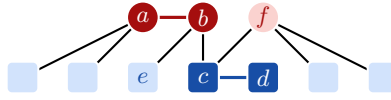
      \begin{ctikz}
        \defc{a/0/1,b/1/1,c/1/0,d/2/0,e/0/0,f/2/1,g/-1/0,h/-2/0,i/3/0,j/4/0}
        \drawe{a/b}{c/d}{b/c,b/e,c/f,a/g,a/h,f/i,f/j}
        \drawv{a,b}{f}{c,d}{e}
        \drawvnl{}{}{}{g,h,i,j}
      \end{ctikz}
    \end{kase}
    \item Suppose $e \nsim f$, every neighbor of $f$ outside $\sset{c,d}$ is bluish, and $e \sim a$ or $f \sim d$. We may assume that $e \sim a$ up to symmetry. Clearly, the other neighbors of $d$ are reddish. We may assume that the third neighbor of $e$ is reddish, because otherwise the current case reduces to \ref{1-neg}. Let $g$ be the third neighbor of $a$. Clearly, $g$ is bluish. 
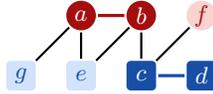
\label{item:path-2-case-4}
    \begin{ctikz}
      \defc{a/0/1,b/1/1,f/2/1,e/0/0,c/1/0,d/2/0,g/-1/0}
      \drawe{a/b}{c/d}{a/e,b/e,b/c,c/f,a/g}
      \drawv{a,b}{f}{c,d}{e,g}
    \end{ctikz}
    \begin{kase}
      \item Suppose $d$ and $e$ share a reddish neighbor. Depending on whether $f \sim d$, there are two possibilities. In each case, $\sset{d,e}$ is a negative tail reducer by \cref{fig:ntr-h}.
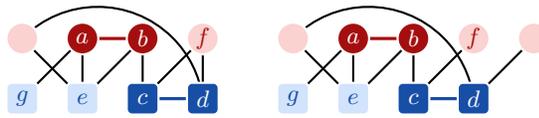
      \begin{ctikz}
        \begin{scope}
          \defc{a/0/1,b/1/1,f/2/1,e/0/0,c/1/0,d/2/0,g/-1/0,h/-1/1}
          \arcthreeminusone{h}
          \drawe{a/b}{c/d}{a/e,b/e,b/c,c/f,d/f,a/g,e/h}
          \drawv{a,b}{f}{c,d}{e,g}
          \drawvnl{}{h}{}{}
        \end{scope}
        \begin{scope}[shift={(4.5,0)}]
          \defc{a/0/1,b/1/1,f/2/1,e/0/0,c/1/0,d/2/0,g/-1/0,h/-1/1,i/3/1}
          \arcthreeminusone{h}
          \drawe{a/b}{c/d}{a/e,b/e,b/c,c/f,a/g,e/h,d/i}
          \drawv{a,b}{f}{c,d}{e,g}
          \drawvnl{}{h,i}{}{}
        \end{scope}
      \end{ctikz}
      \item Suppose $g$ has another red neighbor, say $h$. 
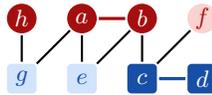
\label{4.2-neg}
      \begin{ctikz}
        \defc{a/0/1,b/1/1,f/2/1,e/0/0,c/1/0,d/2/0,g/-1/0,h/-1/1}
        \drawe{a/b}{c/d}{a/e,b/e,b/c,c/f,a/g,g/h}
        \drawv{a,b,h}{f}{c,d}{e,g}
      \end{ctikz}
      \begin{kase}
        \item Suppose $h$ has a blue neighbor, say $i$. Flip the cut preserver $\sset{h,i}$. Then $M$ contains the path $cd,ab,gh$ as a subgraph, and so a subset of $G(ab,\paththreesubgraphneg)$ is a negative tail reducer by \cref{lem:path-3-subgraph-neg}.
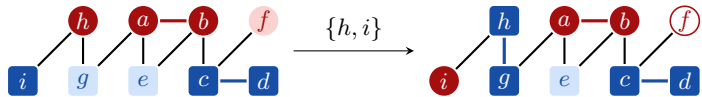
        \begin{ctikz}
          \begin{scope}
            \defc{a/0/1,b/1/1,f/2/1,e/0/0,c/1/0,d/2/0,g/-1/0,h/-1/1,i/-2/0}
            \drawe{a/b}{c/d}{a/e,b/e,b/c,c/f,a/g,g/h,h/i}
            \drawv{a,b,h}{f}{c,d,i}{e,g}
          \end{scope}
          \draw[-stealth, thin] (2.5,0.5) -- (4.5,0.5) node[midway, above, font=\footnotesize] {$\sset{h,i}$};
          \begin{scope}[shift={(7,0)}]
            \defc{a/0/1,b/1/1,f/2/1,e/0/0,c/1/0,d/2/0,g/-1/0,h/-1/1,i/-2/0}
            \drawe{a/b}{c/d,g/h}{a/e,b/e,b/c,c/f,a/g,h/i}
            \drawv{a,b,i}{}{c,d,g,h}{e}
            \node[a-vertex] at (f) {$f$};
          \end{scope}
        \end{ctikz}
        \item Suppose $h$ has no blue neighbor. Depending on whether $h \sim e$, there are two possibilities. In each case, $\sset{a,h}$ is a tail reducer by \cref{fig:ntr-b,fig:ntr-d}.
        \begin{ctikz}
          \begin{scope}
            \defc{a/0/1,b/1/1,f/2/1,e/0/0,c/1/0,d/2/0,g/-1/0,h/-1/1}
            \drawe{a/b}{c/d}{a/e,b/e,b/c,c/f,a/g,g/h,h/e}
            \drawv{a,b,h}{f}{c,d}{e,g}
          \end{scope}
          \begin{scope}[shift={(5.5,0)}]
            \defc{a/0/1,b/1/1,f/2/1,e/0/0,c/1/0,d/2/0,g/-1/0,h/-1/1,i/-2/0}
            \drawe{a/b}{c/d}{a/e,b/e,b/c,c/f,a/g,g/h,h/i}
            \drawv{a,b,h}{f}{c,d}{e,g}
            \drawvnl{}{}{}{i}
          \end{scope}
        \end{ctikz}
      \end{kase}
      \item Suppose $g \sim f$.
      \begin{kase}
        \item If $f \sim d$, then $\sset{a,f}$ is a negative tail reducer by \cref{fig:ntr-h}. \label{4.3-neg}
        \begin{ctikz}
          \defc{a/0/1,b/1/1,f/2/1,e/0/0,c/1/0,d/2/0,g/-1/0}
          \arcthreeone{g}
          \drawe{a/b}{c/d}{a/e,b/e,b/c,c/f,d/f,a/g}
          \drawv{a,b}{f}{c,d}{e,g}
        \end{ctikz}
        \item If $f \nsim d$, then $\sset{a,b,f}$ is a negative tail reducer by \cref{fig:ntr-s}.
        \begin{ctikz}
          \defc{a/0/1,b/1/1,f/2/1,e/0/0,c/1/0,d/2/0,g/-1/0,h/3/0}
          \arcthreeone{g}
          \drawe{a/b}{c/d}{a/e,b/e,b/c,c/f,h/f,a/g}
          \drawv{a,b}{f}{c,d}{e,g,h}
        \end{ctikz}
      \end{kase}
      \item Suppose $e$ and $d$ do not share reddish neighbors, $g$ has no other red neighbors, and $g \nsim f$.
      \begin{kase}
        \item If $g$ shares a reddish neighbor with $e$, then $\sset{c,e,g}$ is a negative tail reducer by \cref{fig:ntr-v,fig:degree-ntr-v} when $g$ is of degree $3$ or $2$, respectively. \label{4.4.1-neg}
        \begin{ctikz}
          \defc{a/0/1,b/1/1,f/2/1,e/0/0,c/1/0,d/2/0,g/-1/0,h/-1/1,l/-2/1}
          \drawe{a/b}{c/d}{a/e,b/e,b/c,c/f,a/g,e/h,g/h,g/l}
          \drawv{a,b}{f}{c,d}{e,g}
          \drawvnl{}{h,l}{}{}
        \end{ctikz}
        \item If $g$ shares no reddish neighbors with $e$, and $g$ is of degree at most $2$, then $\sset{c,e,g}$ is a negative tail reducer by \cref{fig:degree-ntr-x-two,fig:degree-ntr-x-one}. \label{4.4.2-neg}
        \begin{ctikz}
          \defc{a/0/1,b/1/1,f/2/1,e/0/0,c/1/0,d/2/0,g/-1/0,h/-1/1,l/-2/1}
          \drawe{a/b}{c/d}{a/e,b/e,b/c,c/f,a/g,e/h,g/l}
          \drawv{a,b}{f}{c,d}{e,g}
          \drawvnl{}{h,l}{}{}
        \end{ctikz}
        \item Suppose $g$ shares no reddish neighbors with $e$, $g$ is of degree $3$, and $g$ shares a reddish neighbor with $d$. \label{4.4.3-neg}
        \begin{kase}
          \item If $d \sim f$, then $\sset{c,d,e,g}$ is a negative tail reducer by \cref{fig:ntr-ah}.
          \begin{ctikz}
            \defc{a/0/1,b/1/1,f/2/1,e/0/0,c/1/0,d/2/0,g/-1/0,h/-1/1,i/3/1,j/-2/1}
            \arcfourone{g}
            \drawe{a/b}{c/d}{a/e,b/e,b/c,c/f,d/f,a/g,e/h,g/j,d/i}
            \drawv{a,b}{f}{c,d}{e,g}
            \drawvnl{}{h,i,j}{}{}
          \end{ctikz}
          \item If $d \nsim f$, and $d$ and $g$ share two reddish neighbors, then $\sset{c,d,e,g}$ is a negative tail reducer by \cref{fig:ntr-ai}.
          \begin{ctikz}
            \defc{a/0/1,b/1/1,f/2/1,e/0/0,c/1/0,d/2/0,g/-1/0,h/-1/1,i/-2/1,j/3/1}
            \arcfourminusone{i}
            \arcfourone{g}
            \drawe{a/b}{c/d}{a/e,b/e,b/c,c/f,a/g,e/h,g/i,d/j}
            \drawv{a,b}{f}{c,d}{e,g}
            \drawvnl{}{h,i,j}{}{}
          \end{ctikz}
          \item Suppose $d \nsim f$, and $d$ and $g$ share exactly one reddish neighbor, say $h$. We may assume that $h$ is of degree $3$ because otherwise $\sset{a,b,h}$ is a negative tail reducer by \cref{fig:degree-ntr-l}. Let $i$ be the third neighbor of $h$. 
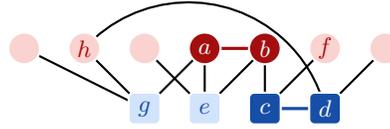
\label{4.4.3.3-neg}
          \begin{ctikz}
            \defc{a/0/1,b/1/1,f/2/1,e/0/0,c/1/0,d/2/0,g/-1/0,i/-1/1,h/-2/1,j/-3/1,k/3/1}
            \arcfourminusone{h}
            \drawe{a/b}{c/d}{a/e,b/e,b/c,c/f,a/g,e/i,g/h,g/j,d/k}
            \drawv{a,b}{f,h}{c,d}{e,g}
            \drawvnl{}{i,j,k}{}{}
          \end{ctikz}
          \begin{kase}
            \item Suppose $i$ is blue, and $i$ has a red neighbor, say $j$. Flip the cut preserver $\sset{i,j}$. Then $M$ contains the path $ab,cd,hi$ as a subgraph, and so a subset of $G(cd, \paththreesubgraphneg)$ is a negative tail reducer by \cref{lem:path-3-subgraph-neg}. \label{4.4.3.3.1-neg}
            \begin{ctikz}
              \defc{a/0/1,b/1/1,f/2/1,e/0/0,c/1/0,d/2/0,g/-1/0,x/-1/1,h/-2/1,y/-3/1,k/3/1,i/-2/0,j/-4/1}
              \arcfourminusone{h}
              \drawe{a/b}{c/d}{a/e,b/e,b/c,c/f,a/g,e/x,g/h,g/y,d/k,h/i,i/j}
              \drawv{a,b,j}{f,h}{c,d,i}{e,g}
              \drawvnl{}{y,k,x}{}{}
              \begin{scope}[shift={(8,-2)}]
                \draw[-stealth, thin] (-6.5,0.5) -- (-4.5,0.5) node[midway, above, font=\footnotesize] {$\sset{i,j}$};
                \defc{a/0/1,b/1/1,f/2/1,e/0/0,c/1/0,d/2/0,g/-1/0,x/-1/1,h/-2/1,y/-3/1,k/3/1,i/-2/0,j/-4/1}
                \arcfourminusone{h}
                \drawe{a/b,h/i}{c/d}{a/e,b/e,b/c,c/f,a/g,e/x,g/h,g/y,d/k,i/j}
                \drawv{a,b,i,h}{}{c,d,j}{e,g}
                \node[a-vertex] at (f) {$f$};
                \node[a-vertex] at (k) {};
                \node[a-vertex] at (x) {};
                \node[a-vertex] at (y) {};
              \end{scope}
            \end{ctikz}
            \item Suppose $i$ has no red neighbors. Depending on how the reddish neighbors of $i$ overlap with those of $c,d,e,g$, there are many possibilities. \label{4.4.3.3.2-neg}
            \begin{kase}
              \item If $i$ and $d$ share another reddish neighbor, then $\sset{d,i}$ is a negative tail reducer by \cref{fig:ntr-e,fig:degree-ntr-e}.
              \begin{ctikz}
                \defc{a/0/1,b/1/1,f/2/1,e/0/0,c/1/0,d/2/0,g/-1/0,x/-1/1,h/-2/1,y/-3/1,k/3/1,i/-2/0,z/-4/1}
                \arcfourminusone{h}
                \arcfiveone{i}
                \drawe{a/b}{c/d}{a/e,b/e,b/c,c/f,a/g,e/x,g/h,g/y,d/k,h/i,i/z}
                \drawv{a,b}{f,h}{c,d}{e,g}
                \node[b-vertex] at (i) {$i$};
                \drawvnl{}{y,k,x,z}{}{}
              \end{ctikz}
              \item Suppose $i$ and $g$ share another reddish neighbor. Depending on whether $i$ and $c$ share the reddish neighbor $f$, there are two possibilities. In each case, $\sset{c,g,i}$ is a negative tail reducer by \cref{fig:ntr-ac,fig:ntr-w,fig:degree-ntr-ac}.
              \begin{ctikz}
                \begin{scope}
                  \defc{a/0/1,b/1/1,f/2/1,e/0/0,c/1/0,d/2/0,g/-1/0,x/-1/1,h/-2/1,y/-3/1,k/3/1,i/-2/0,z/-4/1}
                  \arcfourminusone{h}
                  \drawe{a/b}{c/d}{a/e,b/e,b/c,c/f,a/g,e/x,g/h,g/y,d/k,h/i,i/z,i/y}
                  \drawv{a,b}{f,h}{c,d}{e,g}
                  \node[b-vertex] at (i) {$i$};
                  \drawvnl{}{y,k,x,z}{}{}
                \end{scope}
                \begin{scope}[shift={(7.5,0)}]
                  \defc{a/0/1,b/1/1,f/2/1,e/0/0,c/1/0,d/2/0,g/-1/0,x/-1/1,h/-2/1,y/-3/1,k/3/1,i/-2/0}
                  \arcfourminusone{h}
                  \arcfourone{i}
                  \drawe{a/b}{c/d}{a/e,b/e,b/c,c/f,a/g,e/x,g/h,g/y,d/k,h/i,i/y}
                  \drawv{a,b}{f,h}{c,d}{e,g}
                  \node[b-vertex] at (i) {$i$};
                  \drawvnl{}{k,x,y}{}{}
                \end{scope}
              \end{ctikz}
              \item If $i$ shares no other reddish neighbor with $d$ or $g$, and $i$ and $c$ share the reddish neighbor $f$, then $\sset{c,g,i}$ is a negative tail reducer by \cref{fig:ntr-ad,fig:degree-ntr-ad}.
              \begin{ctikz}
                \defc{a/0/1,b/1/1,f/2/1,e/0/0,c/1/0,d/2/0,g/-1/0,x/-1/1,h/-2/1,y/-3/1,k/3/1,i/-2/0,z/-4/1}
                \arcfourminusone{h}
                \arcfourone{i}
                \drawe{a/b}{c/d}{a/e,b/e,b/c,c/f,a/g,e/x,g/h,g/y,d/k,h/i,i/z}
                \drawv{a,b}{f,h}{c,d}{e,g}
                \node[b-vertex] at (i) {$i$};
                \drawvnl{}{k,x,y,z}{}{}
              \end{ctikz}
              \item If $i$ shares no other reddish neighbor with $d$ or $g$, and $i$ and $e$ share a reddish neighbor, then $\sset{d,e,i}$ is a negative tail reducer by \cref{fig:ntr-ae,fig:degree-ntr-ae}.
              \begin{ctikz}
                \defc{a/0/1,b/1/1,f/2/1,e/0/0,c/1/0,d/2/0,g/-1/0,x/-1/1,h/-2/1,y/-3/1,k/3/1,i/-2/0,z/-4/1}
                \arcfourminusone{h}
                \drawe{a/b}{c/d}{a/e,b/e,b/c,c/f,a/g,e/x,g/h,g/y,d/k,h/i,i/x,i/z}
                \drawv{a,b}{f,h}{c,d}{e,g}
                \node[b-vertex] at (i) {$i$};
                \drawvnl{}{k,x,y,z}{}{}
              \end{ctikz}
              \item If $i$ shares no reddish neighbors with $c,d,e,g$ except $h$, then $\sset{c,d,e,g,i}$ is a negative tail reducer by \cref{fig:ntr-ap,fig:degree-ntr-ap-two,fig:degree-ntr-ap-one}.
              \begin{ctikz}
                \defc{a/0/1,b/1/1,f/2/1,e/0/0,c/1/0,d/2/0,g/-1/0,x/-1/1,h/-2/1,y/-3/1,k/3/1,i/-2/0,z/-4/1,p/-5/1}
                \arcfourminusone{h}
                \drawe{a/b}{c/d}{a/e,b/e,b/c,c/f,a/g,e/x,g/h,g/y,d/k,h/i,i/p,i/z}
                \drawv{a,b}{f,h}{c,d}{e,g}
                \node[b-vertex] at (i) {$i$};
                \drawvnl{}{p,k,x,y,z}{}{}
              \end{ctikz}
            \end{kase}
            \item Suppose $i$ is bluish, and $i$ has a red neighbor, say $j$.
            \begin{ctikz}
              \defc{a/0/1,b/1/1,f/2/1,e/0/0,c/1/0,d/2/0,g/-1/0,x/-1/1,h/-2/1,y/-3/1,k/3/1,i/-4/0,j/-4/1}
              \arcfourminusone{h}
              \drawe{a/b}{c/d}{a/e,b/e,b/c,c/f,a/g,e/x,g/h,g/y,d/k,h/i,i/j}
              \drawv{a,b,j}{f,h}{c,d}{e,g,i}
              \drawvnl{}{y,k,x}{}{}
            \end{ctikz}
            \begin{kase}
              \item Suppose $j$ has a blue neighbor, say $k$. \ce{b} implies that $f \nsim k$ because otherwise $\sset{b,c,f,j,k}$ is a cut enhancer. Flip the cut preserver $\sset{j,k}$. \label{4.4.3.3.3.1-neg}
              \begin{ctikz}
                \begin{scope}
                  \defc{a/0/1,b/1/1,f/2/1,e/0/0,c/1/0,d/2/0,g/-1/0,x/-1/1,h/-2/1,y/-3/1,w/3/1,i/-4/0,j/-4/1,k/-5/0}
                  \arcfourminusone{h}
                  \drawe{a/b}{c/d}{a/e,b/e,b/c,c/f,a/g,e/x,g/h,g/y,d/w,h/i,i/j,j/k}
                  \drawv{a,b,j}{f,h}{c,d,k}{e,g,i}
                  \drawvnl{}{y,w,x}{}{}
                \end{scope}
                \begin{scope}[shift={(6,-2)}]
                  \draw[-stealth, thin] (-7.5,.5) -- (-5.5,.5) node[midway, above, font=\footnotesize] {$\sset{j,k}$};
                  \defc{a/0/1,b/1/1,f/2/1,e/0/0,c/1/0,d/2/0,g/-1/0,x/-1/1,h/-2/1,y/-3/1,w/3/1,i/-4/0,j/-4/1,k/-5/0}
                  \arcfourminusone{h}
                  \drawe{a/b}{c/d,i/j}{a/e,b/e,b/c,c/f,a/g,e/x,g/h,g/y,d/w,h/i,j/k}
                  \drawv{a,b,k}{f,h}{c,d,j,i}{e,g}
                  \node[a-vertex] at (w) {};
                  \node[a-vertex] at (x) {};
                  \node[a-vertex] at (y) {};
                \end{scope}
              \end{ctikz}
              \begin{kase}
                \item If $k$ and $d$ share a neighbor, say $l$, then $M$ contains the path $ab,cd,kl$ as a subgraph, and so a subset of $G(cd,\paththreesubgraphneg)$ is a negative tail reducer by \cref{lem:path-3-subgraph-neg}.
                \item If $k$ and $e$ share a neighbor, then $e$ has a red neighbor, and so the current case reduces to \ref{1-neg}.
                \item If $k$ and $g$ share a neighbor, then $g$ has a red neighbor, and so the current case reduces to \ref{4.2-neg}.
                \item If $k$ shares no neighbor with $d, e, g$, then the current case reduces to \ref{4.4.3.3.1-neg} or \ref{4.4.3.3.2-neg}, where $i$ is blue.
              \end{kase}
              \item Suppose $j$ has no blue neighbors. Let $k$ be the red neighbor of $j$. \label{4.4.3.3.3.2-neg}
              \begin{ctikz}
                \defc{a/0/1,b/1/1,f/2/1,e/0/0,c/1/0,d/2/0,g/-1/0,x/-1/1,h/-2/1,y/-3/1,z/3/1,i/-4/0,j/-4/1,k/-5/1}
                \arcfourminusone{h}
                \drawe{a/b,j/k}{c/d}{a/e,b/e,b/c,c/f,a/g,e/x,g/h,g/y,d/z,h/i,i/j}
                \drawv{a,b,j,k}{f,h}{c,d}{e,g,i}
                \drawvnl{}{x,y,z}{}{}
              \end{ctikz}
              We may assume that $k$ has a blue neighbor because otherwise $jk$ is isolated in $M$, and so a subset of $G(jk, 4)$ is a negative tail reducer by \cref{lem:isolated-neg-4}. Let $l$ be a blue neighbor of $k$. We may assume that $k \nsim i$ because otherwise the current case reduces to \ref{4.4.3.3.3.1-neg}. Let $m$ be the third neighbor of $j$, and let $n$ be the blue neighbor of $l$. Clearly, $m$ is bluish. We may assume that $n \nsim k$ because otherwise $\sset{k}$ is a negative tail reducer by \cref{fig:ntr-a}.
              \begin{ctikz}
                \defc{a/0/1,b/1/1,f/2/1,e/0/0,c/1/0,d/2/0,g/-1/0,x/-1/1,h/-2/1,y/-3/1,z/3/1,i/-4/0,j/-4/1,k/-5/1,l/-5/0,m/-3/0,n/-6/0}
                \arcfourminusone{h}
                \drawe{a/b,j/k}{c/d,l/n}{a/e,b/e,b/c,c/f,a/g,e/x,g/h,g/y,d/z,h/i,i/j,k/l,j/m}
                \drawv{a,b,j,k}{f,h}{c,d,l,n}{e,g,i,m}
                \drawvnl{}{x,y,z}{}{}
              \end{ctikz}
              Flip the cut preserver $\sset{k,l}$. Using the argument in \labelcref{4.4.3.3.3.1-neg}, we may assume that $l$ shares no neighbor with $c,d,e,g$.\footnote{If $l$ and $c$ share a neighbor, then their common neighbor is $f$, and $\sset{b,c,f,k,l}$ is a cut enhancer by \ce{b} before the flip. If $l$ and $d$ share a neighbor, say $u$, then $lu$ is red after the flip and $ab,cd,lu$ is a path in $M$, so \cref{lem:path-3-subgraph-neg} applies. If $l$ shares a neighbor with $e$ or $g$, then $e$ or $g$ has a red neighbor after the flip, reducing to \ref{1-neg} or \ref{4.2-neg}, respectively.}
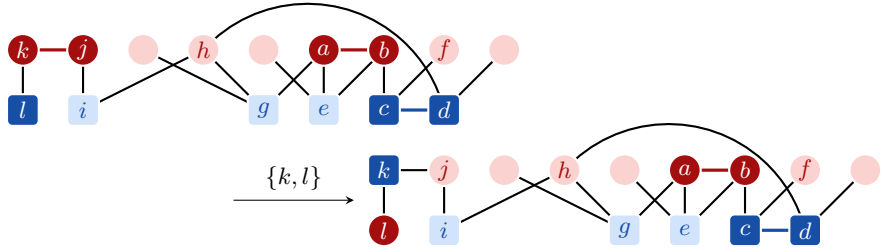
              \begin{ctikz}
                \node at (-9,.5) {$\cdots$};
                \draw[-stealth, thin] (-8.5,.5) -- (-6.5,.5) node[midway,above,font=\footnotesize] {$\sset{k,l}$};
                \defc{a/0/1,b/1/1,f/2/1,e/0/0,c/1/0,d/2/0,g/-1/0,x/-1/1,h/-2/1,y/-3/1,z/3/1,i/-4/0,j/-4/1,k/-5/0,l/-5/1,m/-3/0,n/-6/0}
                \arcfourminusone{h}
                \drawe{a/b}{c/d}{a/e,b/e,b/c,c/f,a/g,e/x,g/h,g/y,d/z,h/i,i/j,k/l,j/k,j/m,l/n}
                \drawv{a,b,l}{j,f,h}{c,d,k}{e,g,i,n}
                \node[b-vertex] at (m) {$m$};
                \drawvnl{}{x,y,z}{}{}
              \end{ctikz}
              We may assume that $i$ still has a red neighbor, because otherwise the current case reduces to \ref{4.4.3.3.2-neg}. Let $o$ be the red neighbor of $i$. By a symmetric argument, we may assume that $o$ has no blue neighbor.
              \begin{kase}
                \item Suppose $o \nsim l$. Let $p$ be the red neighbor of $o$. \ce{a} implies that $p \nsim k$ because otherwise $\sset{k,l,p}$ is a cut enhancer. We may assume that $p$ has a blue neighbor because otherwise $op$ is isolated in $M$, and so a subset of $G(op, 4)$ is a negative tail reducer by \cref{lem:isolated-neg-4}. Let $q$ be a blue neighbor of $p$. Then $q \nsim o$ because otherwise $\sset{q}$ is a negative tail reducer by \cref{fig:ntr-a}.
                \begin{ctikz}
                  \defc{a/0/1,b/1/1,f/2/1,e/0/0,c/1/0,d/2/0,g/-1/0,x/-1/1,h/-2/1,y/-3/1,z/3/1,i/-4/0,j/-4/1,k/-5/0,l/-5/1,m/-3/0,n/-6/0,o/-6/1,p/-7/1,q/-7/0}
                  \arcfourminusone{h}
                  \drawe{a/b,o/p}{c/d}{a/e,b/e,b/c,c/f,a/g,e/x,g/h,g/y,d/z,h/i,i/j,k/l,j/k,j/m,l/n,p/q,i/o}
                  \drawv{a,b,l,o,p}{j,f,h}{c,d,k,q}{e,g,i,n}
                  \node[b-vertex] at (m) {$m$};
                  \drawvnl{}{x,y,z}{}{}
                \end{ctikz}
                Flip the cut preserver $\sset{p,q}$. Using the argument in \ref{4.4.3.3.3.1-neg}, we may assume that $q$ shares no neighbor with $c,d,e,g$. Then the current case reduces to \ref{4.4.3.3.2-neg}.
                \item If $m \nsim k$ and $o \sim l$, then a subset of $\sset{j,l,o}$ is a negative tail reducer by \cref{lem:pentagon-neg-versatile}.
                \begin{ctikz}
                  \defc{a/0/1,b/1/1,f/2/1,e/0/0,c/1/0,d/2/0,g/-1/0,x/-1/1,h/-2/1,y/-3/1,z/3/1,i/-4/0,j/-4/1,k/-5/0,l/-5/1,m/-3/0,n/-6/0,o/-6/1}
                  \arcfourminusone{h}
                  \drawe{a/b,l/o}{c/d}{a/e,b/e,b/c,c/f,a/g,e/x,g/h,g/y,d/z,h/i,i/j,k/l,j/k,j/m,l/n,i/o}
                  \drawv{a,b,l,o}{j,f,h}{c,d,k}{e,g,i,m,n}
                  \drawvnl{}{x,y,z}{}{}
                \end{ctikz}
                \item If $m \sim k$ and $o \sim l$, then $\sset{a,b,h,j}$ is a negative tail reducer by \cref{fig:ntr-aj}.
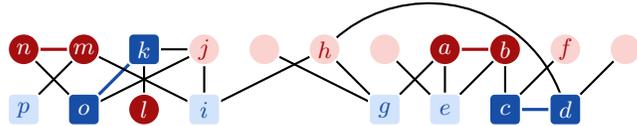
                \begin{ctikz}
                  \defc{a/0/1,b/1/1,f/2/1,e/0/0,c/1/0,d/2/0,g/-1/0,x/-1/1,h/-2/1,y/-3/1,z/3/1,i/-4/0,j/-4/1,k/-5/0,l/-5/1,m/-3/0,n/-6/0,o/-6/1}
                  \arcfourminusone{h}
                  \drawe{a/b,l/o}{c/d}{a/e,b/e,b/c,c/f,a/g,e/x,g/h,g/y,d/z,h/i,i/j,k/l,j/k,j/m,l/n,i/o}
                  \draw[blue-edge] (k) arc [start angle=233.968, end angle=306.032, radius=1.7];
                  \drawv{a,b,l,o}{j,f,h}{c,d,k,m}{e,g,i,n}
                  \drawvnl{}{x,y,z}{}{}
                \end{ctikz}
              \end{kase}
            \end{kase}
          \end{kase}
        \end{kase}
        \item Suppose the reddish neighbors of $d,e,g$ do not overlap, and $g$ is of degree $3$.
        \begin{kase}
          \item Suppose $d \sim f$. Let $h$ be the third neighbor of $f$. Recall from \ref{item:path-2-case-4} that $e \nsim f$, the third neighbor of $e$ is reddish, and $h$ is bluish. Temporarily flip the cut preserver $\sset{b,c}$. We may assume, after the temporary flip, that $h$ has two reddish neighbors, $h \nsim a$, and the reddish neighbors of $d,e,h$ do not overlap, because otherwise the current case reduces to \ref{4.2-neg}, \ref{4.3-neg} or \labelcref{4.4.1-neg,4.4.2-neg,4.4.3-neg}.\footnote{Under the relabeling $(a,b,c,d,e,f,g)\mapsto(f,c,b,e,d,a,h)$, the configuration after the flip is that of \ref{item:path-2-case-4}. Thus \ref{4.2-neg} handles the case that $h$ has another red neighbor, and \ref{4.3-neg} handles $h\sim a$ and $a\sim e$, the latter of which holds already. Finally, \ref{4.4.1-neg}, \ref{4.4.2-neg}, and \ref{4.4.3-neg} handle, respectively, a reddish neighbor shared by $h$ and $d$, the case where $h$ is of degree at most $2$, and a reddish neighbor shared by $h$ and $e$.}
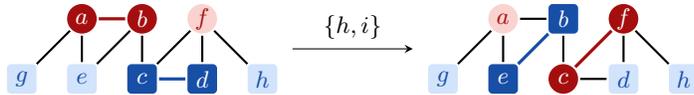
          \begin{ctikz}
            \begin{scope}
              \defc{a/0/1,b/1/1,f/2/1,e/0/0,c/1/0,d/2/0,g/-1/0,h/3/0}
              \drawe{a/b}{c/d}{a/e,b/e,b/c,c/f,d/f,a/g,f/h}
              \drawv{a,b}{f}{c,d}{e,g,h}
            \end{scope}
            \draw[-stealth, thin] (3.5,0.5) -- (5.5,0.5) node[midway, above, font=\footnotesize] {$\sset{b,c}$};
            \begin{scope}[shift={(7,0)}]
              \defc{a/0/1,b/1/1,f/2/1,e/0/0,c/1/0,d/2/0,g/-1/0,h/3/0}
              \drawe{c/f}{e/b}{a/b,a/e,b/c,c/d,d/f,a/g,f/h}
              \drawv{f,c}{a}{b,e}{d,g,h}
            \end{scope}
          \end{ctikz}
          Translating these assumptions back to the configuration before the temporary flip, we may assume that $h$ has $2$ other reddish neighbors, and the reddish neighbors of $d,e,h$ other than $f$ do not overlap.
          \begin{kase}
            \item If $g$ and $h$ share reddish neighbors, then $\sset{c,g,h}$ is a negative tail reducer by \cref{fig:ntr-w,fig:ntr-ad}.
            \begin{ctikz}
              \begin{scope}
                \defc{a/0/1,b/1/1,f/2/1,e/0/0,c/1/0,d/2/0,i/-1/1,h/3/0,g/-1/0,j/3/1,k/-2/1,l/4/1}
                \arcfiveminusone{k}
                \arcfiveone{g}
                \drawe{a/b}{c/d}{a/e,b/e,b/c,c/f,d/f,e/i,f/h,a/g,d/j,k/g,h/l}
                \drawv{a,b}{f}{c,d}{e,h,g}
                \drawvnl{}{i,j,k,l}{}{}
              \end{scope}
              \begin{scope}[shift={(8.5,0)}]
                \defc{a/0/1,b/1/1,f/2/1,e/0/0,c/1/0,d/2/0,i/-1/1,h/3/0,g/-1/0,j/3/1,k/-2/1,l/4/1,m/-3/1}
                \arcfiveminusone{k}
                \drawe{a/b}{c/d}{a/e,b/e,b/c,c/f,d/f,e/i,f/h,a/g,d/j,k/g,h/l,g/m}
                \drawv{a,b}{f}{c,d}{e,h,g}
                \drawvnl{}{i,j,k,l,m}{}{}
              \end{scope}
            \end{ctikz}
            \item If $g$ and $h$ share no reddish neighbors, then $\sset{c,d,e,g,h}$ is a negative tail reducer by \cref{fig:ntr-ao}.
            \begin{ctikz}
              \defc{a/0/1,b/1/1,f/2/1,e/0/0,c/1/0,d/2/0,i/-1/1,h/3/0,g/-1/0,j/3/1,k/-2/1,l/4/1,m/-3/1,n/5/1}
                \drawe{a/b}{c/d}{a/e,b/e,b/c,c/f,d/f,e/i,f/h,a/g,d/j,k/g,h/l,g/m,h/n}
                \drawv{a,b}{f}{c,d}{e,h,g}
                \drawvnl{}{i,j,k,l,m,n}{}{}
            \end{ctikz}
          \end{kase}
          \item If $d \nsim f$, then $\sset{c,d,e,g}$ is a negative tail reducer by \cref{fig:ntr-am}.
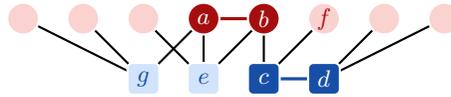
          \begin{ctikz}
            \defc{a/0/1,b/1/1,f/2/1,e/0/0,c/1/0,d/2/0,g/-1/0,h/-1/1,i/-2/1,j/-3/1,k/3/1,l/4/1}
            \drawe{a/b}{c/d}{a/e,b/e,b/c,c/f,a/g,e/h,g/i,g/j,d/k,d/l}
            \drawv{a,b}{f}{c,d}{e,g}
            \drawvnl{}{h,i,j,k,l}{}{}
          \end{ctikz}
        \end{kase}
      \end{kase}
    \end{kase}
  \end{kase}

  Finally notice that the distance between a vertex of $G(cf,\paththreesubgraphneg)$ and $ab$ is at most $\pathtwoneg$ in \ref{1-neg}, and it is the furthest distance possible between any vertex in a cut preserver or a negative tail reducer and $\alpha \in \sset{ab,cd}$.
\end{proof}

\begin{proof}[Proof of \cref{lem:neg-tail-14}]
  Let $\alpha$ be a vertex of $M$.
  \begin{kase}
    \item If $M(\alpha,4)$ contains multiple edges, then a subset of $G(\alpha,10)$ is a negative tail reducer by \cref{lem:no-multi-neg-2}.
    \item If $M(\alpha,4)$ is a simple graph of maximum degree at least $3$, then a subset of $G(\alpha,10)$ is a negative tail reducer by \cref{lem:m-deg-2-neg-2}.
    \item If $M(\alpha,4)$ is a simple graph of maximum degree at most $2$, then $M(\alpha,4)$ must be an isolated vertex, a path (on at least $2$ vertices), or an even cycle (on at least $4$ vertices), and so a subset of $G(\alpha,\pathtwoneg)$ is a negative tail reducer by \cref{lem:isolated-neg-4,lem:path-5-neg-8,lem:cycle-4-neg-10,lem:path-4-neg,lem:path-3-neg,lem:path-2-neg}. \qedhere
  \end{kase}
\end{proof}

\section{Positive fraction} \label{sec:pos-frac}

To prove the positive fraction result \cref{thm:main2} for connected subcubic graphs, we follow Mohar's strategy in \cite{M16} to combine `local' tail reducers.

To make use of \cref{lem:mohar} when the graph is not bipartite, we introduce the following result so that \cref{lem:mohar} is applicable as long as the graph is `locally' bipartite.

\begin{lemma} \label{lem:bipartite-relax}
  For every vertex subset $A$ of a graph $G$, and every induced subgraph $H$ of $G$, if all the neighbors of $V(H) \setminus A$ in $G$ are in $A$, then
  \[
    t_H^\pm(A') - t_H^\pm(A' \setminus C) = t_G^\pm(A) - t_G^\pm(A \setminus C)
  \]
  for every vertex subset $C$ of $H$, where $A' = V(H) \cap A$.
\end{lemma}

\begin{proof}
  We shall only prove the statement for the positive tail estimate. The proof for the negative tail estimate is similar. Set $B = V(G) \setminus A$ and $B' = V(H) \setminus A$, and fix a vertex subset $C$ of $H$. Based on the definition of positive tail estimates, it suffices to show that
  \[
    t_{H[B' \cup C]}^+ - t_{H[B']}^+ = t_{G[B \cup C]}^+ - t_{G[B]}^+.
  \]
  Since all the neighbors of $V(H) \setminus A$ in $G$ are in $A$, all the connected components of $H[B']$ are trivial, and moreover the nontrivial connected components of $H[B' \cup C]$, denoted by $H_1, \dots, H_k$, are precisely the nontrivial connected components of $G[B \cup C]$ that have not appeared in $G[B]$. Therefore both sides of the equation are equal to $t^+_{H_1} + \dots + t^+_{H_k}$.
\end{proof}

The next result lets us combine `local' tail reducers.

\begin{lemma} \label{lem:local-to-global}
  For every vertex subset $A$ of a graph $G$, and every family of vertex subsets $C_1, \dots, C_k$ of $G$, if each connected component of $G-A$ has diameter at most $d$, and $C_1, \dots, C_k$ are pairwise at distance at least $d + 3$, then \[
    t^\pm_G(A) - t^\pm_G\left(A \setminus \bigcup_{i=1}^k C_i\right) = \sum_{i = 1}^k t^\pm_G(A) - t^\pm_G(A \setminus C_i).
  \]
\end{lemma}

\begin{proof}
  We shall only prove the statement for the positive tail estimate. Set $C = C_1 \cup \dots \cup C_k$, and $B = V(G) \setminus A$. It suffices to show that
  \[
    t^+_{G[B \cup C]} - t^+_{G[B]} = \sum_{i=1}^k t^+_{G[B \cup C_i]} - t^+_{G[B]}.
  \]

  Let $\mathcal{H}$ (respectively, $\mathcal{H}'$) be the set of the connected components of $G[B]$ (respectively, $G[B \cup C]$) that are within distance $1$ from $C$ in $G$. Note that
  \[
    t^+_{G[B \cup C]} - t^+_{G[B]} = \abs{A \cap C} + \sum_{H \in \mathcal{H}} t^+_H - \sum_{H' \in \mathcal{H}'}t^+_{H'}.
  \]
  For each $i \in \sset{1, \dots, k}$, let $\mathcal{H}_i$ (respectively, $\mathcal{H}_i'$) be the set of the connected components of $G[B]$ (respectively, $G[B \cup C_i]$) that are within distance $1$ from $C_i$ in $G$. Note that
  \[
    t^\pm_G(A) - t^\pm_G(A \setminus C_i) = \abs{A \cap C_i} + t^+_{G[B]} - t^+_{G[B \cup C_i]} = \abs{A \cap C_i} + \sum_{H \in {\mathcal{H}_i}} t^+_H - \sum_{H' \in \mathcal{H}_i'}t^+_{H'}.
  \]
  Finally, since $C_1, \dots, C_k$ are pairwise at distance at least $d + 3$, we know that $\abs{A \cap C} = \sum_{i=1}^k \abs{A \cap C_i}$, and moreover $\mathcal{H}_1, \dots, \mathcal{H}_k$ (respectively, $\mathcal{H}_1', \dots, \mathcal{H}_k'$) form a partition of $\mathcal{H}$ (respectively, $\mathcal{H}'$).
\end{proof}

We finally have all the ingredients for the proof of our second main theorem.

\begin{proof}[Proof of \cref{thm:main2}]
  With hindsight, we choose $\eta = 2^{-56}$. Let $G$ be a connected subcubic graph on $n$ vertices, and let $(A, B)$ be a maximum cut of $G$.

  We may assume that no subgraph of $G$ can be isomorphic to the Heawood graph. Indeed, suppose for a moment that a subgraph of $G$ is isomorphic to the Heawood graph. Since $G$ is connected subcubic and the Heawood graph is $3$-regular, it must be the case that $G$ itself is isomorphic to the Heawood graph. In particular, $n = 14$, in which case the theorem is vacuously true because $[\tfrac12 (1-\eta)(n+1), \tfrac12 (1+\eta)(n+1)]$ contains no integers.

  \begin{claim} \label{claim:1}
    There exist vertices $v_1, \dots, v_j$ of $G$ and edges $\alpha_{j+1}, \dots, \alpha_k$ of $G[A] \cup G[B]$ such that $v_1, \dots, v_j, \alpha_{j+1}, \dots, \alpha_k$ are pairwise at distance at least $38$ in $G$, each $G(v_i, 18)$ is a subgraph of the bipartite graph $G[A, B]$ for $i \in \sset{1, \dots, j}$, and $k \ge \eta n$.
  \end{claim}

  \begin{claimproof}[Proof of \cref{claim:1}]
    Starting with $V_0 = V(G)$, recursively we pick an edge $\alpha_i$ (if any) of $G[A \cap V_{i-1}] \cup G[B \cap V_{i-1}]$, and then we obtain $V_i$ from $V_{i-1}$ by removing the vertices within distance $37$ from $\alpha_i$ in $G$. Let $\alpha_1, \dots, \alpha_j$ be the edges we have picked. Set $U := V(G) \setminus V_{j}$, and let $U'$ be the set of vertices within distance $17$ from $U$ in $G$.

    Starting with $W_0 = V(G) \setminus U'$, we now recursively pick a vertex $v_i$ of $G[W_{i-1}]$, and then we obtain $W_i$ from $W_{i-1}$ by removing the vertices within distance $38$ from $v_i$ in $G$. Let $v_1, \dots, v_k$ be the vertices we have picked.

    Clearly, $\alpha_1, \dots, \alpha_j, v_1, \dots, v_k$ are pairwise at distance at least $37$ in $G$. Since $G-U$ is a subgraph of the bipartite graph $G[A, B]$, and each $v_i$ is at distance more than $18$ from $U$, we know that each $G(v_i, 17)$ is a subgraph of $G[A, B]$.

    Finally, we are left to prove that $j + k \ge \eta n$. Because $\abs{U'} \le \abs{G(\alpha_1, 37+17)} + \dots + \abs{G(\alpha_j, 37+17)} \le 2^{56}j$, and $\abs{W_0} \le \abs{G(v_1, 37)} + \dots + \abs{G(v_k, 37)} \le 2^{39}k$, we know that $n = \abs{U'} + \abs{W_0} \le 2^{56}j + 2^{39}k$, and so $j + k \ge 2^{-56}n$.
  \end{claimproof}

  \begin{claim} \label{claim:2}
    There exist vertex subsets $C_1, \dots, C_k$ of $G$ such that each $C_i$ decreases $t^+_G(A)$ or $t^+_G(B)$, and $C_1, \dots, C_k$ are pairwise at distance at least $4$ in $G$.
  \end{claim}

  \begin{claimproof}[Proof of \cref{claim:2}]
    Let $v_1, \dots, v_j, \alpha_{j+1}, \dots, \alpha_k$ be the vertices and edges as in \cref{claim:1}.

    For $i \in \sset{1,\dots,j}$, since $G(v_i, 17)$ is not the Heawood graph, we can apply \cref{lem:mohar} to the graph $G(v_i, 17)$ and the vertex $v_i$ to obtain a vertex subset $C_i$ of $G(v_i, 17)$ such that $C_i$ decreases $t^+_{G(v_i, 17)}(A)$ or $t^+_{G(v_i, 17)}(B)$. Since $G(v_i, 18)$ is bipartite, we know from \cref{lem:bipartite-relax} that $C_i$ also decreases $t^+_G(A)$ or $t^+_G(B)$.

    For $i \in \sset{j+1, \dots, k}$, since $(A, B)$ is a maximum cut, we can apply \cref{lem:pos-tail-16} to the graph $G$, and the edge $\alpha_i$ to obtain a positive tail reducer $C_i$ that is a subset of $G(\alpha_i, \abstail)$.

    Since $v_1, \dots, v_j, \alpha_{j+1}, \dots, \alpha_k$ are pairwise at distance at least $38$, the vertex subsets $C_1, \dots, C_k$ are pairwise at distance at least $4$ in $G$.
  \end{claimproof}

  Let $C_1, \dots, C_k$ be the vertex subsets as in \cref{claim:2}. Set $I_A := \dset{i}{C_i \text{ decreases }t^+_G(A)}$ and $I_B = \dset{i}{C_i \text{ decreases }t^+_G(B)}$. From \cref{lem:max-cut}, we know that each connected component of $G[A] \cup G[B]$ has diameter at most $1$, and moreover $t^+_G(A) = \abs{A}$ and $t^+_G(B) = \abs{B}$. We can apply \cref{lem:local-to-global} separately to $\dset{C_i}{i \in I_A}$ and $\dset{C_i}{i \in I_B}$ to obtain
  \[
    t_G^+\left(A \setminus \bigcup_{i \in I_A} C_i\right) \le \abs{A} - \abs{I_A} \text{ and }t_G^+\left(B \setminus \bigcup_{i \in I_B} C_i\right) \le \abs{B} - \abs{I_B},
  \]
  which implies via \cref{lem:key} that \[
    t^+_G \le \tfrac12 (\abs{A} - \abs{I_A} + \abs{B} - \abs{I_B}) \le \tfrac12 (n - k) \le \tfrac12 (1-\eta)n < \tfrac12 (1-\eta)(n+1).
  \]

  Finally, one can also show that $t^-_G < \tfrac12 (1-\eta)(n+1)$ through a symmetric argument.
\end{proof}

\section{Further remarks}

In this paper, we showed that the median eigenvalues of every connected graph of maximum degree at most three, except for the Heawood graph, are at most $1$ in absolute value. This result contributes to a broader line of research on bounding the median eigenvalues of graphs of bounded maximum degree.

A significant related result is due to Mohar and Tayfeh-Rezaie~\cite{MT15}, who proved that for every $d \ge 3$, the median eigenvalues of every connected \emph{bipartite} graph $G$ of maximum degree $d$ are at most $\sqrt{d-2}$ in absolute value, unless $G$ is the incidence graph of a projective plane of order $d-1$, in which case the median eigenvalues are equal to $\pm\sqrt{d-1}$. This generalizes an earlier result of Mohar~\cite{M16} and raises the question of optimality of the bound $\sqrt{d-2}$.

\begin{problem} \label{prob:infinite}
  For every $d \ge 3$ and every $\eps > 0$, are there infinitely many connected bipartite graphs $G$ of maximum degree $d$ whose median eigenvalues have absolute value greater than $\sqrt{d-2} - \eps$?
\end{problem}

For $d = 3$, this question was resolved positively by a construction of Guo and Mohar~\cite{GM14}: they constructed an infinite family of connected bipartite cubic graphs, known as the Guo--Mohar graphs, whose median eigenvalues are $-1$ and $+1$. Recent work \cite{GR26,HJ26} gives a classification of all connected subcubic graphs whose median eigenvalues satisfy $\lambda_\ell = -1$ and $\lambda_h = +1$: they are the Guo--Mohar graphs and nine of the thirteen sporadic graphs from \cite{GR26}, together with the Huang--Jiang graphs and seven of the eight sporadic graphs from \cite{HJ26}. We point out that the Koll\'ar--Sarnak graphs, first introduced in \cite{KS21} and identified in \cite{GR26} as the other infinite family of connected cubic graphs with no eigenvalues in $(-1,1)$, have both median eigenvalues equal to $-1$. We refer the reader to \cite[Section~6.3]{HJ26} for further discussion.

However, \cref{prob:infinite} remains wide open for $d \ge 4$. Mohar and Tayfeh-Rezaie also observed that the incidence graph of a symmetric $2$-design with parameters $(v,d,2)$, also known as a biplane of order $d-2$, would provide a $d$-regular bipartite graph of order $2v$ whose median eigenvalues are exactly $\pm \sqrt{d-2}$, where $v = d(d-1)/2+1$. The existence of infinitely many biplanes is a long-standing open problem in design theory (see \cite[Chapter~15.8]{H86}).

If we remove the bipartite constraint, earlier work by Mohar \cite{M15} established that for every $d \ge 3$, the median eigenvalues of every connected graph of maximum degree $d$ are at most $\sqrt{d}$ in absolute value, and specifically for $d = 3$, this bound can be improved to $\sqrt{d-1}$. We reiterate the conjecture that the same improvement should hold for every $d \ge 4$.

Very recently, Acharya, Jiang and Zhang \cite{AJZ26} resolved Mohar's conjecture for all but finitely many values of $d$. In particular, they proved that the median eigenvalues of any graph of maximum degree $d$ are bounded above by $\sqrt{d-1}$. They also proved the corresponding lower bound $-\sqrt{d-1}$ in three cases: when the graph is triangle-free, when $d-1$ is a perfect square, or when $d \ge 75$.

\section*{Acknowledgements} We thank Xingxing Yu for inspirational discussions on planar subcubic graphs, and are grateful to Patrick Fowler, Gordon Royle, and Shengtong Zhang for their helpful remarks on an earlier draft of this paper. The third author is grateful to Boris Bukh for granting him access to OpenAI's GPT-5.6 Sol model in Codex, which was used in the Lean formalization.

\bibliographystyle{plain}
\bibliography{subcubic}

\appendix
\crefalias{section}{appendix}

\section{Computer-assisted proofs} \label{sec:app}

\begin{proof}[Proof of \cref{prop:ptr,prop:ntr}]
  Recall Sylvester's law of inertia: if $A$ is a symmetric matrix, for any invertible matrix $S$, the number of positive eigenvalues, also known as \emph{the positive index of inertia}, of $D = SAS^\intercal$ is constant. This result is particularly useful when $D$ is diagonal, as the positive index of inertia is equal to the number of positive diagonal entries of $D$.

  Note that the positive and negative tails of $G$ are equal to the positive index of inertia of $A_G - I$ and $- A_G - I$ respectively. Therefore the core of our algorithm needs to compute the positive index of inertia of a symmetric matrix $A$ with rational entries. Given a symmetric matrix $A$ of order $n$ with rational entries, we compute its positive index of inertia $p$ through the $LDL^\intercal$ decomposition.

  The algorithm initializes a counter $p$ to track positive diagonal elements and processes each row $i$ from $1$ to $n$. For each row $i$, it identifies the first non-zero element $A_{ij}$. If no such $j$ exists, the algorithm continues to the next row. Otherwise, it ensures the diagonal element $A_{ii}$ is non-zero by adding row and column $j$ to row and column $i$ if necessary. If $A_{ii}$ is still zero after these operations, which happens if and only if $A_{ii} = 0$ and $A_{jj} = -2A_{ij}$, the algorithm swaps row $i$ and row $j$, as well as column $i$ and column $j$. After swapping, $A_{ii}$ becomes $A_{jj}$, which is nonzero. If $A_{ii} > 0$, $p$ is incremented. The algorithm then eliminates off-diagonal elements in column $i$ by subtracting a multiple of row $i$ from subsequent rows, and eliminates off-diagonal elements in row $i$ by subtracting a multiple of column $i$ from subsequent columns. Finally, it returns $p$. The pseudocode is described in \cref{alg:ldlt}.

  \begin{algorithm}[t]
    $p \gets 0$\;
    \For{$i \leftarrow 1$ \KwTo $n$}{
      $j \gets$ the smallest $j \ge i$ such that $A_{ij} \neq 0$\;
      \If{$j$ does not exist}{
        Continue to the next iteration\;
      }
      \If{$i < j$}{
        $A_{i*} \gets A_{i*} + A_{j*}$\tcc*[r]{Add row $j$ to row $i$}
        $A_{*i} \gets A_{*i} + A_{*j}$\tcc*[r]{Add col $j$ to col $i$}
        \If{$A_{ii} = 0$}{
          $A_{i*}, A_{j*} \gets A_{j*}, A_{i*}$\tcc*[r]{Swap row $i$ and row $j$}
          $A_{*i}, A_{*j} \gets A_{*j}, A_{*i}$\tcc*[r]{Swap col $i$ and col $j$}
        }
      }
      \If{$A_{ii} > 0$}{
        $p \gets p+1$\;
      }
      \For{$j \leftarrow i+1$ \KwTo $n$}{
        $A_{j*} \gets A_{j*} - (A_{ji}/A_{ii})A_{i*}$\tcc*[r]{Subtract a multiple of row $i$ from row $j$}
        $A_{*j} \gets A_{*j} - (A_{ij}/A_{ii})A_{*i}$\tcc*[r]{Subtract a multiple of col $i$ from col $j$}
      }
    }
    \KwRet{$p$}
    \caption{The $LDL^\intercal$ decomposition.} \label{alg:ldlt}
  \end{algorithm}

  Coming back to the description of the main program, we describe the input and output format.

  \paragraph{Input.} Each line represents a graph by two strings. The first string is the label of the graph, which ends with \texttt{+z} or \texttt{-z}. Here the \texttt{+} or \texttt{-} sign indicates whether we would like to compute positive or negative tail, and the number \texttt{z} stands for the number of red or reddish vertices. The second string is of the form \texttt{v[1]v[2]...v[2e-1]v[2e]}, which lists the edges.

  \paragraph{Output.} For each graph $F$, output one line containing \texttt{x y}, where \texttt{x} is the label of $F$, which ends with \texttt{+z} or \texttt{-z}, and \texttt{y} is the positive or negative tail of $F$ depending on the label.

  \medskip

  Our implementation is straightforward. We avoid floating-point errors by representing every number as a rational number. The actual code, written in Ruby, is available as the ancillary file \texttt{anc/tails.rb} in the arXiv version of this paper. Since the output has \texttt{z} greater than \texttt{y} on each line, we obtain a proof of \cref{prop:ptr,prop:ntr}. We provide the input in \cref{tab:in} for the convenience of anyone who wants to program independently, and we also provide the output in \cref{tab:out} for cross-checking.

  \buildtabledata{\inputtablecontents}{anc/input.txt}{2}
  \begin{table}[t]
    {\ttfamily
    \begin{tabular}{ll}
      \inputtablecontents
    \end{tabular}}
    \caption{Input.} \label{tab:in}
  \end{table}

  \buildtabledata{\outputtablecontents}{anc/output.txt}{10}
  \begin{table}
    {\ttfamily
    \begin{tabular}{llllllllll}
      \outputtablecontents
    \end{tabular}}
    \caption{Output.} \label{tab:out}
  \end{table}
\end{proof}

\section{Lean formalization} \label{sec:lean-details}

The formalization covers \cref{lem:pos-tail-16,lem:neg-tail-14} and all lemmas in \cref{sec:pos,sec:neg} needed for their proofs. Because our aim is to verify the case analysis rather than formalize all surrounding infrastructure, we make two simplifications. First, instead of formalizing the maximum cut $(A, B)$, we assume that the edges of $G[A] \cup G[B]$ form a matching, as established by \cref{lem:max-cut}; we call such a cut a \emph{matching cut}. We also assume that the cut admits none of the cut enhancers in \cref{fig:cut-enhancers,fig:degree-cut-enhancers}. Second, the formalized statements assert the existence of one of the tail reducers specified in \cref{fig:ptr,fig:degree-ptr,fig:ntr,fig:degree-ntr}, but omit its distance bound.

More precisely, the positive-tail and negative-tail case analyses are formalized as the following theorem.
\begin{theorem} \label{thm:lean-ptr} \label{thm:lean-ntr}
  For every matching cut $(A, B)$ of a subcubic graph $G$, if $G[A]$ contains at least one edge, then there exist a sequence of matching cuts
  \[
    (A_0,B_0), (A_1,B_1), \dots, (A_k,B_k)
  \]
  and a sequence of edges $a_0b_0, a_1b_1, \dots, a_{k-1}b_{k-1}$ of $G$, with $(A_0,B_0)=(A,B)$, such that, for each $i<k$, the vertex $a_i$ is red and $b_i$ is blue with respect to $(A_i,B_i)$, and
  \[
    A_{i+1}=A_i \oplus \sset{a_i,b_i}
    \quad \text{and} \quad
    B_{i+1}=B_i \oplus \sset{a_i,b_i}.
  \]
  Moreover, with all colors interpreted with respect to either $(A_k,B_k)$ or $(B_k,A_k)$, one of the graphs in \cref{fig:cut-enhancers,fig:degree-cut-enhancers,fig:degree-absolute-tail-reducer}, or in \cref{fig:ptr,fig:degree-ptr} (respectively, \cref{fig:ntr,fig:degree-ntr}), is an induced subgraph of $G$ satisfying the indicated degree constraints when present. \qed
\end{theorem}

The Lean source files are included in the ancillary directory \texttt{anc/lean} with the arXiv submission.

\end{document}